# Differential Calculus, Manifolds and Lie Groups over Arbitrary Infinite Fields

W. Bertram, H. Glöckner and K.-H. Neeb

**Abstract.** We present an axiomatic approach to finite- and infinite-dimensional differential calculus over arbitrary infinite fields (and, more generally, suitable rings). The corresponding basic theory of manifolds and Lie groups is developed. Special attention is paid to the case of mappings between topological vector spaces over non-discrete topological fields, in particular ultrametric fields or the fields of real and complex numbers. In the latter case, a theory of differentiable mappings between general, not necessarily locally convex spaces is obtained, which in the locally convex case is equivalent to Keller's $C_c^k$-theory.

**Contents:**

# Introduction

Differential calculus, as it is usually understood, may be seen as a formalism that allows to "pass to the limit" in certain difference quotients and to prove basic rules for this procedure – maybe the most important of these rules is the chain rule for $C^k$-maps, which makes it possible to define differential manifolds and Lie groups and many other important notions. Here, we present an approach to differential calculus on finite- and infinite-dimensional spaces which keeps close to this basic idea, but substantially widens the scope of differential calculus and the range of phenomena to which the familiar ideas of differential calculus can be applied. Thus, we can work over arbitrary infinite fields (and, more generally, suitable rings): none of the specific properties of the real or complex fields are needed for our considerations. In the real and complex cases, we also profit: it becomes possible to consider $C^k$-maps also between non-locally convex spaces. All these cases are treated in a uniform way, making it unnecessary to distinguish between "finite-" and "infinite-dimensional differential calculus" or "real" and "ultrametric differential calculus," and so on, for the basic theory. In a way, we have distilled those aspects and core results of differential calculus which are independent of any particular properties of the ground field.

This may sound utopian. In order to avoid misunderstandings, let us stress that we do not make such an assertion about *integral* calculus: in general, we will not be able to integrate even the simplest differential equation. The task of defining integral calculi leads to interesting and often difficult problems which highly depend on the given context. Differential calculus, however, appears as a part of analysis whose basic results are completely general, relying on elementary linear algebra and topology only, not on particular properties of the real or complex base field. This sharpens the eye for the peculiarities of the real, complex, and ultrametric cases, as it facilitates to distinguish clearly between the general, basic theory and those results, re-formulations and simplifications which indeed depend on specific properties of the ground field or of the topological vector spaces involved.

Let us explain our approach now, and describe the contents of the paper. Our basic idea may be formulated as follows: "You shall never separate the limit of a difference quotient from the difference quotient itself!" In other words, the limit of the difference quotient alone may be of little use, but it becomes useful if it is considered as the continuous extension of the difference quotient map. Note that a map $f\colon \mathbb{R}^n \supseteq U \to \mathbb{R}^m$ is of class $C^1$ (in the usual sense) if and only if the difference quotient map

$$(x, v, t) \mapsto \frac{f(x + tv) - f(x)}{t}$$

admits a continuous extension onto some neighbourhood in $U \times \mathbb{R}^n \times \mathbb{R}$ of the set determined by the condition $t = 0$. More precisely:

(∗) There exists a continuous map

$$f^{[1]}\colon U \times \mathbb{R}^n \times \mathbb{R} \supseteq U^{[1]} := \{(x,v,t)\mid x \in U,\, x + tv \in U\} \to \mathbb{R}^m$$

such that $f(x + tv) - f(x) = t \cdot f^{[1]}(x, v, t)$ for all $(x, v, t) \in U^{[1]}$.

In fact, if (∗) holds, then the differential of $f$ at $x$ is given by $df(x)v = f^{[1]}(x,v,0)$, and conversely, if $f$ is $C^1$, then the map $f^{[1]}$ defined by

$$f^{[1]}: U^{[1]} \to \mathbb{R}^m, \quad (x,v,t) \mapsto \begin{cases} \frac{f(x+tv)-f(x)}{t} & \text{if } t \in \mathbb{R}^\times \\ df(x)v & \text{if } t = 0 \end{cases} \qquad (1)$$

is continuous; this follows from the Fundamental Theorem of Calculus which provides, locally, an integral representation

$$f^{[1]}(x,v,t) = \int_0^1 df(x+stv)v \, ds. \qquad (2)$$

The same argument can be used to characterize $C^1$-maps in the Michal-Bastiani sense, also known as Keller's $C^1_c$-maps, between open subsets of locally convex spaces. Here, the main point is that the Fundamental Theorem of Calculus in its form (2) remains valid, regardless of completeness questions, because the value of the integral in (2), a priori taken in the completed space, belongs to the space itself since it is actually given by (1).

Motivated by these considerations, replacing $\mathbb{R}^n$ and $\mathbb{R}^m$ with topological vector spaces, we now take the characterization by (∗) as *definition* of the class $C^1$ of continuously differentiable mappings. This has several advantages: first of all, this is a reasonable definition even in the case of general (i.e., not necessarily locally convex) topological vector spaces, where the classical definitions stop to make sense. In a way, since no Fundamental Theorem of Calculus is available beyond locally convex spaces, we have incorporated its most relevant consequences in our definition of $C^1$-maps. Secondly, and much more importantly, the structure of the base field does not play a specific role in Condition (∗) because we only need to know what "continuous maps" are, that is, it works for general topological vector spaces over topological fields and even for topological modules over topological rings having a dense group of invertible elements. But we can (and will) interpret this way of defining differentiability in a still more general way: it may be considered as a machine that produces out of a class of mappings called "$C^0$" a new class called "$C^1$" (and admitting "differentials"). By no means is it necessary that $C^0$ should always mean "continuous"; any other class which shares certain basic formal properties with the class of continuous mappings will do (this idea is present already in [38]). These formal properties are used to define the notion of a "$C^0$-concept" (Section 1). Then (Section 2) we define the class $C^1$ by requiring the existence of a map $f^{[1]}$ of class $C^0$ such that (∗) holds. It turns out that a crucial property of $C^0$-maps $f: \mathbb{K} \supseteq U \to F$ is that their value at one point is already determined by their values on the complement of this point. This is true for continuous maps from non-discrete topological fields $\mathbb{K}$ into separated spaces, but it is is also true, for instance, if $\mathbb{K}$ is an infinite field and $C^0$ denotes the class of *rational* maps defined on Zariski-open parts of finite-dimensional $\mathbb{K}$-vector spaces, where "rational" is taken in the naive sense of being a quotient of a vector-valued polynomial and a scalar-valued polynomial. In this way we get a differential calculus for rational maps over arbitrary (infinite) fields which, in contrast to the well-known purely algebraic approach (cf. e.g. [34] or the appendix of [41]), is not merely "formal". A new feature in the "rational case" is that we have to use as underlying topology the Zariski topology; but the Zariski topology on a product space is not the product of the respective Zariski topologies (it is strictly finer). Thus in our general definition of $C^0$-concepts, we will allow product spaces to carry topologies that may be finer than the product topology.

The basic differentiation rules are proved in Section 3 and, for higher order derivatives, in Section 4. The proofs are "algebraic" in nature and in this way become often simpler and more transparent even than the usual proofs in $\mathbb{R}^n$ because we avoid the repeated use of the Mean Value Theorem (or of the Fundamental Theorem) which are no longer needed once they are incorporated in (∗). Of course, the class $C^2$ is defined to be the class of $C^1$-mappings such that $f^{[1]}$ (and not only $df$!) is again $C^1$; we then let $f^{[2]} := (f^{[1]})^{[1]}$ and so on, leading to the class $C^\infty$ and to iterated maps $f^{[k]}$ for all $k \in \mathbb{N}$. The maps $f^{[k]}$, in spite of their simple definition, are very complicated objects: they have $2^{k+1}-1$ arguments, and there exist non-trivial relations between their values on certain sub-diagonals. It is not exaggerated to say that a good deal of differential geometry can be understood as the study of the invariance properties of the maps $f^{[2]}$ and $f^{[3]}$.

In Section 5 we discuss Taylor's formula. It is truly remarkable that one can prove a "Taylor formula" in arbitrary characteristic (Theorem 5.1 and 5.4). That is, every $C^k$-map admits a finite expansion ("développement limité", cf. [17])

$$f(x+th) = \sum_{j=0}^{k} t^j a_j(x,h) + t^k R(x,h,t)$$

with a remainder term $R$ of class $C^0$ and taking the value 0 for $t = 0$ (Theorem 5.1). The interesting problem is then to identify the coefficients $a_j(x,h)$. They satisfy the relation $j! a_j(x,h) = d^j f(x)(h,\ldots,h)$ with the $j$-th differential $d^j f$ of $f$; if the characteristic of $\mathbb{K}$ is zero, this implies Taylor's formula in the usual form. In the general case, $a_j(x,\cdot)$ is a vector-valued form of degree $j$ (Theorem 5.4) which is given by a (non-canonical, in case of positive characteristic) polynomial expression (Theorem 5.6). The definition of forms and polynomial maps between modules, together with basic facts, are provided in an appendix (Appendix A). As far as we know, the purely algebraic problem of clarifying the precise relations between the latter concepts has never been fully investigated.

In Section 6 we look at the special case of *curves* $f \colon \mathbb{K} \supseteq U \to F$. In this case all arguments $x, v, t$ of the map $f^{[1]}(x,v,t)$ belong to the same space $\mathbb{K}$, and putting $v = 1$, the formalism may be simplified. We then get a characterization of $C^k$-maps which precisely generalizes the one of $C^k$-functions of one $p$-adic variable by W. Schikhof [58]. It follows that in the $p$-adic case Schikhof's and our definitions are equivalent (Proposition 6.9).

The main examples of $C^0$-concepts and associated differential calculi are (Section 7):

  (a) The case of Hausdorff topological vector spaces over $\mathbb{K} = \mathbb{R}$ or $\mathbb{C}$. It is proved that, if the range of $f$ is a locally convex topological vector space, our class $C^k$ coincides with the class of $C^k$-maps in the sense of Michal-Bastiani (Proposition 7.4). Furthermore, in the complex case a mapping with locally convex range is of class $C^\infty$ (over $\mathbb{C}$) if and only if it is complex analytic in the usual sense (Proposition 7.7).

(b) Generalizing (a), the case of Hausdorff topological vector spaces over a non-discrete topological field $\mathbb{K}$; in particular, the case of an ultrametric field such as the $p$-adic numbers. In the case of a function of one variable our definition of $C^k$-maps is equivalent to the one given by Schikhof in [58, Defn. 27.19]. As shown in [58], analytic functions between open subsets of a complete ultrametric field are of class $C^\infty$. We generalize this result to analytic mappings from open subsets of ultrametric normed vector spaces to locally convex spaces (Proposition 7.20). We also point out relations between $C^k$-maps and strictly differentiable mappings.

(c) The rational case – as explained above; we may take for $C^0$ the class of rational maps between finite-dimensional $\mathbb{K}$-vector spaces. "Intermediate cases" between the rational case (Zariski-topology) and the topological case (product topology on products) are also conceivable.

Once the differential calculus has been developed, there is no problem to define *smooth manifolds* and to establish the basic features of differential calculus on manifolds (Part II, Section 8). There are, however, some slight deviations from the usual theory: first of all, when working over base *rings*, one should avoid maximal atlasses and better work in the category of *manifolds with atlas*. Second, we can no longer use the algebra $C^\infty(M, \mathbb{K})$ of smooth functions on a manifold $M$ to define geometric objects since it may, even locally, be reduced to the constants (see Remark 8.1). Thus vector fields should be defined "geometrically" and not as differential operators. A systematic discussion of differential geometry will be given in a subsequent paper [8]. In particular, it will be proved there that the tangent bundle $TM$ of $M$ is, in a natural way, a manifold over the ring $\mathbb{K}[x]/(x^2)$ (the "dual numbers" over $\mathbb{K}$). Thus, a "synthetic differential geometry for manifolds" in the sense of [47] is possible in our context.

*Lie groups* are defined as usual to be groups in the category of manifolds in question (Section 9). The Lie algebra of a Lie group is defined and it is proved that in this way we get a functor into $C^0$-Lie algebras over $\mathbb{K}$ (Theorem 9.1). We also define *symmetric spaces* (following [40]); the general theory (which includes an analogue of Theorem 9.1) will be developed in [8] and [9]. Differential geometric aspects of Lie groups will be discussed in [8], and the main classes of infinite-dimensional Lie groups over topological fields are constructed in [27]; see Section 13 for an overview.

The constructions in [27] rely on results which are specific to the case of differential calculus over topological fields which are established in Part III of this paper. In particular, we transfer ideas of the convenient differential calculus of Frölicher, Kriegl and Michor (see [21], [36]) to the non-locally convex or ultrametric case. Notably, it is shown that a map $f: U \to F$ on an open subset of a metrizable topological vector space $E$ over $\mathbb{K} = \mathbb{R}$ or an ultrametric field $\mathbb{K}$, with values in any topological $\mathbb{K}$-vector space $F$, is of class $C^k$ if and only if $f \circ c \colon \mathbb{K}^{k+1} \to F$ if of class $C^k$ for each smooth map $c \colon \mathbb{K}^{k+1} \to U$ (Theorem 12.4). Note that neither $E$ nor $F$ needs to be locally convex here. For the real locally convex special case, compare [35]; and apparently, this result has also been inspired by Souriau's theory of diffeological spaces ([60], [42]). The utility of the preceding result, and its high potential for applications, is illustrated in [26], where it is used to deduce a generalized implicit function theorem for mappings on metrizable spaces over ultrametric fields, which in turn facilitates to construct a Lie group structure on diffeomorphism groups of $\sigma$-compact, finite-dimensional smooth manifolds over local fields (of arbitrary characteristic) [27]. Cf. [43] and further works by S. V. Ludkovsky for diffeomorphism groups in characteristic 0. In the final section, we describe (without proofs) examples of Lie groups over topological fields. All of the major constructions of (infinite-dimensional) Lie groups familiar from the real or complex locally convex case turn out to be extremely robust, and can be used just as well to produce Lie groups over arbitrary topological fields, valued fields, or at least arbitrary local fields.

**Further topics and open problems.** The present work is basic for several subsequent developments and raises many questions. We already mentioned the problem of defining suitable *integral calculi*. For the $p$-adic case the reader may consult [58]. Another integration problem is the *integration of Lie algebras*: which $C^0$-Lie algebras over $\mathbb{K}$ can be integrated to Lie groups? This question generalizes the enlargibilty problem of Banach-Lie algebras solved by van Est and Korthagen [20]; see also [50] for results in the real locally convex case. It is remarkable that in the category of *Jordan algebraic structures* the corresponding integration problem can be solved under very general assumptions [9]. Thus, via the Jordan-Lie functor (cf. [6], [7]) one gets a rich family of examples of Lie structures (Lie groups and symmetric spaces) over $\mathbb{K}$, and it would be very interesting to know "how big" this family is inside the Lie category.

Another interesting problem would be to define differential calculi over *non-commutative* base fields or -rings. In Sections 1 and 2, $\mathbb{K}$ could be non-commutative, but commutativity is crucial for product- and quotient rule in Section 3. Thus, if $\mathbb{K}$ is non-commutative, in general not even the polynomials on $\mathbb{K}$ would be differentiable. However, it might be possible to compensate non-commutativity by an additional "twistor"-structure and to define a modified differential calculus; this would be rather special but may be interesting in connection with exceptional geometries.

Another task (which we leave to the reader) would be to define "pointwise" $C^0$-concepts and a notion of maps that are differentiable "at a point" and to prove the analogues of the usual differentiation rules in this set-up (cf. [2] for an extensive comparative study of remainder conditions and notions of differentiability at a point).

As already mentioned, in our framework we are naturally lead to a more algebraic approach to *differential geometry* where e.g. dual numbers can be used [8]. Moreover, if one wants to generalize the theory of manifolds over $\mathbb{K}$ to include also "singular manifolds" and to include all mapping spaces ("Cartesian closedness"), it seems possible to develop a theory of "smooth toposes over $\mathbb{K}$" generalizing much of the real theory from [47].

**Acknowledgements.** The first named author owes much of the inspiration that lead to the present approach to lectures by his former teacher H. S. Holdgrün (the reader may notice this by comparing with the book [32]), and to suggestions by his collegue Yanis Varouchas who died before the final version of this text was ready. We would like to dedicate this work to his memory.

**Notation.** In the following, $\mathbb{K}$ denotes a commutative ring with unit 1 (for a first reading, the reader is well-advised to think of $\mathbb{K}$ as a field; for this reason we use the notation $\mathbb{K}$ instead of the possibly more appropriate $R$).

# Part I: Differential calculus

# 1 $C^0$-concepts

A *topologized* $\mathbb{K}$-*module* is a $\mathbb{K}$-module $E$, equipped with a topology $\mathcal{T}(E)$ (which need not have any particular properties). We shall frequently omit the word "topologized," for convenience. The base ring $\mathbb{K}$ is also assumed to be equipped with a topology.

**Definition 1.1** A $C^0$-*concept over* $\mathbb{K}$ consists of the following data:

(a) a class $\mathcal{M}$ of topologized $\mathbb{K}$-modules (usually denoted by $E, F, \ldots$) is given, such that $\mathbb{K} \in \mathcal{M}$;

(b) for all $E, F \in \mathcal{M}$ and each open set $U \in \mathcal{T}(E)$, a subset $C^0(U, F)$ of the set $C(U, F)$ of continuous functions from $U$ to $F$ is given;

(c) a functional class is given which associates to each pair $(E_1, E_2)$ of topologized $\mathbb{K}$-modules $E_1, E_2 \in \mathcal{M}$ a topology $\mathcal{T}(E_1 \times E_2)$ on the direct product $E_1 \times E_2$ (which need not be the product topology), such that $(E_1 \times E_2, \mathcal{T}(E_1 \times E_2)) \in \mathcal{M}$.

These data are required to satisfy the following three groups of axioms:

**I. Basic axioms not involving products.**

(I.1) If two $C^0$-maps are composable, then their composition is $C^0$. The identity map $E \to E$ is $C^0$, and, more generally, all inclusion maps $U \hookrightarrow E$, $U \in \mathcal{T}(E)$, are $C^0$. This implies that restrictions of $C^0$-maps to open subsets are $C^0$.

(I.2) All translations and dilations, and hence all affine maps of the form

$$E \to E, \quad x \mapsto rx + b, \quad r \in \mathbb{K}, \, b \in E,$$

are $C^0$. This includes the assumption that the topology is invariant under translations and proper dilations.

(I.3) For any $v, x \in E$, the affine map $\mathbb{K} \to E$, $t \mapsto tv + x$ is $C^0$.

(I.4) The set $\mathbb{K}^\times$ is open in $\mathbb{K}$, and the inversion map $i: \mathbb{K}^\times \to \mathbb{K}$ is $C^0$.

(I.5) (Locality.) If $f: U \to F$ is a mapping such that $f|_{U_i} \in C^0(U_i, F)$ for some open cover $(U_i)_{i \in I}$ of $U$, then $f \in C^0(U, F)$.

**II. Basic axioms involving products.**

(II.1) The projections $\mathrm{pr}_i : E_1 \times E_2 \to E_i$ are $C^0$, and for all $x \in E_1, y \in E_2$, the mappings
$$E_1 \to E_1 \times E_2, \quad v \mapsto (v, y), \qquad E_2 \to E_1 \times E_2, \quad w \mapsto (x, w)$$
are $C^0$.

(II.2) If $f_i: U_i \to F_i$ are $C^0$, $i = 1, 2$, then $f_1 \times f_2 : U_1 \times U_2 \to F_1 \times F_2$ is $C^0$.

(II.3) The diagonal map
$$\delta: E \to E \times E, \quad \delta(x) := (x, x)$$
is $C^0$.

(II.4) The exchange maps $E \times F \to F \times E$ and the natural isomorphisms $(E_1 \times E_2) \times E_3 \cong E_1 \times (E_2 \times E_3)$ are $C^0$ in both directions.

(II.5) The structural maps of the $\mathbb{K}$-module $E$,
$$E \times E \to E, \quad (x, y) \mapsto x + y,$$
$$\mathbb{K} \times E \to E, \quad (r, v) \mapsto rv$$
are $C^0$.

**III. Determination Axiom.** A $C^0$-map $f: U \to F$ on an open subset $U \subseteq \mathbb{K}$ is uniquely determined by its values on $U \cap \mathbb{K}^\times$. (In view of the preceding requirements this is equivalent to: If $f|_{U \cap \mathbb{K}^\times} = 0$, then $f = 0$. If $\mathbb{K}$ is a field, then these properties are equivalent to: A $C^0$-map $f: \mathbb{K} \supseteq U \to F$ is uniquely determined by its values on the complement of one point).

Of course, the topology on a topologized module $E \in \mathcal{M}$ need not be determined by the underlying abstract module: various topologies on a given $\mathbb{K}$-module may turn it into an element of $\mathcal{M}$.

**Remark 1.2** Some immediate consequences of the basic axioms I and II are:

(a) partial maps obtained by fixing one or several arguments of a $C^0$-map $f: E_1 \times \cdots \times E_k \supseteq U \to F$ are again $C^0$;

(b) a mapping $f = (f_1, f_2) : E \supseteq U \to F_1 \times F_2$ is $C^0$ if and only if so are its coordinate

functions $f_i := \operatorname{pr}_i \circ f \colon U \to F_i$ for $i \in \{1,2\}$ (using that $f = (f_1 \times f_2) \circ \delta|_U$, where $\delta \colon E \to E \times E$ is the diagonal map);

(c) sums and multiples of $C^0$-maps having values in the same space are again $C^0$; in particular,
$$L(E,F) := \operatorname{Hom}_{\mathbb{K}}(E,F) \cap C^0(E,F)$$
is a $\mathbb{K}$-module (which, however, we shall not try to turn into a member of $\mathcal{M}$). Moreover, if $E \cong \mathbb{K}^n$, then every linear map $f \colon E \to F$ is of the form $f = \sum_{i,j=1}^n f_{ij}\, \iota_i \circ p_j$ with projections $p_j \colon \mathbb{K}^n \to \mathbb{K}$ and linear maps $\iota_i \colon \mathbb{K} \to E$ and hence is $C^0$; thus $L(\mathbb{K}^n, F) = \operatorname{Hom}(\mathbb{K}^n, F)$. In particular, $\operatorname{GL}(n, \mathbb{K}) \subseteq L(\mathbb{K}^n, \mathbb{K}^n)$, and it follows that the notion of a $C^0$-map on $\mathbb{K}^n$ does not depend on the basis chosen.

**Remark 1.3** Not all of the axioms of a $C^0$-concept are independent of each other. For instance, (I.2) and (I.3) are implied by the axioms of group II. However, we prefer to mention these axioms explicitly since they do not involve products.

The Locality Axiom I.5 will hardly be used in our development of the general theory; it could be omitted if in Sections 8 and 9, the reader accepted to work entirely in the category of manifolds with atlas (instead of using maximal atlasses).

**Remark 1.4** The Determination Axiom III is the key axiom for our differential calculus. Some easy consequences are:

(a) The topology on $\mathbb{K}$ and hence on all other spaces cannot be discrete. In fact, if $\{x_0\}$ were open in $\mathbb{K}$, then so would be $U := \{0\}$, and now $f \colon U \to \mathbb{K}$, $x \mapsto 0$ and $g \colon U \to \mathbb{K}$, $x \mapsto 1$ would be distinct $C^0$-functions on $U$ whose restriction to $U \cap \mathbb{K}^\times = \emptyset$ coincides, contradicting Axiom III. The same argument shows that $\mathbb{K} \setminus \mathbb{K}^\times$ contains no non-empty open subset, i.e., $\mathbb{K}^\times$ is dense in $\mathbb{K}$.

(b) Assume $U \subseteq E$ is a neighbourhood of $0$. Then $U$ is *absorbing* in the sense that $\mathbb{K} U = E$. (In fact, for $x \in E$ consider the set $\{t \in \mathbb{K} \colon tx \in U\}$ which, in virtue of I, is a neighbourhood of $0$ in $\mathbb{K}$. As we have just remarked, it must then contain an element $t \in \mathbb{K}^\times$. Therefore $x = t^{-1}(tx) \in \mathbb{K} U$). In particular, proper submodules of $E$ are never open in $E$, nor do they contain open subsets.

(c) Assume $f \colon E \supseteq U \to F$ is $C^0$ and $H \subset E$ is a proper submodule. Then $f$ is uniquely determined by its values on $U \setminus H$. This follows by applying III to $g(t) := f((1-t)x + ty)$, where $x$ is an arbitrary point of $U \cap H$ and $y$ an arbitrary point of $U \setminus H$.

The following definitions allow us to use in our general context some notions that are familiar from ordinary analysis:

**Definition 1.5**

(i) A subset $U_1 \subseteq U$ is called $C^0$-*dense* if, for all $F \in \mathcal{M}$ and all $f \in C^0(U, F)$, $f|_{U_1} = 0$ implies $f = 0$.

(ii) If $x_0 \in U$ and $f \colon U \to F$ is $C^0$, we say that the function $f|_{U \setminus \{x_0\}}$ *admits a limit at* $x_0$ which, by definition, is $f(x_0)$, and we use the notation

$$\lim_{x \to x_0} f(x) := f(x_0).$$

The $C^0$-map $f$ being uniquely determined by its restriction $f|_{U \setminus \{x_0\}}$ by the Determination Axiom, we see that $\lim_{x \to x_0} f(x)$ is uniquely determined by $f|_{U \setminus \{x_0\}}$ (and thus well-defined). More generally, if $U_1 \subseteq U$ is $C^0$-dense and $x_0 \in U \setminus U_1$, we use the same notation. Then $\lim_{x \to x_0} f(x)$ is uniquely determined by $f|_{U_1}$.

Note that, even when we are dealing with non-Hausdorff spaces, Definition 1.5 (ii) makes it possible to define a unique limit in the situations we are interested in (one of the various limits in the sense of general topology is picked out). However, we do not want to discuss the notion of limits of "arbitrary" maps, nor do we want to define for such maps the notion of being "$C^0$ at one point".

**Remark 1.6** (Fields versus rings.) All basic results of differential calculus from Sections 2, 3 and 4 hold for base rings as well as for base fields. The only exceptions are Lemma 2.6, Lemma 3.9 (linking the continuous extension of difference quotients and partial difference quotients), and Lemma 4.9 saying that being $C^k$ is a local property.

Therefore, when defining manifolds and Lie groups over base rings, $C^k$-properties must in principle be checked for all charts of a given atlas $\mathcal{A}$ and not only for some suitable sub-atlas. This is no problem as long as one works in the category of manifolds with atlas.

## 2 From $C^0$ to $C^1$

We fix for the rest of the paper a $C^0$-concept over $\mathbb{K}$ and assume henceforth without further mention that all $\mathbb{K}$-modules $E, F, \ldots$ belong to the class $\mathcal{M}$. Generically, a non-empty open subset of $E$ will be denoted by $U$.

**Definition 2.1** Let $f \in C^0(U, F)$. We write $f \in C^1(U, F)$ or say that $f$ *is of class $C^1$* if there exists a $C^0$-map

$$f^{[1]} \colon U \times E \times \mathbb{K} \supseteq U^{[1]} := \{(x, v, t) \mid x \in U, x + tv \in U\} \to F$$

such that

$$f(x + tv) - f(x) = t \cdot f^{[1]}(x, v, t)$$

whenever $(x, v, t) \in U^{[1]}$. (Note that $U^{[1]}$ is open since $(x, v, t) \mapsto x + tv$ is $C^0$ by the axioms of group II.) Put another way (when $\mathbb{K}$ is a field): the $C^0$-map

$$U^{[1]} \cap (E \times E \times \mathbb{K}^\times) \to F, \quad (x, v, t) \mapsto \frac{f(x + tv) - f(x)}{t}$$

admits a $C^0$-continuation "for $t = 0$", i.e. onto $U^{[1]}$. Taking $t = 0$, we define the *differential of $f$ at $x$* by

$$df(x)\colon E \to F, \quad v \mapsto df(x)v := f^{[1]}(x, v, 0).$$

(By the Determination Axiom III, the $C^0$-map $f^{[1]}$ is uniquely determined by $f$, and hence $df(x)$ is well-defined.) We shall also use the notation $df(x, v) := df(x)(v)$.

**Proposition 2.2** *If $f \in C^1(U, F)$ and $x \in U$, then $df(x)\colon E \to F$ is a $\mathbb{K}$-linear $C^0$-map.*

**Proof.** Since $f^{[1]}$ is $C^0$, it follows from Axiom II that the partial map $df(x) = f^{[1]}(x, \cdot, 0)$ also is $C^0$. Let us show that it is $\mathbb{K}$-linear: *Homogeneity:*

$$sf^{[1]}(x, rv, s) = f(x + s(rv)) - f(x) = f(x + srv) - f(x) = srf^{[1]}(x, v, sr);$$

both sides are $C^0$ in $s$; for $s \in \mathbb{K}^\times$ we divide by $s$ and get

$$f^{[1]}(x, rv, s) = rf^{[1]}(x, v, sr).$$

By uniqueness of the $C^0$-continuation, this still holds for $s = 0$ and then gives $df(x)(rv) = rdf(x)v$.

*Additivity:* For $s \in \mathbb{K}$ sufficiently close to 0, we have

$$\begin{aligned} sf^{[1]}(x, v + w, s) &= f(x + s(v + w)) - f(x) \\ &= f(x + sv + sw) - f(x + sv) + f(x + sv) - f(x) \\ &= sf^{[1]}(x + sv, w, s) + sf^{[1]}(x, v, s). \end{aligned}$$

By the same argument as above it follows that

$$f^{[1]}(x, v + w, s) = f^{[1]}(x + sv, w, s) + f^{[1]}(x, v, s),$$

which for $s = 0$ gives $df(x)(v + w) = df(x)w + df(x)v$. $\square$

**Example 2.3** Consider an affine map $f(x) = Ax + b$, where $A$ is linear and of class $C^0$. In this case $sf^{[1]}(x, v, s) = A(x + sv) + b - (Ax + b) = sAv$, whence

$$f^{[1]}(x, v, s) = Av \quad \text{and} \quad df(x)v = Av. \tag{3}$$

**Notation 2.4** The $C^0$-map $f^{[1]}$ whose existence is required by Definition 2.1 will be called the *difference quotient map of $f$*. Depending on the context, the following notation involving the difference symbol $\Delta$ will also be useful:

$$\frac{\Delta f}{\Delta}(x, v, t) := \Delta_{v,t} f(x) := f^{[1]}(x, v, t).$$

For $t = 1$, $\Delta_{v,t}$ is the *(ordinary) difference operator*

$$\Delta_v f(x) := \Delta_{v,1} f(x) = f(x + v) - f(x) = f^{[1]}(x, v, 1)$$

(cf. [17] and Appendix A). The *directional derivative* in direction $v$

$$\partial_v f(x) := \Delta_{v,0} f(x) = f^{[1]}(x,v,0) = df(x)v$$

is defined for all $x \in U$ and $v \in E$. Finally, we define the *tangent map* by

$$Tf : U \times E \to F \times F, \quad (x,v) \mapsto (f(x), df(x)v)$$

and the *extended tangent map* by

$$\hat{T}f : U^{[1]} \to F \times F \times \mathbb{K}, \quad (x,v,t) \mapsto (f(x), f^{[1]}(x,v,t), t).$$

For $t = 0$, $\hat{T}f$ contains $Tf$ as a partial map, and for $t = 1$, it contains the ordinary difference map as a partial map. Finally, note that for differentiable *curves*, i.e. for $C^1$-maps of open sets of $\mathbb{K}$ into $F$, Definition 2.1 can be simplified since it suffices here to take $v = 1$ – see Section 6 where also a simplified notation for this case is introduced.

**Remark 2.5** Our definition of differentiability is related to more common ones (as e.g. in [37]) as follows: let us say that a $C^0$-map $h \mapsto f(h)$, defined on an open neighbourhood of the origin in $E$, *is $O(h)$* if $f(0) = 0$. Now assume that $f$ is $C^1$ and write

$$\begin{aligned} f(x+th) &= f(x) + tf^{[1]}(x,h,t) \\ &= f(x) + tdf(x)h + t(f^{[1]}(x,h,t) - f^{[1]}(x,h,0)) \end{aligned}$$

where $t \mapsto f^{[1]}(x,h,t) - f^{[1]}(x,h,0)$ clearly is $O(t)$. Thus we may say that

$$f(x+th) = f(x) + tdf(x)h + tO(t),$$

which is in close analogy with the usual definition of differentiability. Iterating this argument, we are lead in a straightforward way to Taylor's formula (see Section 5).

It is useful to observe that in Definition 2.1, the set $U^{[1]}$ may be replaced with any smaller open neighbourhood of $U \times E \times \{0\}$, provided $\mathbb{K}$ is a field:

**Lemma 2.6** *Let $\mathbb{K}$ be a field, $f \in C^0(U,F)$, and suppose that there exists a $C^0$-map $g : P \to F$, defined on an open neighbourhood $P$ of $U \times E \times \{0\}$ in $U^{[1]}$, such that*

$$g(x,y,t) = \frac{f(x+ty) - f(x)}{t} \tag{4}$$

*for all $(x,y,t) \in P$, $t \neq 0$. Then $f$ is of class $C^1$, and $g = f^{[1]}|_P$.*

**Proof.** Using Eqn. (4) also for $(x,y,t) \in U^{[1]} \setminus P \subseteq U^{[1]} \cap (U \times E \times \mathbb{K}^\times)$, we extend $g$ to a mapping $U^{[1]} \to F$, also denoted $g$. Then $g$ is $C^0$ on the open set $P$ (by hypothesis), and also $C^0$ on the open set $Q := \{(x,y,t) \in U^{[1]} : t \in \mathbb{K}^\times\}$, as it is given by (4) there (and $f$ is $C^0$). Since $U^{[1]} = P \cup Q$, we deduce with the Locality Axiom I.5 that $g : U^{[1]} \to F$ is $C^0$, and thus $f$ is $C^1$, with $f^{[1]} = g$. □

# 3 Chain rule and other differentiation rules

In this section, we establish the familiar basic rules of calculus in our general setting.

First, let us show that compositions of composable $C^1$-maps are $C^1$. To this end, assume that $f\colon U \to V \subseteq F$, $g\colon V \to H$ are defined on open subsets $U \subseteq E$ and $V \subseteq F$, where $E, F, H \in \mathcal{M}$. Then we have:

**Proposition 3.1 (Chain Rule)** *If $f$ and $g$ are of class $C^1$, then also the composition $g \circ f\colon U \to H$ is of class $C^1$. We have $(f(x), f^{[1]}(x,y,t), t) \in V^{[1]}$ for all $(x,y,t) \in U^{[1]}$, and*
$$(g \circ f)^{[1]}(x,y,t) = g^{[1]}(f(x), f^{[1]}(x,y,t), t). \tag{5}$$
*In particular, $d(g \circ f)(x) = dg(f(x)) \circ df(x)$ for all $x \in U$.*

**Proof.** Let $(x,y,t) \in U^{[1]}$. Then $f(x+ty) = f(x) + t f^{[1]}(x,y,t)$, hence $(f(x), f^{[1]}(x,y,t), t) \in V^{[1]}$ and
$$g(f(x+ty)) - g(f(x)) = g(f(x) + t f^{[1]}(x,y,t)) - g(f(x)) = t \cdot g^{[1]}(f(x), f^{[1]}(x,y,t), t),$$
where $U^{[1]} \to H$, $(x,y,t) \mapsto g^{[1]}(f(x), f^{[1]}(x,y,t), t)$ is $C^0$ by Axioms I and II. Thus $g \circ f$ is $C^1$, with $(g \circ f)^{[1]}$ given by (5). Taking $t = 0$ in (5), the final formula follows. $\square$

Using Notation 2.4, Equation (5) is tantamount to
$$\hat{T}(g \circ f) = \hat{T}g \circ \hat{T}f, \tag{6}$$
and for $t = 0$ we get the chain rule in the form
$$T(g \circ f) = Tg \circ Tf. \tag{7}$$

**3.2 Linearity.** If $f$ and $g$ are $C^1$, with values in the same space, then $f + g$ and $\lambda f$ ($\lambda \in \mathbb{K}$) are $C^1$, and
$$(f + g)^{[1]} = f^{[1]} + g^{[1]}, \quad (\lambda f)^{[1]} = \lambda f^{[1]}.$$
This follows from a trivial calculation, together with Axiom II which ensures that sum and multiples of $C^0$-maps are $C^0$.

**3.3 Product rule.** If $b\colon E_1 \times E_2 \to F$ is bilinear and $C^0$, then it is $C^1$ with
$$db(x_1, x_2)(h_1, h_2) = b(x_1, h_2) + b(h_1, x_2).$$
This follows from
$$b(x_1 + th_1, x_2 + th_2) - b(x_1, x_2) = t(b(h_1, x_2) + b(x_1, h_2) + t b(h_1, h_2)).$$

Combined with the chain rule, we can now differentiate "products of functions" $b \circ (g_1 \times g_2)$ in the usual way, and we see that $C^0$-polynomials are $C^1$. In particular, if $f, g : U \to \mathbb{K}$ are scalar-valued functions, then

$$(fg)^{[1]}(x, v, t) = f^{[1]}(x, v, t)g(x) + f(x)g^{[1]}(x, v, t) + tf^{[1]}(x, v, t)g^{[1]}(x, v, t)$$

which for $t = 0$ gives
$$d(fg)(x)v = f(x)dg(x)v + g(x)df(x)v.$$

Using Notation 2.4, we see that $\partial_v$ is a derivation in the sense that

$$\partial_v(fg) = (\partial_v f)g + f(\partial_v g).$$

**3.4 Example: Polynomial maps.** For the squaring map $q : \mathbb{K} \to \mathbb{K}$, $x \mapsto x^2$, we get

$$q^{[1]}(x, h, t) = \frac{(x + th)^2 - x^2}{t} = 2xh + th^2,$$

whence $dq(x) = 2x$ (note, however, that $q^{[1]}$ is a polynomial of degree 3). Similarly, for $g(x) = x^n$ we have $dg(x) = nx^{n-1}$ (note that this is zero if $n = 0$ in $\mathbb{K}$). Moreover, we see that, for any polynomial map $p : \mathbb{K} \to \mathbb{K}$, $p^{[1]} : \mathbb{K}^3 \to \mathbb{K}$ is again a polynomial map (which in general is of *higher* degree than $p$). For the definition of polynomial maps between general $\mathbb{K}$-modules, we refer to Definition A.5 in the appendix. If $m \colon E^k \to F$ is a $\mathbb{K}$-multilinear map of class $C^0$, then it is of class $C^1$, and the polynomial map $f \colon E \to F$ obtained by restricting $m$ to the diagonal is again $C^1$ (recall that the diagonal map $\delta$ is $C^0$ by Axiom II.3, and hence $C^1$). In particular, polynomial maps $\mathbb{K}^n \to F$ are always $C^1$ (cf. Remark 1.2 (c)). If the domain $E$ is not a free module, then no general statements of this kind can be made.

**3.5 Quotient rule.** Recall from Axiom I and II that the multiplication map $\mathbb{K} \times \mathbb{K} \to \mathbb{K}$ and the inversion map $i \colon \mathbb{K}^\times \to \mathbb{K}^\times$ are $C^0$. Then the inversion map is in fact $C^1$:

$$\begin{aligned} i(x + tv) - i(x) &= (x + tv)^{-1} - x^{-1} = x^{-1}(x - (x + tv))(x + tv)^{-1} \\ &= -x^{-1}tv(x + tv)^{-1} = -t\, x^{-1}v(x + tv)^{-1} \end{aligned}$$

(note that we have used from the commutativity of $\mathbb{K}$ only that $t$ belongs to the center of $\mathbb{K}$), and thus
$$i^{[1]}(x, v, t) = -x^{-1}v(x + tv)^{-1},$$

and hence $di(x)v = -x^{-1}vx^{-1}$. (These arguments go through more generally for any unital $\mathbb{K}$-algebra with open unit group and such that product and inversion are $C^0$: if this is the case, then product and inversion are actually $C^1$.)

**3.6 Example: Rational maps.** We say that $f : U \to F$ is a *rational map* if $f(x) = q(x)^{-1}p(x)$ with polynomial maps $p : E \to F$, $q : E \to \mathbb{K}$ (in the sense defined in Appendix A) such that $q(U) \subseteq \mathbb{K}^\times$. Product, chain and quotient rule imply that rational maps are $C^1$ if nominator and denominator are $C^1$ and are differentiated in the usual way. In particular, rational maps $f \colon \mathbb{K}^n \supseteq U \to F$ are always $C^1$.

**3.7 (Direct Products).** If $f_i \colon U_i \to F_i$, $i = 1, 2$, are $C^1$, then so is $f_1 \times f_2 \colon U_1 \times U_1 \to F_1 \times F_2$, with

$$(f_1 \times f_2)^{[1]}((x_1, x_2), (v_1, v_2), s) = \left((f_1)^{[1]}(x_1, v_1, s), (f_2)^{[1]}(x_2, v_2, s)\right). \tag{8}$$

This readily follows from the definitions and Axiom II.

**3.8 Rule on partial derivatives.** Assume $f \colon U \to F$ is a $C^1$-map, where $U$ is open in $E_1 \times E_2$. Then it follows from Axiom II that the partial maps $f(x_1, \cdot)$ and $f(\cdot, x_2)$ are also $C^1$. Then

$$\begin{aligned}
f(x_1 + tv_1, x_2 + tv_2) - f(x_1, x_2) &= f(x_1 + tv_1, x_2 + tv_2) - f(x_1 + tv_1, x_2) \\
&\quad + f(x_1 + tv_1, x_2) - f(x_1, x_2) \\
&= t\frac{\Delta_1 f}{\Delta_1}(x_1, x_2, v_1, t) + t\frac{\Delta_2 f}{\Delta_2}(x_1 + tv_1, x_2, v_2, t),
\end{aligned}$$

where we denote by $\frac{\Delta_i f}{\Delta_i}$ for $i \in \{1, 2\}$ the partial difference quotient maps of $f$. Dividing by $t \in \mathbb{K}^\times$ and letting $t \to 0$, we get the usual rule on partial derivatives:

$$df(x_1, x_2)(v_1, v_2) = d_1 f(x_1, x_2) v_1 + d_2 f(x_1, x_2) v_2,$$

where $d_i f(x_1, x_2, v_i) := \frac{\Delta_i f}{\Delta_i}(x_1, x_2, v_i, 0)$ for $i \in \{1, 2\}$.

For the converse we need the Locality Axiom and stronger assumptions on $\mathbb{K}$:

**Lemma 3.9** *Assume $\mathbb{K}$ is a field and $f \colon U \to F$ is a $C^0$-map, where $U$ is open in $E_1 \times E_2$. Let $U_1 := \{(x_1, x_2, v_1, t) \in U \times E_1 \times \mathbb{K} \colon (x_1 + tv_1, x_2) \in U\}$ and $U_2 := \{(x_1, x_2, v_2, t) \in U \times E_2 \times \mathbb{K} \colon (x_1, x_2 + tv_2) \in U\}$. Suppose that there are $C^0$-maps $\frac{\Delta_i f}{\Delta_i} \colon U_i \to F$ for $i \in \{1, 2\}$ extending the partial difference quotient maps of $f$, viz.*

$$\frac{f(x_1 + tv_1, x_2) - f(x_1, x_2)}{t} = \frac{\Delta_1 f}{\Delta_1}(x_1, x_2, v_1, t) \quad \text{for all } (x_1, x_2, v_1, t) \in U_1 \text{ such that } t \neq 0,$$

*and similarly for $\frac{\Delta_2 f}{\Delta_2}$. Then $f$ is a mapping of class $C^1$, and*

$$\frac{\Delta f}{\Delta}((x_1, x_2), (v_1, v_2), t) = \frac{\Delta_1 f}{\Delta_1}(x_1, x_2, v_1, t) + \frac{\Delta_2 f}{\Delta_2}(x_1 + tv_1, x_2, v_2, t) \tag{9}$$

*for all $(x_1, x_2, v_1, v_2, t) \in P := \{(x_1, x_2, v_1, v_2, t) \in U^{[1]} \colon x_1 + tv_1 \in U\}$.*

**Proof.** Define $g(x_1, x_2, v_1, v_2, t)$ by the right hand side of (9) for $(x_1, x_2, v_1, v_2, t) \in P$. Then $P$ is an open neighbourhood of $U \times E_1 \times E_2 \times \{0\}$ in $U^{[1]}$, and for all $(x_1, x_2, v_1, v_2, t) \in P$ such that $t \neq 0$, we have

$$\begin{aligned}
f(x_1 + tv_1, x_2 + tv_2) - f(x_1, x_2) &= f(x_1 + tv_1, x_2 + tv_2) - f(x_1 + tv_1, x_2) \\
&\quad + f(x_1 + tv_1, x_2) - f(x_1, x_2) \\
&= tg(x_1, x_2, v_1, v_2, t)
\end{aligned}$$

(cf. **3.8**). Lemma 2.6 shows that $f$ is $C^1$, with $f^{[1]}|_P = g$. □

# 4 Higher differentials: from $C^1$ to $C^2$ and $C^\infty$

**Definition 4.1** Let $f\colon E \supseteq U \to F$ be of class $C^1$. We say that $f$ is $C^2(U, F)$ or *of class $C^2$* if $f^{[1]}$ is $C^1$, in which case we define $f^{[2]} := (f^{[1]})^{[1]}\colon U^{[2]} \to F$, where $U^{[2]} := (U^{[1]})^{[1]}$. Inductively, we say that $f$ is $C^{k+1}(U, F)$ or *of class $C^{k+1}$* if $f$ is of class $C^k$ and $f^{[k]}\colon U^{[k]} \to F$ is of class $C^1$, in which case we define $f^{[k+1]} := (f^{[k]})^{[1]}\colon U^{[k+1]} \to F$ with $U^{[k+1]} := (U^{[k]})^{[1]}$. We also use the notation $\frac{\Delta^k f}{\Delta^k} := f^{[k]}$. The map $f$ is called *smooth* or *of class $C^\infty$* if it is of class $C^k$ for each $k \in \mathbb{N}_0$.

**Remark 4.2** Note that $U^{[k+1]} = (U^{[1]})^{[k]}$ for each $k \in \mathbb{N}_0$, and that $f$ is of class $C^{k+1}$ if and only if $f$ is of class $C^1$ and $f^{[1]}$ is of class $C^k$; in this case, $f^{[k+1]} = (f^{[1]})^{[k]}$. These claims are proved by a trivial induction.

**Example 4.3** Every linear $C^0$-map $\lambda$ is smooth, by a trivial induction based on the fact that $\lambda^{[1]}$ is linear and $C^0$ by Eqn. (3).

We shall show presently that compositions of composable $C^k$-maps are $C^k$. The following auxiliary result will be used:

**Lemma 4.4** *Let $k \in \mathbb{N}$, $f_i\colon E_i \supseteq U_i \to F_i$ be a $C^k$-map for $i \in \{1, 2\}$, $U \subseteq E$ an open subset, and $\lambda = (\lambda_1, \lambda_2)\colon E \to E_1 \times E_2$ a linear $C^0$-map such that $\lambda(U) \subseteq U_1 \times U_2$. Then*

$$g := (f_1 \times f_2) \circ \lambda|_U \colon\ U \to F_1 \times F_2$$

*is of class $C^k$. In particular, $f_1 \times f_2$ is $C^k$ if so are $f_1$ and $f_2$.*

**Proof.** The case $k = 1$ is immediate from **3.7**, Example 4.3 and Prop. 3.1, which also show that $g^{[1]}(x, y, t) = (f_1^{[1]}(\lambda_1(x), \lambda_1(y), t), f_2^{[1]}(\lambda_2(x), \lambda_2(y), t))$ for all $(x, y, t) \in U^{[1]}$. By the preceding formula, $g^{[1]}$ is the composition of the $C^{k-1}$-map $f_1^{[1]} \times f_2^{[1]}$ and the restriction to $U^{[1]}$ of the linear $C^0$-map $E^{[1]} \to E_1^{[1]} \times E_2^{[1]}$, $(x, y, t) \mapsto (\lambda_1(x), \lambda_1(y), t, \lambda_2(x), \lambda_2(y), t)$ and hence of class $C^{k-1}$ by induction. Consequently, $g$ is of class $C^k$. Choosing $\lambda := \mathrm{id}_{E_1 \times E_2}$, the final assertion follows. $\square$

**Proposition 4.5** *If $g$ and $f$ are $C^k$ and composable, then $g \circ f$ is $C^k$.*

**Proof.** We may assume that $k \in \mathbb{N}_0$; the proof is by induction on $k$. The case $k = 0$ is covered by Axiom I, and the case $k = 1$ has been treated in Prop. 3.1. Thus, let $k \geq 2$ now, and suppose that $g$ and $f$ (with all the notation as in Prop. 3.1) are $C^k$ and composable. According to Eqn. (5), $(g \circ f)^{[1]}$ can be written as a composition of $g^{[1]}$, $f^{[1]}$ and diagonal maps $\delta$:

$$(g \circ f)^{[1]} = g^{[1]} \circ (f \times f^{[1]} \times \mathrm{id}_{\mathbb{K}}) \circ (\delta_E \times \mathrm{id}_E \times \delta_{\mathbb{K}})$$

(where the second and last mapping have to be restricted and co-restricted properly, of course). Since $\delta_E \times \mathrm{id}_E \times \delta_{\mathbb{K}}$ is a linear $C^0$-map (cf. Remark 1.2) and $f$, $f^{[1]}$ and $\mathrm{id}_{\mathbb{K}}$ are of class $C^{k-1}$, Lemma 4.4 shows that $(f \times f^{[1]} \times \mathrm{id}_{\mathbb{K}}) \circ (\delta_E \times \mathrm{id}_E \times \delta_{\mathbb{K}})$ is of class $C^{k-1}$. Since $(g \circ f)^{[1]}$ is a composition of the $C^{k-1}$-map $g^{[1]}$ and the mapping just recognized as $C^{k-1}$, the induction hypothesis shows that $(g \circ f)^{[1]}$ is $C^{k-1}$. Thus $g \circ f$ is of class $C^k$. □

The arguments from the preceding proof also show that the iterated tangent maps $\hat{T}^k f := \hat{T} \cdots \hat{T} f$ are $C^0$. Moreover, one can prove that $f^{[k]}$ is obtained by composing $\hat{T}^k$ with a projection on some factor, and this yields further information on the maps $f^{[k]}$ (see [8]). – Our next aim is to define the $k$-th order differential $d^k f(x)$ of $f$ at $x$. If $f : U \to F$ is $C^2$, then the maps $df : U \times E \to F$ and $\partial_v f = df(\cdot)v : U \to F$, being partial maps of $f^{[1]}$, are $C^1$. In particular, for $v, w \in E$ fixed,

$$\partial_v(\partial_w f) : U \to F$$

is defined and is of class $C^0$.

**Lemma 4.6 ("Schwarz' Lemma")** *If $f$ is $C^2$, then, for all $v, w \in E$,*

$$\partial_v \partial_w f = \partial_w \partial_v f.$$

**Proof.** Recall the notation $\Delta_{v,t} f(x)$ (Notation 2.4). For $s, t \in \mathbb{K}^\times$,

$$\begin{aligned}
\Delta_{v,s}\left(\Delta_{w,t} f\right)(x) &= \frac{\Delta f^{[1]}(\cdot, w, t)}{\Delta}(x, v, s) \\
&= \frac{f^{[1]}(x+sv, w, t) - f^{[1]}(x, w, t)}{s} \\
&= \frac{f(x+sv+tw) - f(x+sv) - f(x+tw) + f(x)}{st}.
\end{aligned}$$

It follows that for all $t = s \in \mathbb{K}^\times$,

$$\Delta_{v,t}\left(\Delta_{w,t} f\right)(x) = \Delta_{w,t}\left(\Delta_{v,t} f\right)(x).$$

Since both sides are, for fixed $x, v, w$, $C^0$-maps of $t$, it follows that also for $t = 0$ we have equality, whence $\partial_v \partial_w f(x) = \partial_w \partial_v f(x)$. □

Note that if $f$ is $C^n$ ($n \geq 1$), then $\partial_v f = f^{[1]}(\bullet, v, 0)$ is $C^{n-1}$. This facilitates:

**Definition 4.7** If $f$ is $C^2$, we let

$$(d^2 f(x))(v, w) := d(df)(x, w) \cdot (v, 0) = \partial_v \partial_w f(x),$$

and if $f$ is $C^n$, $n \geq k$,

$$(d^k f(x))(v_1, \ldots, v_k) := \partial_{v_1} \ldots \partial_{v_k} f(x)$$

and

$$(d^k f(x))v := (d^k f(x))(v, \ldots, v).$$

**Lemma 4.8** $(d^k f)(x) \colon E^k \to F$ *is a symmetric multilinear $C^0$-map.*

**Proof.** By Lemma 4.6, $d^2 f(x)$ is symmetric, and inductively it follows that also $d^k f(x)$ is symmetric. Since $\partial_h$ is linear in $h$, it follows that $d^k f(x)$ is linear in the first argument, and the symmetry implies that it is linear in all arguments. It remains to show that $d^k f(x)$ is $C^0$. In order to prove this, it suffices to show that $d^k f(x)$ is a partial map of $f^{[k]}$. In fact, this is a straightforward consequence of the definitions. For later use, let us make this more explicit: for $k = 2$ and $s, t \in \mathbb{K}^\times$, we have

$$\Delta_{v,s}(\Delta_{w,t} f)(x) = \frac{f(x + sv + tw) - f(x + sv) - f(x + tw) + f(x)}{st}$$
$$= \frac{1}{s}(f^{[1]}(x + sv, w, t) - f^{[1]}(x, w, t))$$
$$= f^{[2]}((x, w, t), (v, 0, 0), s).$$

For $s = t = 0$, this gives

$$d^2 f(x)(v, w) = \partial_v \partial_w f(x) = f^{[2]}((x, w, 0), (v, 0, 0), 0),$$

proving our claim for $k = 2$. By induction, we get for $k \in \mathbb{N}$:

$$(\Delta_{v_1, t_1} \cdots \Delta_{v_k, t_k} f)(x) = \frac{1}{\prod_i t_i} \left( f\left(x + \sum_i t_i v_i\right) - \sum_j f\left(x + \sum_{i; i \neq j} t_i v_i\right) + \right.$$
$$\left. \sum_{j,l: j < l} f\left(x + \sum_{i: i \neq j, i \neq l} t_i v_i\right) - \ldots + (-1)^{k-1} \sum_i f(x + t_i v_i) + (-1)^k f(x) \right).$$

In particular, if all $v_i$ and $t_i$ coincide, we get

$$(\Delta_{v,t})^k f(x) = \frac{1}{t^k} \sum_{j=0}^{k} \binom{k}{j} (-1)^j f(x + (k-j)tv)$$

and for $t = 0$ we may write:

$$\partial_v^k f(x) = \lim_{t \to 0} \frac{f(x + ktv) - k f(x + (k-1)tv) + \ldots + (-1)^k f(x)}{t^k}.$$

All these expressions are partial maps of $f^{[k]}$ and hence are $C^0$. □

It is also possible to express $d^2 f$ via the map

$$d(df) \colon (U \times E) \times (E \times E) \to F.$$

In fact, $df = f^{[1]}(\cdot, \cdot, 0)$, and hence

$$ddf(x, w) \cdot (v, h) = \frac{\Delta f^{[1]}}{\Delta}((x, w, 0), (v, h, 0), 0),$$

and thus

$$\partial_v \partial_w f(x) = ddf(x, w) \cdot (v, 0).$$

Finally we discuss locality of the property of being $C^k$.

**Lemma 4.9** *Assume $\mathbb{K}$ is a field and $f\colon U \to F$ a mapping, defined on an open subset $U$ of $E$. Let $k \in \mathbb{N}_0 \cup \{\infty\}$. If there is an open cover $(U_i)_{i \in I}$ of $U$ such that $f|_{U_i}\colon U_i \to F$ is $C^k$ for each $i \in I$, then $f$ is of class $C^k$.*

**Proof.** It suffices to prove the assertion for $k \in \mathbb{N}_0$; the proof is by induction. The case $k = 0$ is incorporated in our setting of differential calculus (Locality Axiom I.5).

*Induction step.* Suppose the assertion holds for some $k \in \mathbb{N}_0$, and suppose that $f|_{U_i}$ is of class $C^{k+1}$ for all $i \in I$, for some open cover $(U_i)_{i \in I}$ of $U$. Then $f$ is of class $C^k$, by induction. We now define a mapping $g\colon U^{[1]} \to F$, which will turn out to be $f^{[1]}$. Let $(x, y, t) \in U^{[1]}$. If $t \neq 0$, we define

$$g(x, y, t) := \tfrac{1}{t}(f(x + ty) - f(x)). \tag{10}$$

If $t = 0$ and $x \in U_i$, we define

$$g(x, y, 0) := d(f|_{U_i})(x, y). \tag{11}$$

If $x \in U_i \cap U_j$ here for certain $i, j \in I$, then $(f|_{U_i})^{[1]}|_{(U_i \cap U_j)^{[1]}} = (f|_{U_i \cap U_j})^{[1]} = (f|_{U_j})^{[1]}|_{(U_i \cap U_j)^{[1]}}$ apparently, entailing that $d(f|_{U_i})(x, y) = d(f|_{U_i \cap U_j})(x, y) = d(f|_{U_j})(x, y)$. Consequently, $g$ is well-defined. We have $g|_{U_i^{[1]}} = (f|_{U_i})^{[1]}$, showing that $g$ is of class $C^k$ on the open set $U_i^{[1]}$, for each $i \in I$. On the other hand, the map $f$ being of class $C^k$ by induction, it is clear from Equation (10) that $g$ is of class $C^k$ on the open subset $W := \{(x, y, t) \in U^{[1]}\colon t \in \mathbb{K}^\times\}$ of $U^{[1]}$. As the sets $U_i^{[1]}$, together with $W$, form an open cover of $U^{[1]}$, by induction the mapping $g$ is of class $C^k$ and thus $C^0$ in particular. Thus $f$ is of class $C^1$, with $f^{[1]} = g$ of class $C^k$. As a consequence, $f$ is of class $C^{k+1}$. $\square$

## 5 Taylor's Formula

Assume $f\colon E \supseteq U \to F$ is a map of class $C^k$. A first order Taylor expansion of the form $f(x + th) = f(x) + tdf(x)h + tO(t)$ has already been given in Remark 2.5. Let us iterate the calculation given there: we fix $x \in U$ and $h \in E$; then, for all $t$ such that $x + th \in U$,

$$\begin{aligned}
f(x + th) &= f(x) + tf^{[1]}(x, h, t) \\
&= f(x) + tdf(x)h + t(f^{[1]}(x, h, t) - f^{[1]}(x, h, 0)) \\
&= f(x) + tdf(x)h + t^2 f^{[2]}((x, h, 0), (0, 0, 1), t) \\
&= f(x) + tdf(x)h + t^2 f^{[2]}((x, h, 0), (0, 0, 1), 0) \\
&\quad + t^2 \left(f^{[2]}((x, h, 0), (0, 0, 1), t) - f^{[2]}((x, h, 0), (0, 0, 1), 0)\right).
\end{aligned}$$

The last term is a product of $t^2$ and a $C^0$-map taking the value $0$ at $t = 0$; following our convention from Remark 2.5, we say that it is $t^2 O(t)$. It is clear that this procedure can be iterated $k$ times, and we get the following $k$-th order expansion of $f$ (compare with the "développement limité" in Chapter 7 of [17]):

**Theorem 5.1** *If $f\colon U \to F$ is $C^k$, then for all $(x, h, t) \in U^{[1]}$, the following holds:*

$$f(x + th) = \sum_{j=0}^{k} t^j a_j(x, h) + t^k R_{k+1}(x, h, t),$$

*where $a_j\colon U \times E \to F$ is of class $C^{k-j}$ and $R_{k+1}\colon U^{[1]} \to F$ is of class $C^0$ and takes the value $0$ for $t = 0$. An expansion of $f$ with these properties is unique. Moreover, $a_j(x, h)$ is homogeneous of degree $j$ in $h$.*

**Proof.** Existence of the expansion has already been proved in the case $k = 2$, and iterating the argument in the obvious way, we get an expansion of the desired form for arbitrary $k$. The uniqueness assertion is covered by the following lemma. The expansion being unique, we readily deduce that $a_j(x, th) = t^j a_j(x, h)$ for all $t \in \mathbb{K}$. $\square$

**Lemma 5.2** *Suppose that $I \subseteq \mathbb{K}$ is an open zero-neighbourhood, $k \in \mathbb{N}_0$, and $f\colon I \to F$ a mapping of the form*

$$f(t) := \sum_{j=0}^{k} t^j a_j + t^k R(t) \quad \text{for } t \in I, \tag{12}$$

*where $a_0, \ldots, a_k \in F$ and where $R\colon I \to F$ is $C^0$ with $R(0) = 0$. Then the elements $a_0, a_1, \ldots, a_k$ are uniquely determined by $f$. If $f$ is homogeneous of degree $p \leq k$ in the sense that $f(st) = s^p f(t)$ for all $t \in I$ and $s \in \mathbb{K}$ such that $st \in I$, then $R = 0$ and $a_j = 0$ for all $j \neq p$.*

**Proof.** Suppose that also $f(t) = \sum_{j=0}^{k} t^j a'_j + t^k R'(t)$; then $\sum_{j=0}^{k} t^j (a_j - a'_j) + t^k (R(t) - R'(t))$ is an expansion for the function $I \to F$, $t \mapsto 0$. To prove uniqueness of the expansion, it therefore suffices to assume that $f = 0$ and show that this entails that $a_j = 0$ for all $j = 0, 1, \ldots, k$, whenever we have an expansion (12). Taking $t = 0$, we obtain $0 = f(0) = a_0$. Given $t \in I \cap \mathbb{K}^\times$, dividing both sides of (12) by $t$ we obtain

$$0 = \sum_{j=0}^{k-1} t^j a_{j+1} + t^{k-1} R(t),$$

which then actually holds for all $t \in I$, by the Determination Axiom. By induction, the preceding equation entails that $a_1, \ldots, a_k = 0$.

To prove the final assertion, suppose that $f$ is homogeneous of degree $p$ and admits the expansion (12). Given $t \in I$, there exists an open zero-neighbourhood $J \subseteq I$ such that $J \cdot t \subseteq I$. Then

$$\sum_{j=0}^{k} s^j t^j a_j + s^k t^k R(st) = f(st) = t^p f(s) = \sum_{j=0}^{k} s^j t^p a_j + s^k t^p R(s)$$

for all $s \in J$, whence $t^j a_j = t^p a_j$ for all $j = 0, \ldots, k$, by uniqueness of the expansion of the function $J \to F$, $s \mapsto f(st)$. Thus, for fixed $j \in \{0, \ldots, k\}$, we have $g(t) := t^j a_j = t^p a_j$ for all $t \in I$. If $j \neq p$, then the uniqueness of the expansion for $g$ entails that $a_j = 0$. Thus

$$f(t) = t^p a_p + t^k R(t). \tag{13}$$

Given $t \in I \cap \mathbb{K}^\times$, set $J := \frac{1}{t} I$; then $J$ is an open zero-neighbourhood in $\mathbb{K}$. We have $s^p f(t) = f(st) = s^p t^p a_p + s^k t^k R(st)$ for all $s \in J$ and hence

$$f(t) = t^p a_p + s^{k-p} t^k R(st)$$

for all $s \in J \cap \mathbb{K}^\times$. By the Determination Axiom, both sides actually coincide for all $s \in J$. Comparing this with (13), we find that $s^{k-p} t^k R(st) = t^k R(t)$ for all $s \in J$. Setting $s = 0$, we obtain $0 = 0^{k-p} t^k R(0) = t^k R(t)$, whence $R(t) = 0$. The map $R$ being $C^0$, the Determination Axiom shows that $R(t) = 0$ for all $t \in I$, whence $f(t) = t^p a_p$. $\square$

Our next task is to identify and describe the coefficients $a_j(x, h)$. For $j = 0, 1, 2$, we have the following:

**Proposition 5.3** *Under the assumptions of Theorem 5.1, $a_0(x, h)$ is constant and $a_1(x, h)$ is linear in $h$ and $a_2(x, h)$ is an $F$-valued quadratic form in $h$. More precisely,*

$$\begin{aligned} a_0(x, h) &= f(x), \\ a_1(x, h) &= df(x)h, \\ a_2(x, h_1 + h_2) - a_2(x, h_1) - a_2(x, h_2) &= d^2 f(x)(h_1, h_2), \end{aligned}$$

*and, in particular,* $2 a_2(x, h) = d^2 f(x)(h, h)$.

**Proof.** The determination of $a_0$ is clear and the determination of $a_1$ is given by the calculation preceding Theorem 5.1. Recall from [14] that an $F$-valued quadratic form is a map $q \colon E \to F$ which is homogeneous of degree two and such that $q(h_1 + h_2) - q(h_1) - q(h_2)$ is $\mathbb{K}$-bilinear in $h_1, h_2$. Since we already know that $a_2(x, \cdot)$ is homogeneous of degree two, it suffices to prove the third equality stated in the proposition since the right hand side is bilinear in $h_1, h_2$ by Lemma 4.8. We use the expansion from Theorem 5.1 and write $O(t)$ for terms which are $C^0$ and take the value 0 for $t = 0$:

$$\begin{aligned} &f(x + t(h_1 + h_2)) - f(x + th_1) - f(x + th_2) + f(x) \\ &= t\,(df(x)(h_1 + h_2) - df(x)h_1 - df(x)h_2) + t^2(a_2(x, h_1 + h_2) - a_2(x, h_1) - a_2(x, h_2)) \\ &\quad + t^2 O(t) \\ &= t^2(a_2(x, h_1 + h_2) - a_2(x, h_1) - a_2(x, h_2)) + t^2 O(t). \end{aligned}$$

Thus, for $t \in \mathbb{K}^\times$,

$$\frac{f(x + t(h_1 + h_2)) - f(x + th_1) - f(x + th_2) + f(x)}{t^2}$$
$$= a_2(x, h_1 + h_2) - a_2(x, h_1) - a_2(x, h_2) + O(t).$$

For $t = 0$, the $C^0$-extension of the left hand side equals $d^2 f(x)(h_1, h_2)$ (cf. proof of La. 4.6), and this implies our claim for $a_2$. As we already know that $a_2(x, h)$ is homogeneous quadratic in $h$, we get for $h_1 = h_2 = h$ the relation $2a_2(x, h) = d^2 f(x)(h, h)$. $\square$

The preceding arguments can be generalized to arbitrary degree. The definition and some elementary facts on *F-valued forms of degree j* are given in an appendix (Appendix A).

**Theorem 5.4** *In the situation of Theorem 5.1, the maps $a_j(x, h)$ are F-valued forms in $h$ of degree $j$. The linearization of $a_j(x, \cdot)$ of order $j$ satisfies the relation*

$$(L_j a_j(x, \cdot))(0; h_1, \ldots, h_j) = d^j f(x)(h_1, \ldots, h_j),$$

*and in particular*

$$j!\, a_j(x, h) = d^j f(x)(h, \ldots, h).$$

*Therefore, if $2, 3, \ldots, k$ are invertible in $\mathbb{K}$, we have the Taylor expansion*

$$f(x + th) = \sum_{j=0}^{k} \frac{t^j}{j!} d^j f(x)(h, \ldots, h) + t^k R_{k+1}(x, h, t).$$

**Proof.** The claim is proved by induction on $j$, the cases $j = 0, 1, 2$ being already proved. After performing a translation, we may assume that $x = 0$. We let $b_j(h) := a_j(0, h)$ and re-write the expansion from Theorem 5.1 as

$$f(th) = \sum_{i=0}^{k} t^i b_i(h) + t^k R_{k+1}(0, h, t) = \sum_{i=0}^{k} b_i(th) + t^k R_{k+1}(0, h, t).$$

We apply the linearization operator $L_j$ with respect to the variable $h$ to this equation:

$$L_j f(0; th_1, \ldots, th_j) = \sum_{i=0}^{k} L_j b_i(0; th_1, \ldots, th_j) + t^k O(t).$$

For $i < j$, by induction hypothesis $b_i$ is a form of degree $i$ and hence $L_j b_i$ vanishes, by Proposition A.3. The remaining terms give

$$L_j f(0; th_1, \ldots, th_j) = L_j b_j(0; th_1, \ldots, th_j) + t^j O(t) = t^j (L_j b_j(0; h_1, \ldots, h_j) + O(t)).$$

We divide by $t^j$ and let $t = 0$. On the left hand side we find precisely the expression for $d^j f(0)(h_1, \ldots, h_j)$ (cf. proof of Lemma 4.8). From Lemma 4.8, we know that this is $\mathbb{K}$-linear in all arguments, hence so is the right hand side, and finally

$$d^j f(0)(h_1, \ldots, h_j) = L_j b_j(0; h_1, \ldots, h_j).$$

According to Proposition A.3, $b_j(h)$ is therefore a form in $h$ of degree $j$, and using Lemma A.4 we see that $d^j f(0)(h, \ldots, h) = j! \, b_j(h)$. If the integers $2, 3, \ldots, k$ are invertible, Taylor's formula in its usual form is now an immediate consequence. □

Note that as a corollary we get a collection of non-trivial identities for the higher-order difference quotient maps $f^{[k]}$. Namely, on the one hand we have

$$(d^2 f(x))(v, w) = f^{[2]}((x, w, 0), (v, 0, 0), 0)$$

(see proof of Lemma 4.8). On the other hand, the calculation preceding Theorem 5.1 shows that

$$a_2(x, h) = f^{[2]}((x, h, 0), (0, 0, 1), 0),$$

and thus Proposition 5.3 implies the non-trivial identity

$$f^{[2]}((x, h, 0), (h, 0, 0), 0) = 2 f^{[2]}((x, h, 0), (0, 0, 1), 0).$$

In a similar way, Theorem 5.4 implies identities for the higher $f^{[k]}$'s which, however, are too complicated to be written out here.

In case of characteristic zero, Theorem 5.4 implies in particular that $a_j(x, h)$ is a polynomial in $h$ which can be defined by a *symmetric* multilinear form (namely $\frac{1}{j!} d^j f(x)$). We are going to prove now that, even in the case of positive characteristic, $a_j(x, \cdot)$ still is a polynomial map (as defined in A.5); however, in the general case there seems to be no *canonical* way to write it as a polynomial.

**Lemma 5.5** *Assume $q : E \to F$ is homogeneous of degree $k$, i.e. $q(tx) = t^k q(x)$ for all $t \in \mathbb{K}$, $x \in E$. If $q$ is $C^k$, then $q$ is a homogeneous polynomial map of degree $k$.*

**Proof.** Let $(e_i)_{i \in I}$ be a system of generators of $E$, $\tilde{E}$ the free module with basis $(e_i)_{i \in I}$ and $\phi : \tilde{E} \to E$ the surjection defined by $\phi(e_i) = e_i$. We write $x \in \tilde{E}$ as $x = \sum_i t_i e_i$ (finite sum). Then it is proved by induction on the number of non-zero terms in this expression that (using the multi-index notation from A.5) $q(\phi(\sum_i t_i e_i)) = \sum_{|\alpha|=k} t^\alpha a_\alpha$ with coefficients $a_\alpha$ depending on the system of generators, but not on $x$. For simplicity of notation, let us assume that $x = t_1 e_1 + t_2 e_2$, the general case being similar. Then we get, using Theorem 5.1 repeatedly, first to expand $q$ around $t_2 e_2$, then to expand each $a_i$ around $(0, e_1)$:

$$
\begin{aligned}
q(t_2 e_2 + t_1 e_1) &= \sum_{i=0}^{k} t_1^i a_i(t_2 e_2, e_1) + t_1^k R_{k+1}(t_2 e_2, e_1, t_1) \\
&= \sum_{i=0}^{k} t_1^i \left( \sum_{j=0}^{k-i} t_2^j a_{ij} + t_2^{k-i} R_{k-i+1}((0, e_1), (e_2, 0), t_2) \right) + t_1^k R_{k+1}(t_2 e_2, e_1, t_1) \\
&= \sum_{i+j \leq k} t_1^i t_2^j a_{ij} + \sum_{i+j=k} t_1^i t_2^j R_{k-i+1}((0, e_1), (e_2, 0), t_2) + t_1^k R_{k+1}(t_2 e_2, e_1, t_1),
\end{aligned}
$$

for suitable $a_{ij} \in F$. Now replace $t_1e_1+t_2e_2$ with $s(t_1e_1+t_2e_2)$ where $s \in \mathbb{K}$. By assumption, the left hand side is homogeneous of degree $k$ and hence will simply be multiplied by a factor $s^k$. On the right hand side we have a sum of homogeneous terms of degree $\ell$, $\ell = 0, \ldots, k$, and the remainder term which is $s^k O(s)$. By Lemma 5.2, all terms on the right hand side except for the term of degree $k$ vanish. Thus $s^k q(t_1e_1 + t_2e_2) = s^k \sum_{i+j=k} t_1^i t_2^j a_{ij}$; letting $s = 1$, the claim follows. □

**Theorem 5.6** *Assume $f\colon E \supseteq U \to F$ is of class $C^{2k}$ and $x, x + h \in U$.*

(a) *The "regular part" $\sum_j a_j(x, h)$ in the expansion*

$$f(x+h) = \sum_{j=0}^{k} a_j(x, h) + R_{k+1}(x, h, 1)$$

*from Theorem 5.1 is a polynomial map in $h$ of degree at most $k$.*

(b) *$f$ is (the restriction to $U$ of) a polynomial map of degree at most $k$ if and only if $f$ coincides with the regular part of its Taylor expansion about one (and hence any) point $x \in U$.*

(c) *Assume $U$ is a neighbourhood of $0$ and $f \in C^k(U, F)$ is homogeneous of degree $p$, $p \in \mathbb{N}$, $p \leq k$, i.e. $f(th) = t^p f(h)$ whenever $h, th \in U$. Then $f$ is (the restriction of) a homogeneous polynomial of degree $p$.*

**Proof.** (a) Using 5.5, we see that each term $a_j(x, \cdot)$ (which is homogeneous of degree $j$ and of class $C^{2k-j}$ and hence of class $C^j$) is a polynomial map as defined in A.5, and hence the regular part is polynomial.

(b) If $f$ is polynomial, then $f(x+h)$ is again a polynomial $p(h)$ in $h$, and $f(x+h) = p(h) + 0$ is an expansion satisfying the properties of Theorem 5.1. By uniqueness of the expansion, it follows that $f$ is given by its Taylor expansion about $x$. Conversely, if a map $f$ is given by its regular part, then by (a) it is polynomial.

(c) By Lemma 5.2, we have $f(h) = a_p(0, h)$. Thus $a_p(0, \cdot)|_U$ is of class $C^k$ and hence so is $a_p(0, \cdot)$, using that $U$ is absorbing (Remark 1.4) and $a_p(0, x) = s^p a_p(0, s^{-1}x)$ for all $s \in \mathbb{K}^\times$ and $x \in sU$, whence $a_p(0, \cdot)|_{sU}$ is $C^k$. Now Lemma 5.5 shows that $a_p(0, h)$ is indeed polynomial in $h$. □

# 6  $C^k$-curves

In this section, we provide a simpler, equivalent description of the $C^k$-property for curves, which mimics a definition used in Schikhof's book [58] on ultrametric calculus, formulated there for mappings between subsets of an ultrametric field. The alternative description will be extremely useful later, when we need to construct smooth curves with very special properties. When $\mathbb{K}$ is a ring, the alternative condition is still necessary for being $C^k$, but to prove equivalence we shall need to assume that $\mathbb{K}$ is a field.

A *curve of class* $C^k$ is a map $f: U \to E$ of class $C^k$, where $E \in \mathcal{M}$ and $U \subseteq \mathbb{K}$ is an open subset. For curves, our formalism simplifies: the map $f^{[1]}(x, h, t)$ may be replaced with $f^{[1]}(x, 1, t)$. More precisely, we have, after performing an affine change of coordinates:

**Lemma 6.1** *A curve* $f: U \to E$ *is of class* $C^1$ *if and only if there exists a map* $f^{<1>}: U \times U \to E$ *of class* $C^0$ *such that, for all* $s, t \in U$,

$$f(s) - f(t) = f^{<1>}(s, t) \cdot (s - t). \tag{14}$$

*Then* $df(t) = f^{<1>}(t, t)$, *and* $df: U \to E$ *is a curve of class* $C^0$.

**Proof.** Assume $f$ is $C^1$. We let $f^{<1>}(s, t) := f^{[1]}(t, 1, s-t)$; then $f^{<1>}: U \times U \to E$ is well-defined and $C^0$, and (14) is nothing but the relation $f(t+(s-t)) - f(t) = (s-t)f^{[1]}(t, 1, s-t)$. Conversely, given $f^{<1>}$ satisfying (14), we let $F(x, h, t) := hF(x, 1, th) := hf^{<1>}(x+th, x)$; this is $C^0$, and (14) implies that $f$ is of class $C^1$ with $f^{[1]} = F$.

It follows that $f^{<1>}(t, t) = f^{[1]}(t, 1, 0) = df(t)1$, and, identifying $df(t) \in \mathrm{Hom}_{\mathbb{K}}(\mathbb{K}, E)$ with its value $df(t)1 \in E$, $df: U \to E$ is a curve of class $C^0$. $\square$

Next, we wish to generalize the characterization of Lemma 6.1 to curves $f$ of class $C^k$, where $k < \infty$. Having defined $f^{<1>}$ as above, we define by induction maps

$$f^{<j>}: U^{j+1} \to E$$

of class $C^{k-j}$ by:

$$\begin{aligned} f^{<j+1>}(t_1, \ldots, t_{j+1}) &:= \frac{\Delta f^{<j>}(\cdot, t_3, t_4, \ldots, t_{j+1})}{\Delta}(t_1, 1, t_2 - t_1) \\ &= \frac{f^{<j>}(t_1, t_3, \ldots, t_{j+1}) - f^{<j>}(t_2, t_3, \ldots, t_{j+1})}{t_1 - t_2}, \end{aligned}$$

where the second expression is valid only if $t_1 - t_2 \in \mathbb{K}^\times$, whereas the first expression is valid in general. By induction, one can prove the following explicit formula for $f^{<k>}$ (cf. [58, Exercise 29.A]): If $t_i - t_j \in \mathbb{K}^\times$ for $i \neq j$, then

$$f^{<k>}(t_1, \ldots, t_{k+1}) = \sum_{j=1}^{k+1} \frac{f(t_j)}{\prod_{i \neq j}(t_j - t_i)}.$$

Then the following holds:

**Proposition 6.2** *If* $f: \mathbb{K} \supseteq U \to E$ *is a curve of class* $C^k$, *then there exists a map*

$$f^{<k>}: U^{k+1} \to E$$

of class $C^0$ such that, on the open subset $U^{>k<}$ of $U^{k+1}$ defined by the condition $\forall i \neq j$: $t_j - t_i \in \mathbb{K}^\times$, we have

$$f^{<k>}(t_1, \ldots, t_{k+1}) = \sum_{j=1}^{k+1} \frac{f(t_j)}{\prod_{i \neq j}(t_j - t_i)}. \tag{15}$$

Moreover, for any map $f^{<k>}$ having these properties, the relation

$$d^j f(t) = j!\, f^{<j>}(t, \ldots, t) \tag{16}$$

holds, and, if $0 \in U$, the expansion of $f$ at the origin given by Theorem 5.1 can be written

$$f(t) = \sum_{j=0}^{k} t^j f^{<j>}(0, \ldots, 0) + t^k R_{k+1}(0, 1, t). \tag{17}$$

**Proof.** Existence of $f^{<k>}$ is proved by the inductive definition given above. (We do not claim that $f^{<k>}$ is unique; this is so if $\mathbb{K}$ is a field since then $U^{>k<}$ is $C^0$-dense in $U^{k+1}$, as is easily seen.) In order to prove (17), we simply re-write the calculation given before Theorem 5.1 by using the maps $f^{<j>}$:

$$\begin{aligned}
f(t) &= f(0) + tf^{[1]}(0, 1, t) = f(0) + tf^{<1>}(t, 0) \\
&= f(0) + tf^{<1>}(0, 0) + t(f^{<1>}(t, 0) - f^{<1>}(0, 0)) \\
&= f(0) + tf^{<1>}(0, 0) + t^2 f^{<2>}(t, 0, 0) \\
&= f(0) + tf^{<1>}(0, 0) + t^2 f^{<2>}(0, 0, 0) + t^2(f^{<2>}(t, 0, 0) - f^{<2>}(0, 0, 0)) \\
&= f(0) + tf^{<1>}(0, 0) + t^2 f^{<2>}(0, 0, 0) + t^3 f^{<3>}(t, 0, 0, 0)
\end{aligned}$$

and so on, leading to (17) with

$$R_{k+1}(0, 1, t) = f^{<k>}(t, 0, \ldots, 0) - f^{<k>}(0, 0, \ldots, 0)$$

which clearly is $O(t)$. It follows that $a_j(0, 1) = f^{<j>}(0, \ldots, 0)$ and, by translation, $a_j(t, 1) = f^{<j>}(t, \ldots, t)$. Comparing with Theorem 5.4, we get (16).[1] Finally, assume $\tilde{f}^{<k>}$ is another map having the properties of $f^{<k>}$ mentioned in the proposition. Fix some vector $v \in U^{>k<}$. Then $tv \in U^{>k<}$ for all $t \in \mathbb{K}^\times$ sufficiently close to 0 (using that $t_i - t_j \in \mathbb{K}^\times$ entails $tt_i - tt_j = t(t_i - t_j) \in \mathbb{K}^\times$). Let $\alpha(t) := f^{<k>}(tv) - \tilde{f}^{<k>}(tv)$. Then $\alpha$ is $C^0$ on some neighbourhood of the origin in $\mathbb{K}$ and vanishes for all invertible $t$ in this neighbourhood since $f^{<k>}|_{U^{>k<}} = \tilde{f}^{<k>}|_{U^{>k<}}$. By the Determination Axiom, $\alpha$ vanishes on the neighbourhood, whence in particular $\alpha(0) = 0$. Thus $f^{<j>}(0, \ldots, 0)$ is indeed independent of $f^{<j>}$. $\square$

In the remainder of this section, we assume that the topologized ring $\mathbb{K}$ is a field. Our goal is to show that, in this case, also the converse of Proposition 6.2 holds: whenever $f^{<k>}$ exists, $f$ will be a $C^k$-curve. Stimulated by W. H. Schikhof's definition of $C^k$-maps between subsets of ultrametric fields [58], we make the following definition:

---

[1] Cf. [58, Thm. 29.5] for a slightly more direct, but essentially not much different proof of (16).

**Definition 6.3** Let $U$ be an open, non-empty subset of $\mathbb{K}$, $E \in \mathcal{M}$ be a topologized $\mathbb{K}$-vector space, and $f\colon U \to E$ a map. For $n \in \mathbb{N}_0$, define

$$U^{>n<} := \{(x_1, \ldots, x_{n+1}) \in U^{n+1} \colon \text{if } i \neq j \text{ then } x_i \neq x_j\},$$

set $f^{>0<} := f\colon U = U^{>0<} \to E$, and recursively define $f^{>n<}\colon U^{>n<} \to \mathbb{K}$ for $n \in \mathbb{N}$ via

$$f^{>n<}(x_1, \ldots, x_{n+1}) := \frac{f^{>n-1<}(x_1, x_3, \ldots, x_{n+1}) - f^{>n-1<}(x_2, x_3, \ldots, x_{n+1})}{x_1 - x_2}.$$

The function $f$ is called a $C^n_{Sch}$-map if $f^{>n<}$ can be extended to a $C^0$-map $f^{<n>}\colon U^{n+1} \to E$. We call $f$ a $C^\infty_{Sch}$-map if it is a $C^n_{Sch}$-map for each $n \in \mathbb{N}_0$.

It is easy to see that $U^{>n<}$ is $C^0$-dense in $U^{n+1}$, and apparently $U^{>n<}$ is open. Thus, if it exists, $f^{<n>}$ is uniquely determined.

**Remark 6.4** To see that $U^{>n<}$ is $C^0$-dense, note that $t + rs \in U^{>n<}$ for all $r \in \mathbb{K}^\times$ sufficiently close to $0$, for each given element $t = (t_1, \ldots, t_{n+1}) \in U^{n+1}$ and $s := (s_1, \ldots, s_{n+1})$ with $s_1, \ldots, s_{n+1} \in \mathbb{K}$ pairwise distinct. Choosing $s_1, \ldots, s_{n+1} \in \mathbb{K}^\times$ here, we can furthermore achieve that $t + rs \notin F$, for any given finite subset $F \subseteq U^{n+1}$; this observation will be useful in the proof of Lemma 6.7.

**Remark 6.5** Our notation differs from Schikhof's, who writes $\nabla^{n+1}U$ for $U^{>n<}$, $\Phi_n f$ for $f^{>n<}$, and $\overline{\Phi}_n f$ for $f^{<n>}$.

The following lemmas will be needed:

**Lemma 6.6** *Let $U$ be an open non-empty subset of $\mathbb{K}$ and $f\colon U \to E$ be a mapping into a $\mathbb{K}$-vector space. Let $n \in \mathbb{N}$. Then the following holds:*

(a) $f^{>n<}\colon U^{>n<} \to E$ *is a symmetric function of its $n+1$ variables.*

(b) *For all $(y_1, \ldots, y_n, x_1, \ldots, x_n) \in U^{>2n-1<}_n \subseteq U^{2n}$, we have*

$$f^{>n-1<}(y_1, \ldots, y_n) - f^{>n-1<}(x_1, \ldots, x_n) = \sum_{j=1}^{n}(y_j - x_j) \cdot f^{>n<}(x_1, \ldots, x_j, y_j, \ldots, y_n). \quad (18)$$

**Proof.** The proof of [58], Lemma 29.2 (ii) and (iii) can be repeated verbatim. □

**Lemma 6.7** *If $f\colon U \to E$ is a $C^n_{Sch}$-map (where $n \in \mathbb{N}_0$), then $f$ is a $C^k_{Sch}$-map for all $k \in \mathbb{N}_0$ such that $k < n$. For all $1 \leq k \leq n$ and $(y_1, \ldots, y_k, x_1, \ldots, x_k) \in U^{2k}$, we have*

$$f^{<k-1>}(y_1, \ldots, y_k) - f^{<k-1>}(x_1, \ldots, x_k) = \sum_{j=1}^{k}(y_j - x_j) \cdot f^{<k>}(x_1, \ldots, x_j, y_j, \ldots, y_k). \quad (19)$$

*Furthermore, $f^{<k>}\colon U^{k+1} \to E$ is a symmetric function of its $k+1$ variables.*

**Proof.** The proof is by induction on $n \in \mathbb{N}_0$; the case $n = 0$ is trivial. Suppose that $n \in \mathbb{N}$ and suppose that the assertions of the lemma hold when $n$ is replaced with $n-1$. Let $f: U \to E$ be a $C^n_{Sch}$-map. Pick $y = (y_1, \ldots, y_n) \in U^{>n-1<}$. Then $D_y := \{x = (x_1, \ldots, x_n) \in U^{>n-1<} : (\forall i, j \in \{1, \ldots, n\}) \; x_i \neq y_j\}$ is easily seen to be $C^0$-dense in $U^{>n-1<}$ (cf. Remark 6.4). The map

$$\psi_y : U^n \to E, \quad \psi_y(x_1, \ldots, x_n) := f^{>n-1<}(y) + \sum_{j=1}^n (x_j - y_j) \cdot f^{<n>}(y_1, \ldots, y_j, x_j, \ldots, x_n)$$

is $C^0$, and $\psi_y|_{D_y} = f^{>n-1<}|_{D_y}$ by (18). If also $z \in U^{>n-1<} \subseteq U^n$, then $D_y \cap D_z$ is $C^0$-dense in $U^n$. Thus $\psi_y|_{D_y \cap D_z} = f^{>n-1<}|_{D_y \cap D_z} = \psi_z|_{D_y \cap D_z}$ entails that $\psi_y = \psi_z$. Thus $\psi := \psi_y$ is independent of the choice of $y \in U^{>n-1<}$. Given $x = (x_1, \ldots, x_n) \in U^{>n-1<}$, there exists $y \in D_x$. Then $x \in D_y$ and thus $f^{>n-1<}(x) = \psi_y(x) = \psi(x)$. We have shown that $\psi|_{U^{>n-1<}} = f^{>n-1<}$. Since $\psi$ is $C^0$, we deduce that $f$ is a $C^{n-1}_{Sch}$-map, and $f^{<n-1>} = \psi$. Thus, for $k = n$, both sides of (19) are $C^0$-functions on $U^{2n}$ which coincide on the $C^0$-dense subset $U^{>2n-1<}$ by (18), and which therefore coincide. Similarly, $f^{<n>}$ is a symmetric function of its $n+1$ variables since $f^{<n>} : U^{n+1} \to E$ is $C^0$ and its restriction $f^{>n<}$ to the $C^0$-dense subset $U^{>n<}$ of $U^{n+1}$ is symmetric. The assertions for $k < n$ follow from the induction hypothesis, since we have already shown that $f$ is a $C^{n-1}_{Sch}$-map. $\square$

We readily deduce:

**Lemma 6.8** *Let $f : U \to E$ be a $C^n_{Sch}$-map, where $n \in \mathbb{N}_0$. Then, for each $k \in \mathbb{N}_0$ such that $k \leq n$, the mapping $f^{<k>} : U^{k+1} \to E$ is of class $C^{n-k}$.*

**Proof.** We show by induction on $j = 0, \ldots, n$ the following claim (thus establishing the lemma):

*Claim. For each $k \in \mathbb{N}_0$ such that $k \leq n - j$, the map $f^{<k>}$ is of class $C^j_\mathbb{K}$.*

The case $j = 0$ is trivial, since all of the mappings $f^{<k>} : U^{k+1} \to E$ are $C^0$ (being $C^0$-extensions of the maps $f^{>k<}$ by definition).

*Induction step.* Suppose that the claim holds for some $j \in \mathbb{N}_0$ such that $j < n$, and suppose that $k \leq n - (j+1) = n - j - 1$. Let $x = (x_1, \ldots, x_{k+1})$, $y = (y_1, \ldots, y_{k+1}) \in \mathbb{K}^{k+1}$ and $t \in \mathbb{K}$ such that $(x, y, t) \in (U^{k+1})^{[1]}$. If $t \neq 0$, then

$$\tfrac{1}{t}(f^{<k>}(x+ty) - f^{<k>}(x)) = \sum_{i=1}^{k+1} y_i \cdot f^{<k+1>}(x_1, \ldots, x_i, x_i + ty_i, \ldots, x_{k+1} + ty_{k+1}) \quad (20)$$

by Equation (19). Note that $p_i : (\mathbb{K}^{k+1})^{[1]} = \mathbb{K}^{2(k+1)+1} \to \mathbb{K}$, $p_i(x, y, t) := y_i$ is a linear $C^0$-map and thus of class $C^\infty$; furthermore,

$$q_i : (\mathbb{K}^{k+1})^{[1]} \to \mathbb{K}^{k+2} : \quad q_i(x, y, t) := (x_1, \ldots, x_i, x_i + ty_i, \ldots, x_{k+1} + ty_{k+1})$$

is a polynomial function for $i = 0, \ldots, k+1$, which apparently satisfies $q_i((U^{k+1})^{[1]}) \subseteq U^{k+2}$ and which is of class $C^\infty$. Define $g\colon (U^{k+1})^{[1]} \to E$ via $g := \sum_{i=0}^{k+1} p_i \cdot (f^{<k+1>} \circ q_i)$. Since $k+1 \leq n-j$, the mapping $f^{<k+1>}$ is of class $C^j$, by induction. The functions $p_i$ and $q_i$ being smooth, we deduce from the Chain Rule and the Product Rule that $g$ is a mapping of class $C^j$ and thus continuous in particular. Now (20) shows that $f^{<k>}$ is of class $C^1$, with $(f^{<k>})^{[1]} = g$ of class $C^j$, and thus $f^{<k>}$ is of class $C^{j+1}$. $\square$

Summing up:

**Proposition 6.9** *Let $E \in \mathcal{M}$ be a topologized $\mathbb{K}$-vector space, $U \subseteq \mathbb{K}$ be an open, non-empty subset, and $f\colon U \to E$ be a map. Let $r \in \mathbb{N}_0 \cup \{\infty\}$. Then $f$ is a $C^r_{Sch}$-map if and only if $f$ is of class $C^r$.*

**Proof.** We may assume that $r \in \mathbb{N}_0$. If $f$ is a $C^r_{Sch}$-map, then $f = f^{<0>}$ is of class $C^r$ by Lemma 6.8. If, conversely, $f$ is of class $C^r$, the $f$ is a $C^r_{Sch}$-map by Proposition 6.2. $\square$

# 7 Main examples of $C^0$-concepts and associated differential calculi

In this section, we describe the main examples of $C^0$-concepts, and relate the resulting notions of $C^k$-maps to the literature. In our general definition of $C^0$-concepts, we considered rings and modules equipped with certain topologies, but did not presume that they are topological rings (resp., topological modules), to obtain a general framework which is as flexible as possible. Now, we turn to the special case where $\mathbb{K}$ and all modules are indeed topological. Of particular interest are the cases where $\mathbb{K}$ is the field of real or complex numbers (Subsection 7.1), or an ultrametric field (Subsection 7.2).

Here, we are using the following terminology: A *topological ring* is a commutative ring $\mathbb{K}$ with unit 1, equipped with a topology making the ring operations continuous; we require also that $\mathbb{K}^\times$ be open in $\mathbb{K}$, and that inversion $\iota\colon \mathbb{K}^\times \to \mathbb{K}$ be continuous. If, furthermore, $\mathbb{K}$ is a field, it is called a *topological field*. A *topological $\mathbb{K}$-module* is a module $E$ over a topological ring $\mathbb{K}$, equipped with a topology making the structure maps $E \times E \to E$ and $\mathbb{K} \times E \to E$ continuous. All topological rings, fields and modules considered in this section are assumed Hausdorff. It is clear that the class $\mathcal{M}$ of (Hausdorff) topological $\mathbb{K}$-modules satisfies Axioms I and II of a $C^0$-concept (taking the product topology on product spaces, and calling continuous mappings $C^0$). Since a continuous mapping from a topological space to a Hausdorff topological space is uniquely determined by its restriction to a dense subset, we readily obtain:

**Lemma 7.1** *Let $\mathbb{K}$ be a topological ring such that $\mathbb{K}^\times \subseteq \mathbb{K}$ is dense. Then Axiom* III *of a $C^0$-concept is satisfied for all $F \in \mathcal{M}$. In particular, this holds if $\mathbb{K}$ is a non-discrete topological field.* $\square$

Another important $C^0$-concept, based on rational maps, is discussed in Subsection 7.3; in this case, the modules are not topological. Also the smooth mappings of convenient differential calculus can be shown to arise from a suitable $C^0$-concept (Subsection 7.4).

## 7.1 Calculus on real or complex topological vector spaces

In this section, we have a closer look at the case where $\mathbb{K}$ is the topological field of real or complex numbers. The $C^0$-concept based on the class $\mathcal{M}$ of Hausdorff real (resp., complex) topological vector spaces (as just described) gives rise to a notion of $C^k$-maps (Definition 4.1). We show that a mapping to a *locally convex* topological vector space is a $C^k$-map in this sense if and only if it is a $C^k$-map in the classical sense of Michal-Bastiani. We also clarify the relation between $C^k$-maps in the complex case and complex analytic mappings.

Throughout this subsection, $\mathbb{K} \in \{\mathbb{R}, \mathbb{C}\}$.

**Definition 7.2** Let $E$ be a topological $\mathbb{K}$-vector space, $F$ a locally convex topological $\mathbb{K}$-vector space, and $k \in \mathbb{N}_0 \cup \{\infty\}$. A map

$$f\colon E \supseteq U \to F$$

on an open subset $U \subseteq E$ is called a $C^k_{MB}$-map (where MB stands for Michal-Bastiani) if it is continuous, the (real, resp., complex) derivatives

$$d^j f(x, v_1, \ldots, v_j) := \partial_{v_1} \cdots \partial_{v_j} f(x)$$

exist for all $j \in \mathbb{N}$ such that $j \leq k$, $x \in U$, $v_1, \ldots, v_j \in E$, and the mappings $d^j f\colon U \times E^j \to F$ so defined are continuous.

**Remark 7.3** In the literature, $C^k$-maps in the Michal-Bastiani sense (also known as Keller's $C^k_c$-maps [33]) are usually defined only if $\mathbb{K} = \mathbb{R}$ and if the domain $E$ is also locally convex. However, all of the basic theory (as described on a few pages in [22], Section 1) is easily seen to work just as well for non-locally convex $E$, and also over the complex field (by trivial adaptations of the proofs); all which matters is the local convexity of the range $F$. We can therefore use results from [22] in the following also if $\mathbb{K} = \mathbb{C}$, or if $E$ is not locally convex.

In particular, the Mean Value Theorem holds for $C^1_{MB}$-maps; $d^j f(x, \bullet) : E^j \to F$ is a symmetric $j$-linear map; compositions of composable $C^k_{MB}$-maps are $C^k_{MB}$; and furthermore a map $f : E \supseteq U \to F$ is of class $C^{k+1}_{MB}$ if and only if it is of class $C^1_{MB}$ and $df \colon E \times E \supseteq U \times E \to F$ is of class $C^k_{MB}$.

Our goal is to show that the concepts of $C^k$-maps and $C^k_{MB}$-maps are equivalent:

**Proposition 7.4** *Let $E$ be a topological $\mathbb{K}$-vector space, $F$ a locally convex topological $\mathbb{K}$-vector space, $k \in \mathbb{N}_0 \cup \{\infty\}$, and $f\colon U \to F$ be a map, where $U$ is an open subset of $E$. Then $f$ is of class $C^k$ if and only if $f$ is a $C^k_{MB}$-map.*

The proof depends on two technical lemmas, dealing with the differentiability of parameter-dependent integrals:

**Lemma 7.5** *Let $E$ be a topological $\mathbb{K}$-vector space, $F$ a locally convex topological $\mathbb{K}$-vector space, $U \subseteq E$ a non-empty open subset, $I \subseteq \mathbb{R}$ an open interval such that $[0,1] \subseteq I$, and $h\colon U \times I \to F$ be a mapping of class $C^k_{MB}$, where $k \in \mathbb{N}_0$. Suppose that the weak integral*

$$g(x) := \int_0^1 h(x,t)\, dt$$

*exists in $F$ for all $x \in U$, and suppose that the weak integrals*

$$\int_0^1 d_1^j h(x, t, y_1, \ldots, y_j)\, dt$$

*exist for all $j = 1, \ldots, k$ and $x \in U$, $y_1, \ldots, y_j \in E$, where $d_1^j h(x, t, y_1, \ldots, y_j) := d^j h((x,t), (y_1, 0), \ldots, (y_j, 0))$. Then $g\colon U \to F$ is a mapping of class $C^k_{MB}$, and*

$$d^j g(x, y_1, \ldots, y_j) = \int_0^1 d_1^j h(x, t, y_1, \ldots, y_j)\, dt$$

*for all $j = 1, \ldots, k$ and $x \in U$, $y_1, \ldots, y_j \in E$.*

**Proof.** The simple proof, based on induction on $j$ and standard estimates, is left to the reader. Full details can also be found in [29]. □

Conversely, we have:

**Lemma 7.6** *Let $E$ be a topological $\mathbb{K}$-vector space, $F$ a locally convex topological $\mathbb{K}$-vector space, $U \subseteq E$ a non-empty open subset, $I \subseteq \mathbb{R}$ an open interval such that $[0,1] \subseteq I$, and $h\colon U \times I \to F$ be a mapping of class $C^1_{MB}$. Suppose that the weak integral*

$$g(x) := \int_0^1 h(x,t)\, dt$$

*exists in $F$ for all $x \in U$, and defines a map $g\colon U \to F$ which is of class $C^1_{MB}$. Then the weak integral $\int_0^1 d_1 h(x,t,y)\, dt$ exists for all $x \in U$ and $y \in E$, and it is given by*

$$\int_0^1 d_1 h(x,t,y)\, dt = dg(x,y)\,.$$

**Proof.** We consider $g$ and $h$ as mappings into the completion $\widetilde{F}$ of $F$. The weak integral in question exists in $\widetilde{F}$. By the preceding lemma, it coincides with $dg(x,y)$ and therefore is an element of $F$. □

**Proof of Proposition 7.4.** We may assume that $k \in \mathbb{N}_0$. Clearly, $C^k$ implies $C^k_{MB}$ (using that $d^j f$ is continuous as a partial map of $f^{[j]}$, for each $j \leq k$). The converse direction is proved by induction on $k$. Assume that $f$ is $C^1_{MB}$ first. Then

$$g\colon U^{[1]} \to F, \quad (x,y,t) \mapsto \begin{cases} \frac{f(x+ty)-f(x)}{t} & \text{if } t \in \mathbb{R}^\times \\ df(x,y) & \text{if } t = 0 \end{cases}$$

is continuous. In fact, since $f$ is $C^1_{MB}$, so is the restriction $g|_{U^{]1[}}$ of $g$ to the open subset $U^{]1[} := \{(x,y,t) \in U^{[1]} : t \neq 0\}$ of $U^{[1]}$, and hence continuous. On the other hand, by the Fundamental Theorem of Calculus (cf. [22, Thm. 1.5]), we have

$$f(x+ty) - f(x) = \int_0^1 df(x+sty, ty)\, ds,$$

for all $(x,y,t)$ in the open neighbourhood $W := \{(x,y,t) \in U^{[1]} : x + [-1,2]ty \subseteq U\}$ of $U \times E \times \{0\}$ in $U^{[1]}$. Thus

$$g(x,y,t) = \int_0^1 df(x+sty, y)\, ds = \int_0^1 h((x,y,t), s)\, ds \qquad (21)$$

for all $(x,y,t) \in W$, where the integral depends continuously on $(x,y,t)$ by Lemma 7.5, applied with the $C^0$-map $h\colon W \times\ ]-1,2[\to F$, $h((x,y,t), s) := df(x+sty, y)$. Thus $g$ is continuous on $U^{[1]} = U^{]1[} \cup W$, and hence $f$ is $C^1$ with $f^{[1]} = g$.

Induction step: Suppose that $f$ is a $C^k_{MB}$-map now, where $k \geq 2$. We define $g$ and $h$ as above and note that $g|_{U^{]1[}}$ is a $C^k_{MB}$-map and hence of class $C^{k-1}$, by induction; furthermore, $h$ is a $C^{k-1}_{MB}$-map.

If $F$ is sequentially complete, then the existence of weak integrals is automatic, whence, in view of Eqn. (21), $g|_W$ is a $C^{k-1}_{MB}$-map by Lemma 7.5 and hence of class $C^{k-1}$, by induction. Now, $f$ being $C^1$ with $g = f^{[1]}$ of class $C^{k-1}$, the map $f$ is of class $C^k$, and we are home.

If $F$ fails to be sequentially complete, in order to be able to apply Lemma 7.5, it remains to show that the weak integrals

$$\int_0^1 d_1^j h((x,y,t), s; w_1, \ldots, w_j)\, ds \qquad (22)$$

exist for all $j = 1, \ldots, k-1$, $(x,y,t) \in U^{[1]}$, $w_1, \ldots, w_j \in E^{[1]}$. Now, the map $g|_{U^{]1[}}$ being of class $C^k_{MB}$, we inductively deduce from Lemma 7.6 that the weak integrals (22) exist, provided $(x,y,t) \in U^{]1[}$. It remains to prove the existence when $x \in U$, $y \in E$, $t = 0$. To enable this, we need to understand what the mappings $d_1^j h$ look like.

If $j = 1$, $(x,y,t) \in W$, $w_1 = (y_1, y_2, t_1) \in E^{[1]}$, and $s \in [0,1]$, using the linearity of $df(x+sty, \bullet)$ we obtain

$$\begin{aligned} d_1 h((x,y,t), s, (y_1, y_2, t_1)) &= d^2 f(x+sty, y, y_1) + std^2 f(x+sty, y, y_2) \\ &\quad + st_1 d^2 f(x+sty, y, y) + df(x+sty, y_2). \end{aligned}$$

Similarly, using the multilinearity of the higher differentials, a simple inductive argument shows that $d_1^j h((x,y,t), s, w_1, \ldots, w_j)$ (where $(x,y,t) \in W$ and $w_i = (y_{2i-1}, y_{2i}, t_i) \in E^{[1]}$ for $i = 1, \ldots, j$) is a sum of terms of the form

$$s^a t^b t_1^{c_1} \cdots t_j^{c_j} d^i f(x+sty, y_{\ell_1}, \ldots, y_{\ell_i}), \qquad (23)$$

where $i \in \{1, \ldots, j\}$, $\ell_1, \ell_2, \ldots \ell_i \in \{0, 1, \ldots, 2j\}$, $y_0 := y$ (for convenience of notation), and $a, b, c_1, \ldots, c_j \in \{0, 1, \ldots, j\}$. If now $t = 0$ in particular (this is the case we are interested in), fixing all variables except for $s$, note that each term (23) is a polynomial (actually, a monomial) in $s$ with coefficients in $F$. Hence $d_1^j h((x, y, 0), s, w_1, \ldots, w_j)$, as a function of $s$, is a polynomial with coefficients in $F$, say equal to $\sum_{n=0}^{j} s^n z_n$ for certain $z_n \in F$. Applying linear functionals, we now readily see that

$$\sum_{n=0}^{j} \left( \int_0^1 s^n ds \right) z_n = \sum_{n=1}^{j} \frac{1}{n+1} z_n \in F$$

satisfies the defining property of the weak integral $\int_0^1 h_1^j((x, y, 0), s, w_1, \ldots, w_j) \, ds$, which consequently exists. We can now complete the proof as in the sequentially complete case. $\square$

In order to distinguish the notions of $C^k$-maps over the real numbers and those over the complex numbers, let us write $C_{\mathbb{R}}^k$ and $C_{\mathbb{C}}^k$ for the moment. Then we have:

**Proposition 7.7** *Let $E$ and $F$ be complex topological vector spaces, where $F$ is locally convex, $U \subseteq E$ be an open subset, and $f : U \to F$ be a mapping. Then the following conditions are equivalent:*

(a) *$f$ is complex analytic in the sense of [10, Defn. 5.6], i.e., $f$ is continuous and, for every $x \in U$, there exists a zero-neighbourhood $V \subseteq E$ such that $x + V \subseteq U$ and certain continuous homogeneous polynomials $\beta_n : E \to F$ of degree $n$, such that $f$ admits the expansion: $f(x + y) = \sum_{n=0}^{\infty} \beta_n(y)$, for all $y \in V$.*

(b) *$f$ is of class $C_{\mathbb{C}}^{\infty}$.*

(c) *$f$ is of class $C_{\mathbb{R}}^{\infty}$, and $df(x, \bullet) : E \to F$ is complex linear for each $x \in U$.*

*If $F$ is Mackey complete, then (a)–(c) are also equivalent to any of the following:*

(d) *$f$ is of class $C_{\mathbb{C}}^1$.*

(e) *$f$ is of class $C_{\mathbb{R}}^1$, and $df(x, \bullet)$ is complex linear for each $x \in U$.*

**Proof.** (a)$\Rightarrow$(b): In view of the characterization of $C_{\mathbb{C}}^{\infty}$-maps given in Proposition 7.2, we can repeat the proof of [22, Prop. 2.4] verbatim to deduce that complex analytic maps are of class $C_{\mathbb{C}}^{\infty}$.

The implication (b)$\Rightarrow$(c) is obvious, and (c)$\Rightarrow$(a) is the content of [22, La. 2.5] (which remains valid for non-locally convex domains).

Now assume that $F$ is Mackey complete (see [36, La. 2.2 and Thm. 2.14]). Then clearly (b) implies (d), and (d) implies (e). To see that (e) implies (a), we argue as in [24, 1.4] (the local convexity of $E$ assumed there is inessential for the arguments). $\square$

For non-locally convex ranges, the situation is totally different. There are $C_{\mathbb{C}}^1$-maps from $\mathbb{C}$ to metrizable, complete non-locally convex spaces which are not $C_{\mathbb{C}}^2$, and there are $C_{\mathbb{C}}^{\infty}$-maps which are not given locally by their Taylor series, around any point. Furthermore, there are compactly supported, non-zero $C^\infty_\mathbb{C}$-maps from $\mathbb{C}$ to suitable non-locally convex spaces (showing that the Identity Theorem for Analytic Functions becomes false), and there are injective $C^\infty_\mathbb{R}$-maps on $\mathbb{R}$ whose derivative vanishes identically (see [28]).

## 7.2 Differential calculus over ultrametric fields

We now consider the $C^0$-concept described at the beginning of Section 7, in the case where $\mathbb{K}$ is a valued field, *i.e.*, a field, equipped with an absolute value $|.|\colon \mathbb{K} \to [0,\infty[$ which we assume non-trivial (meaning that it gives rise to a non-discrete topology on $\mathbb{K}$). Already having dealt with $\mathbb{R}$ and $\mathbb{C}$, we focus on valued fields $(\mathbb{K}, |.|)$ now whose absolute value $|.|$ satisfies the strong triangle inequality, viz. $|x+y| \leq \max\{|x|, |y|\}$ for all $x, y \in \mathbb{K}$; in this case, $(\mathbb{K}, |.|)$ is called an *ultrametric field*.

As a special case of our discussions in Section 6, we know that a function $f\colon \mathbb{K} \supseteq U \to \mathbb{K}$, where $\mathbb{K}$ is an ultrametric field, is $C^k$ in our sense if and only if it is $C^k$ in the usual sense of ultrametric calculus, as considered in Schikhof [58]. In the following, we inspect further important classes of mappings encountered in non-archimedian analysis (like strictly differentiable maps, or analytic maps), and relate them to our concept of $C^k$-maps.

**7.8** Recall that a topological vector space $E$ over an ultrametric field $\mathbb{K}$ is called *locally convex* if every zero-neighbourhood of $E$ contains an open $\mathbb{O}$-submodule of $E$, where $\mathbb{O} := \{t \in \mathbb{K} : |t| \leq 1\}$ is the valuation ring of $\mathbb{K}$. Equivalently, $E$ is locally convex if and only if its vector topology is defined by a family of ultrametric continuous seminorms $\gamma\colon E \to [0,\infty[$ on $E$ (cf. [46] for more information). Let $\mathbb{K}$ be a valued field. We call a topological $\mathbb{K}$-vector space *polynormed* if its vector topology is defined by a family of continuous seminorms (which need not be ultrametric when $\mathbb{K}$ is an ultrametric field). This terminology slightly deviates from the one in Bourbaki [15], where only polynormed spaces over ultrametric fields are considered whose topology arises from a family of continuous *ultrametric* seminorms, and which therefore are precisely the locally convex spaces over such fields in our terminology. We shall not presume that norms, nor Banach spaces, be ultrametric, unless saying so explicitly. For example, $\ell^1(\mathbb{Q}_p)$ is a non-ultrametric (and non-locally convex) Banach space over $\mathbb{Q}_p$. We shall also write $\|.\|_\gamma$ for a continuous seminorm $\gamma$.

**Definition 7.9** Let $\mathbb{K}$ be a valued field, $E$ be a normed $\mathbb{K}$-vector space, $F$ a polynormed $\mathbb{K}$-vector space, $U \subseteq E$ be open, and $f\colon U \to F$ be a map. Given $x \in U$, we say that $f$ is *strictly differentiable at $x$* if there exists a continuous linear map $A \in L(E,F)$ such that, for every $\varepsilon > 0$ and continuous seminorm $\gamma$ on $F$, there exists $\delta > 0$ such that

$$\|f(z) - f(y) - A.(z-y)\|_\gamma < \varepsilon \|z-y\|$$

for all $y, z \in U$ such that $\|z - x\| < \delta$ and $\|y - x\| < \delta$. Here $A$ is uniquely determined; we write $f'(x) := A$. The map $f$ is called *strictly differentiable* if it is strictly differentiable at each $x \in U$.

It can be shown that every strictly differentiable mapping $f$ is of class $C^1$ [26, La. 3.2]. This facilitates the following definition (cf. [26, Rem. 5.2]):

**Definition 7.10** Let $(\mathbb{K}, |.|)$ be a valued field, $E$ be a normed $\mathbb{K}$-vector space, $F$ be a polynormed $\mathbb{K}$-vector space, $U \subseteq E$ an open subset, and $k \in \mathbb{N}_0 \cup \{\infty\}$. A mapping $f: U \to F$ is called $k$ *times strictly differentiable* (or an $SC^k$-*map*, for short) if it is of class $C^k$ and $f^{[j]}: U^{[j]} \to F$ is strictly differentiable, for all $j \in \mathbb{N}_0$ such that $j < k$.

The class of $SC^k$-maps has many useful properties. For example, the Inverse Function Theorem and Implicit Function Theorem hold for $SC^k$-maps between open subsets of Banach spaces over a complete valued field $\mathbb{K}$ ([26]; cf. [15, 1.5.1] when $k = 1$).

We remark that, in the case of ultrametric fields, an Inverse Function Theorem cannot be based on the mere existence and continuity of differentials: the conclusion can fail spectacularly (see [58, Example 26.6], which refutes earlier claims to the contrary in [57]).

By definition, every $SC^k$-map is of class $C^k$. Conversely, the following can be shown (see [26], Remarks 5.3 and 5.4):

**Proposition 7.11** *Let $(\mathbb{K}, |.|)$ be a valued field, $E$ be a normed $\mathbb{K}$-vector space, $F$ be a polynormed $\mathbb{K}$-vector space, $U \subseteq E$ an open subset, and $k \in \mathbb{N}_0 \cup \{\infty\}$. Let $f: U \to F$ be a map. Then the following holds:*

(a) *If $f$ is of class $C^{k+1}$, then $f$ is an $SC^k$-map.*

(b) *If $\mathbb{K}$ is locally compact and $E$ is finite-dimensional, then $f$ is of class $C^k$ if and only if $f$ is an $SC^k$-map.* □

**Remark 7.12** In the real case, it is known that a mapping is $SC^1$ if and only if it is continuously Fréchet differentiable ([14, 2.3.3], cf. also [17, Thm. 3.8.1]), and it is known that every $k$ times continuously Fréchet differentiable mapping is $SC^k$ ([26], appendix).

In non-archimedean analysis, (locally) analytic functions are a well-established and widely used concept (see [15], also [58], [59] for the finite-dimensional case). We recall from [15] the definition of analytic functions, and some background material:

**7.13** Let $\mathbb{K}$ be a complete ultrametric field, $F$ be a locally convex topological $\mathbb{K}$-vector space, $n \in \mathbb{N}$, $E_1, \ldots, E_n$ be ultrametric normed $\mathbb{K}$-vector spaces, and $E := E_1 \times \cdots \times E_n$, equipped with the maximum norm. Given a multi-index $\alpha = (\alpha_i) \in \mathbb{N}_0^n$, we define

$$\alpha(j) := \sup\{i \in \mathbb{N}: 1 \leq i \leq n \text{ and } j > \alpha_1 + \cdots + \alpha_{i-1}\}$$

for $j \in \{1, \ldots, |\alpha|\}$, where $|\alpha| := \alpha_1 + \cdots + \alpha_n$. Thus $\alpha(1), \alpha(2), \ldots, \alpha(|\alpha|)$ is the finite sequence obtained by writing $\alpha_1$ times 1, then $\alpha_2$ times 2, etc.

**7.14** We let $E_\alpha := \prod_{j=1}^{|\alpha|} E_{\alpha(j)}$, equipped with the supremum norm; $\mathcal{L}_\alpha(E_1, \ldots, E_n; F)$ denotes the space of continuous $|\alpha|$-multilinear maps from $E_\alpha$ to $F$. We define

$$p_\alpha \colon E \to E_\alpha, \quad p_\alpha(x) := (p_{\alpha(j)}(x))_{j=1}^{|\alpha|},$$

where $p_i \colon E \to E_i$ is the canonical projection for $i = 1, \ldots, n$. A mapping $f \colon E \to F$ is called a *continuous multi-homogeneous polynomial of multi-degree $\alpha$* if there exists $u \in \mathcal{L}_\alpha(E_1, \ldots, E_n; F)$ such that $f = u \circ p_\alpha$; we let $P_\alpha(E_1, \ldots, E_n; F)$ be the space of all such $f$. If $n = 1$ and thus $E = E_1$, we write $\mathcal{L}_k(E; F)$ instead of $\mathcal{L}_{(k)}(E; F)$ and $P_k(E; F) := P_{(k)}(E; F)$. Then $\delta_k := p_{(k)} \colon E \to E_{(k)} = E^k$ is the k-fold diagonal map. Following Bourbaki, the elements of $P_k(E; F)$ will be called *continuous homogeneous polynomials of degree $k$* in the present section; however, elsewhere in this text we shall use the terminology from Appendix A.[2]

**7.15** Given $u \in \mathcal{L}_\alpha(E_1, \ldots, E_n; F)$ and a continuous seminorm $\gamma$ on $F$, we define[3]

$$\|u\|_\gamma := \inf\{M \geq 0 \colon (\forall (x_1, \ldots, x_k) \in E_\alpha) \; \gamma(u(x_1, \ldots, x_k)) \leq M \cdot \|x_1\| \cdot \ldots \cdot \|x_k\|\},$$

where $k := |\alpha|$. Given $f \in P_\alpha(E_1, \ldots, E_n; F)$, we define

$$\|f\|_\gamma := \inf\{\|u\|_\gamma \colon u \in \mathcal{L}_\alpha(E_1, \ldots, E_n; F) \text{ such that } u \circ p_\alpha = f\}.$$

**7.16** The elements $f$ of $\widehat{P}(E_1, \ldots, E_n; F) := \prod_{\alpha \in \mathbb{N}_0^n} P_\alpha(E_1, \ldots, E_n; F)$ are called *formal series*, and are written in the form $f = \sum_{\alpha \in \mathbb{N}_0^n} f_\alpha$, with $f_\alpha \in P_\alpha(E_1, \ldots, E_n; F)$. Given a formal series $f$, we let $I(f) \subset \,]0, \infty[^n$ be the set of all $R = (R_1, \ldots, R_n) \in \,]0, \infty[^n$ such that, for every continuous seminorm $\gamma$ on $F$, we have

$$\|f_\alpha\|_\gamma R^\alpha \to 0 \quad \text{as } |\alpha| \to \infty,$$

where $R^\alpha := R_1^{\alpha_1} \cdot \ldots \cdot R_n^{\alpha_n}$. The set $I(f)$ is called the *indicatrix of strict convergence of $f$*. If $n = 1$, we identify $\widehat{P}(E; F)$ with $\prod_{k \in \mathbb{N}_0} P_k(E; F)$.

---

[2]Note that, in Bourbaki's terminology, a continuous homogeneous polynomial $f$ of degree $k$ is not merely a homogeneous polynomial of degree $k$ which is continuous, but it is assumed that we can find a *continuous* k-multilinear map $u \colon E^k \to F$ such that $f = u \circ \delta_k$. In positive characteristic, where no polarization formula is available, this requirement might very well be stronger than the mere continuity of $f$.

[3]Our terminology slightly differs from Bourbaki's here, who uses on $\mathcal{L}_\alpha(E_1, \ldots, E_n; F)$ the equivalent seminorms defined via $\|u\|_\gamma' := \sup\{\gamma(u(x_1, \ldots, x_k)) \colon x_j \in E_{\alpha(j)}, \|x_j\| \leq 1\}$ (and the corresponding seminorms $\|.\|_\gamma'$ on $P_\alpha(E_1, \ldots, E_n; F)$). It is easy to see that $\|u\|_\gamma' \leq \|u\|_\gamma \leq a^k \|u\|_\gamma'$, where $a := \inf\{|t| \colon t \in \mathbb{K}^\times, |t| > 1\}$. The use of $\|.\|_\gamma'$ is problematic. For example, the inequality asserted in [15, A.4] is false if we take $E := F := G := \mathbb{Q}_p$ with $\|.\|_E := \|.\|_F := |.|_p$ the usual absolute value on $\mathbb{Q}_p$, $\|.\|_G := p^{-\frac{1}{2}} \cdot |.|_p$, $f \colon E \to G$, $x \mapsto x$ and $g \colon G \to F$, $x \mapsto x$. Then $\|f\|' = p^{-\frac{1}{2}}$, $\|g\|' = \|g \circ f\|' = 1 > \|g\|' \cdot \|f\|' = p^{-\frac{1}{2}}$ (cf. also [55, p. 59]).

**7.17** Given $R \in \,]0, \infty[^n$, let $B(R) := \{(y_1, \ldots, y_n) \in E \colon \|y_i\| \leq R_i \text{ for } i = 1, \ldots, n\,\}$ be the closed polyball of multi-radius $R$ around the origin; given $r \in \,]0, \infty[$, define $B(r) := \{y \in E \colon \|y\| \leq r\}$. Then $B(R)$ and $B(r)$ are open and closed $\mathbb{O}$-submodules of $E$, where $\mathbb{O} := \{t \in \mathbb{K} \colon |t| \leq 1\}$ is the valuation ring of $\mathbb{K}$.

**7.18** A function $f \colon U \to F$ on an open subset $U \subseteq E$ is called *analytic* if, for every $x \in U$, there exists a formal series $f_x = \sum_{\alpha \in \mathbb{N}_0^n} f_{x,\alpha} \in \widehat{P}(E_1, \ldots, E_n; F)$ and $R \in I(f_x)$ such that $x + B(R) \subseteq U$ and $(f_{x,\alpha}(y))_{\alpha \in \mathbb{N}_0^n}$ is a summable family in $F$ for all $y \in B(R)$, with limit

$$f(x+y) = \sum_{\alpha \in \mathbb{N}_0^n} f_{x,\alpha}(y).$$

The summability is automatic if $F$ is sequentially complete. The formal series $f_x$ is uniquely determined.

**7.19** We recall that a formal series in $\widehat{P}(E_1, \ldots, E_n; F)$ as before exists if and only if there exists a formal series $f_x = \sum_{k=0}^{\infty} f_{x,k} \in \widehat{P}(E; F)$ and $r \in I(f_x)$ such that $x + B(r) \subseteq U$ and $(f_{x,k}(y))_{k \in \mathbb{N}_0}$ is summable in $F$ for all $y \in B(r)$, with limit

$$f(x+y) = \sum_{k=0}^{\infty} f_{x,k}(y).$$

It is known that every analytic function between open subsets of a complete ultrametric field $\mathbb{K}$ is of class $C_{\mathbb{K}}^{\infty}$ [58, Cor. 29.11]. It is also known that every analytic mapping from an open subset of a normed $\mathbb{K}$-vector space to a locally convex topological $\mathbb{K}$-vector space is strictly differentiable with analytic derivatives ([15], 4.2.3 and 3.2.4). We show:

**Proposition 7.20** *Let $\mathbb{K}$ be a complete ultrametric field, $E$ an ultrametric normed $\mathbb{K}$-vector space, $F$ a locally convex topological $\mathbb{K}$-vector space, $U \subseteq E$ be open, and $f \colon U \to F$ a $\mathbb{K}$-analytic mapping. Then $f$ is of class $C_{\mathbb{K}}^{\infty}$, and $f^{[j]}$ is $\mathbb{K}$-analytic for each $j \in \mathbb{N}_0$.*

**Proof.** Apparently, we only need to show that every $\mathbb{K}$-analytic mapping is of class $C_{\mathbb{K}}^1$, and that $f^{[1]}$ is $\mathbb{K}$-analytic: the assertion then follows by a trivial induction.

It is well-known (and easy to see) that the directional derivative $df(x, y)$ exists for all $x \in U$ and $y \in E$ (cf. [15], 4.2.3 and 3.2.4). We can therefore define a mapping $g \colon U^{[1]} \to F$ (which will turn out to be $f^{[1]}$) via

$$g(x, y, t) := \begin{cases} \frac{1}{t}(f(x+ty) - f(x)) & \text{if } t \neq 0; \\ df(x, y) & \text{if } t = 0, \end{cases}$$

for all $(x, y, t) \in U^{[1]}$. The map $f$ and inversion $\mathbb{K}^{\times} \to \mathbb{K}$, $t \mapsto \frac{1}{t}$ being analytic, it readily follows that $g$ is analytic on the open set $U^{]1[} := \{(x, y, t) \in U^{[1]} \colon t \neq 0\}$ (cf. [15], 4.2.3, 3.2.7 and 3.2.8). In order that $g$ be analytic, it therefore only remains to show that $g$ is given by a convergent formal series around $(x_0, y_0, 0)$, for all $x_0 \in U$ and $y_0 \in E$. After a translation, we may assume without loss of generality that $0 \in U$ and $x_0 = 0$. There is $r \in \,]0,1[$ such that $B(r) \subseteq U$. After shrinking $r$, we may assume that

$$f(z) = \sum_{k=0}^{\infty} f_k(z) \quad \text{for all } z \in B(r),$$

where $f_k \in P_k(E;F)$ for each $k \in \mathbb{N}_0$, and

$$\lim_{k\to\infty} \|f_k\|_\gamma \, r^k = 0, \tag{24}$$

for every continuous ultrametric seminorm $\gamma$ on $F$. We set $R_1 := r$, $R_2 := \|y_0\| + 1$, $R_3 := \frac{r}{\|y_0\|+1}$, $R := (R_1, R_2, R_3)$. Then the polyball $B(R) \subseteq E \times E \times \mathbb{K}$ is contained in $B(r)^{[1]} \subseteq U^{[1]}$, and $(0, y_0, 0) \in B(R)$. Fix $\gamma$. For each $k \in \mathbb{N}_0$, we find $u_k \in \mathcal{L}_k(E;F)$ such that $\|u_k\|_\gamma \leq \|f_k\|_\gamma + 1$ and $u_k \circ \delta_k = f_k$. Write $z_1 := y_0$. Then

$$g(z_0, z_1 + z_2, t) = \tfrac{1}{t}(f(z_0 + tz_1 + tz_2) - f(z_0)) = \sum_{k=1}^{\infty} \tfrac{1}{t}\left(f_k(z_0 + tz_1 + tz_2) - f_k(z_0)\right)$$

for all $(z_0, z_2, t) \in B(R)$ such that $t \neq 0$, where

$$\begin{aligned}
&\tfrac{1}{t}\left(f_k(z_0 + tz_1 + tz_2) - f_k(z_0)\right) \\
&= \tfrac{1}{t}\left(u_k(z_0 + tz_1 + tz_2, \ldots, z_0 + tz_1 + tz_2) - u_k(z_0, \ldots, z_0)\right) \\
&= \sum_{i \in I_k} t^{\#\{j:\, i_j \neq 0\}-1} u_k(z_{i_1}, \ldots, z_{i_k}) = \sum_{a=1}^{k}\sum_{b=0}^{a} g_{(k-a,b,a-1)}(z_0, z_2, t),
\end{aligned}$$

using the set $I_k := \{i = (i_1, \ldots, i_k) \in \{0,1,2\}^k\} \setminus \{(0,0,\ldots,0)\}$ and the continuous multi-homogeneous polynomial $g_{(k-a,b,a-1)} \in P_{(k-a,b,a-1)}(E, E, \mathbb{K}; F)$ on the product $E \times E \times \mathbb{K}$, of multi-degree $(k-a, b, a-1)$, defined via

$$g_{(k-a,b,a-1)}(z_0, z_2, t) := \sum t^{a-1} u_k(z_{i_1}, \ldots, z_{i_k}) \tag{25}$$

for $(z_0, z_2, t) \in E \times E \times \mathbb{K}$, where the summation is over all $i \in I_k$ such that $\#\{j:\, i_j \neq 0\} = a$ and $\#\{j:\, i_j = 2\} = b$. It is apparent from (25) that $\|g_{(k-a,b,a-1)}\|_\gamma \leq \|u_k\|_\gamma$, whence

$$\|g_{(k-a,b,a-1)}\|_\gamma \leq \|f_k\|_\gamma + 1. \tag{26}$$

For multi-indices $\alpha \in \mathbb{N}_0^3$ which are not of the form $(k-a, b, a-1)$ with $k \in \mathbb{N}$, $1 \leq a \leq k$ and $0 \leq b \leq a$, we define $g_\alpha := 0 \in P_\alpha(E, E, \mathbb{K}; F)$. Note that, if $g_\alpha \neq 0$, then $\alpha = (k-a, b, a-1)$ for uniquely determined $k \in \mathbb{N}$, $a \in \{1, \ldots, k\}$ and $b \in \{0, \ldots, a\}$. Then, in view of (24),

$$\|g_\alpha\|_\gamma R^\alpha = R_1^{k-a} R_2^b R_3^{a-1} \|g_\alpha\|_\gamma = \tfrac{(\|y_0\|+1)^{b-a+1}}{r} \|g_\alpha\|_\gamma r^k \leq \tfrac{\|y_0\|+1}{r}(\|f_k\|_\gamma + 1) r^k \to 0 \tag{27}$$

as $|\alpha| = k + b - 1 \to \infty$, using that $|\alpha| = k + b - 1 \leq 2k - 1$ and thus $k \geq \frac{|\alpha|+1}{2} \to \infty$.

Note that $g_{(k-a,b,a-1)}$ is independent of the choice of $\gamma$ and the $u_k$'s: if $\gamma'$ and maps $u'_k \in \mathcal{L}_k(E;F)$ are given with $u'_k \circ \delta_k = f_k$ and $\|u'_k\|_{\gamma'} \leq \|f_k\|_{\gamma'} + 1$, leading to multi-homogeneous polynomials $g'_{(k-a,b,a-1)}$, then

$$\sum_{a=1}^{k}\sum_{b=0}^{a} t^{a-1} s^b g_{(k-a,b,a-1)}(z_0, z_2, 1) = \tfrac{1}{t}\left(f_k(z_0 + tz_1 + tsz_2) - f_k(z_0)\right)$$

$$= \sum_{a=1}^{k}\sum_{b=0}^{a} t^{a-1} s^b g'_{(k-a,b,a-1)}(z_0, z_2, 1)$$

for all $(z_0, z_2, t) \in B(R)$ and $s \in \mathbb{K}$ with $t \neq 0$ and $|s| \leq 1$, whence $g_{(k-a,b,a-1)}(z_0, z_2, 1) = g'_{(k-a,b,a-1)}(z_0, z_2, 1)$ clearly and thus $g_{(k-a,b,a-1)} = g'_{(k-a,b,a-1)}$, using the multi-homogeneity. Now $g_{(k-a,b,a-1)}$ being independent of the choice of $\gamma$ and the $u_k$'s, we see that (27) holds for all continuous ultrametric seminorms $\gamma$ on $F$, whence $R$ is contained in the indicatrix of strict convergence of $\sum_\alpha g_\alpha$. This entails that the family $(g_\alpha)_\alpha$ is uniformly summable on $B(R)$, to an analytic (and hence continuous) function $\widehat{g}\colon B(R) \to \widetilde{F}$ into the completion $\widetilde{F}$ of $F$ (cf. [15, 4.2.4]). For all $(z_0, z_2, t) \in B(R)$ such that $t \neq 0$, we have

$$g(z_0, z_1 + z_2, t) = \sum_{k=1}^{\infty}\sum_{a=1}^{k}\sum_{b=0}^{a} g_{(k-a,b,a-1)}(z_0, z_2, t) = \sum_{\alpha \in \mathbb{N}_0^3} g_\alpha(z_0, z_2, t) = \widehat{g}(z_0, z_2, t)$$

(the re-ordering is clearly permissable). Letting $t \to 0$ for fixed $z_0$ and $z_2$, the left hand side converges to $df(z_0, z_1+z_2) = g(z_0, z_1+z_2, 0)$, the right hand side to $\widehat{g}(z_0, z_2, 0)$, by continuity of $\widehat{g}$. Thus $g((0, y_0, 0) + (z_0, z_2, t)) = g(z_0, z_1 + z_2, t) = \widehat{g}(z_0, z_2, t) = \sum_\alpha g_\alpha(z_0, z_2, t)$ for all $(z_0, z_2, t) \in B(R)$ is the desired expansion of $g$ around $(0, y_0, 0)$. $\square$

## 7.3 Rational maps

Assume $\mathbb{K}$ is an infinite field and consider the class $\mathcal{M}$ of finite-dimensional vector spaces over $\mathbb{K}$ together with their Zariski-topologies. Let $C^0$-maps be "rational" maps, *i.e.*, maps of the form $f(x) = \frac{p(x)}{q(x)}$ with a vector-valued polynomial $p$ and a scalar-valued polynomial $q$. Then Axioms I and II of a $C^0$-concept are clearly satisfied, and also the Determination Axiom holds since $\mathbb{K}$ is assumed to be infinite, so non-empty Zariski-open sets contain infinitely many points and rational functions are uniquely determined by their values at these points even if we remove one of them. By definition all polynomials are $C^0$. By the product rule, they are then $C^\infty$, and the quotient and chain rule now imply that all rational maps are $C^\infty$. Therefore, in this case, the classes $C^0, C^1, \ldots, C^\infty$ all coincide, and we see that the importance of our construction is not so much to produce new classes of maps but to provide a calculus for a given class. Here we recover (for infinite fields) the well-known differential calculus of rational maps (cf. e.g. [34] and the appendix to [41]).

**Remark 7.21** If $\mathbb{K}$ is algebraically closed and of characteristic 0, then each algebraic $\mathbb{K}$-group is a "rational" Lie group. It would be interesting to know whether this remains true for more general fields of characterictic 0.

## 7.4 Further $C^0$-concepts

It is possible, in principle, to describe the convenient differential calculus of Frölicher, Kriegl and Michor ([21], [36]) by means of a $C^0$-concept, as follows. Define $\mathbb{K} = \mathbb{R}$, let $\mathcal{M}$ be the class of all topologized real vector spaces $(E, \tau_E)$ such that $\tau_E$ is the $c^\infty$-refinement of some Mackey complete locally convex vector topology $\lambda_E$ on $E$, i.e., $\tau_E$ is the final topology on $E$ with respect to the set of all smooth curves $\mathbb{R} \to (E, \lambda_E)$, where smoothness is understood in the usual sense (cf. [36, Defn. 2.12]). Given $(E, \tau_E), (F, \tau_F) \in \mathcal{M}$ and open sets $U \in \tau_E$, $V \in \tau_F$, let $C^0(U, V) := C(U, V)$ be the set of all continuous mappings $U \to V$ with respect to the topologies induced by $\tau_E$ and $\tau_F$. We equip $E \times F$ with the $c^\infty$-refinement of the product topology on $(E, \lambda_E) \times (F, \lambda_F)$; the resulting topology is independent of the choices of $\lambda_E$ and $\lambda_F$, as it is easily seen to be final with respect to the set of all curves $\gamma = (\gamma_1, \gamma_2) \colon \mathbb{R} \to E \times F$ whose coordinate functions $\gamma_1 \colon \mathbb{R} \to (E, \lambda_E)$, $\gamma_2 \colon \mathbb{R} \to (F, \lambda_F)$ are smooth; the smoothness property of the curves $\gamma_1$ and $\gamma_2$ however only depends on $\tau_E$ and $\tau_F$. It is easy to check that all axioms of a $C^0$-concept are satisfied. Furthermore, a map between open subsets of spaces $(E, \tau_E), (F, \tau_F)$ in $\mathcal{M}$ is smooth in the sense provided by this $C^0$-concept if and only if it is a smooth map between the corresponding $c^\infty$-open subsets of $(E, \lambda_E)$ and $(F, \lambda_F)$, in the sense of convenient differential calculus. We suppress the details.

In spite of this possible embedding of convenient differential calculus in our general framework, which may be of theoretical interest, the original formulation of convenient differential calculus is of course preferable for practical purposes.

# Part II: Manifolds and Lie Groups

# 8 Manifolds and bundles

## 8.1 Manifolds modelled on a $C^0$-concept

We fix a base ring $\mathbb{K}$, a $C^0$-concept over $\mathbb{K}$, and a $\mathbb{K}$-module $E \in \mathcal{M}$, called the "model space". A $C^k$-*manifold with atlas* (*modelled on* $E$) (where $k \in \mathbb{N}_0 \cup \{\infty\}$) is a topological space $M$ together with an $E$-*atlas* $\mathcal{A} = \{(\phi_i, U_i) : i \in I\}$. This means that $U_i$, $i \in I$, is a covering of $M$ by open sets, and $\phi_i : M \supseteq U_i \to \phi_i(U_i) \subseteq E$ is a *chart*, i.e. a homeomorphism of the open set $U_i \subseteq M$ onto an open subset $\phi_i(U_i) \subseteq E$, and any two charts $(\phi_i, U_i), (\phi_j, U_j)$ are $C^k$-*compatible* in the sense that

$$\phi_{ij} := \phi_i \circ \phi_j^{-1}|_{\phi_j(U_i \cap U_j)} \colon \phi_j(U_i \cap U_j) \to \phi_i(U_i \cap U_j)$$

and its inverse $\phi_{ji}$ are of class $C^k$. We see no reason to assume that the topology of $M$ is separated (compare [37, p. 23] for this issue); if we want to work with separated manifolds, then we also have to assume that all $\mathbb{K}$-modules belonging to our $C^0$-concept are separated (an assumption which we have not needed so far!).

If the atlas $\mathcal{A}$ is *maximal* in the sense that it contains all compatible charts, then $M$ is called a $C^k$-*manifold* (*modelled on E*).

Let $(M, \mathcal{A}), (N, \mathcal{B})$ be $C^k$-manifolds with atlas modelled on $\mathbb{K}$-modules belonging to a given $C^0$-concept. A map $f\colon M \to N$ is *of class $C^k$* if it is continuous and if, for all choices of charts $(\phi, U)$ of $M$ and $(\psi, W)$ of $N$,

$$\psi \circ f \circ \phi^{-1}\colon \phi(U \cap f^{-1}(W)) \to \psi(W)$$

is of class $C^k$. Then one shows as usual that $C^k$-manifolds with atlas modelled on a given $C^0$-concept form a category. Taking maximal atlasses, we see that $C^k$ manifolds form a category.

If $\mathbb{K}$ is a field, then Lemma 4.9 shows that smoothness of a map may be checked by using any sub-atlas of the given one. If $\mathbb{K}$ is not a field, this may fail. For this reason we distinguish between the categories of manifolds and of manifolds with atlas. For instance, product atlasses and bundle atlasses are in general not maximal, and in case of a base ring statements may become false when using maximal atlasses on products and on bundles.

In particular, all $\mathbb{K}$-modules from the class $\mathcal{M}$ are $C^\infty$-manifolds and we can define *smooth functions on $M$* to be smooth maps $f\colon M \to \mathbb{K}$. The space $C^\infty(M)$ of smooth functions on $M$ may be reduced to the constants, and it may also happen that $C^\infty(U_i)$ is reduced to the constants for all $i$ (e.g. case of topological vector spaces that admit no non-zero continuous linear forms). Therefore it is no longer possible to define differential geometric objects via their action on smooth functions. In the category of smooth manifolds with atlas one can form direct products: given two $C^k$-manifolds $(M, \mathcal{A}), (N, \mathcal{B})$, endow $M \times N$ with the topology generated by the $U_\phi \times W_\psi$, $\phi \in \mathcal{A}, \psi \in \mathcal{B}$, together with the open subsets of the $U_\phi \times W_\psi$ (this is in general finer than the product topology and coincides with the product topology if our $C^0$-concept is "topological"), and the charts are given by the maps $\phi \times \psi$. These charts are again $C^k$-compatible and define an atlas $\mathcal{A} \times \mathcal{B}$. Clearly direct products of smooth maps are then smooth maps. If $\mathbb{K}$ is a field, this remains true also if we complete $\mathcal{A} \times \mathcal{B}$ to a maximal atlas.

It is also possible to define *submanifolds*; details will be given in [8].

## 8.2 Tangent bundle

A point $p \in M$ is described in the form $\phi_i^{-1}(x)$ with $x \in \phi_i(U_i)$ and $i \in I$. In a different chart it is given by $p = \phi_j^{-1}(y)$. In other words, $M$ is the set of equivalence classes $S/\sim$, where

$$S := \{(i, x) \mid x \in \phi_i(U_i)\} \subseteq I \times E,$$

and $(i, x) \sim (j, y)$ if and only if $\phi_i^{-1}(x) = \phi_j^{-1}(y)$. We write $p = [i, x] \in M = S/\sim$.

Next we define an equivalence relation on the set

$$TS := S \times E \subseteq I \times E \times E$$

via $(i, x, v) \sim (j, y, w)$ if and only if $\phi_j \circ \phi_i^{-1}(x) = y$, $d(\phi_j \circ \phi_i^{-1})(x)v = w$. By the chain rule, this is an equivalence relation. We denote equivalence classes by $[i, x, v]$, and let

$$TM := TS/\sim.$$

If $[i, x, v] = [j, y, w]$, then $[i, x] = [j, y]$, and hence the map

$$\pi : TM \to M, \quad [i, x, v] \mapsto [i, x]$$

is well-defined. For $p = [i, x] \in M$, we let

$$T_p M := \pi^{-1}(p) = \{[i, x, v] \in TM \mid v \in E\}.$$

The map

$$T_x \phi_i^{-1} : E \to T_p M, \quad v \mapsto [i, x, v]$$

is a bijection (surjective by definition, injective since the differentials are bijections), and we can use it to define a $\mathbb{K}$-module structure on $T_p M$ which actually does not depend on $(i, x)$. We then use the linear bijection $T_x \phi_i^{-1}$ to transport the topology and $C^0$-manifold strucure from $E$ to $T_p M$, which enables us to speak of $C^0$-maps on $T_p M$ and its open subsets (although $T_p M$ might not be a member of the given class $\mathcal{M}$). The $\mathbb{K}$-module $T_p M$ is called the *tangent space of $M$ at $p$*.

We define an atlas $T\mathcal{A} := (T\phi_i)_{i \in I}$ on $TM$ by

$$TU_i := \pi^{-1}(U_i), \quad T\phi_i : TU_i \to E \times E, \quad [i, x, v] \mapsto (\phi_i(x), v)$$

and give $TM$ the topology generated by the $TU_i$ and their open subsets (note the topology on $TU_i$ may be strictly finer than the product topology). Change of charts is now given by

$$T\phi_{ij} : (x, v) \mapsto (\phi_{ij}(x), d\phi_{ij}(x)v)$$

which is $C^{k-1}$ if $f$ is $C^k$. Thus $(TM, T\mathcal{A})$ is a manifold with atlas.

If $f: M \to N$ is $C^k$ we define its *tangent map* by

$$Tf : TM \to TN, \quad [i, x, v] \mapsto [j, f_{ij}(x), df_{ij}(x)v]$$

where $f_{ij} = \psi_j \circ f \circ \phi_i^{-1}$ (supposed to be defined on a non-empty open set). In other words,

$$Tf = (T\psi_j)^{-1} \circ (f_{ij}, df_{ij}) \circ T\phi_i = (T\psi_j)^{-1} \circ Tf_{ij} \circ T\phi_i.$$

This is well-defined, linear in fibers and $C^{k-1}$. Clearly the functorial rules hold, i.e. we have defined a covariant functor $T$ from the category of $C^k$-manifolds modelled on a $C^0$-concept into the category of $C^{k-1}$-manifolds. (In [8] it will be shown that this functor can be seen as the functor of "scalar extension by dual numbers over $\mathbb{K}$.")

If $f: M \to \mathbb{K}$ is a smooth function, then $T_x f: T_x M \to T_x \mathbb{K} = \mathbb{K}$ gives rise to a function $TM \to \mathbb{K}$, linear in fibers, which we denote, by some abuse of notation, again by $Tf$ or by $df$. The product rule **3.3** implies that $d(fg) = f dg + g df$.

## 8.3 General fiber bundles

General bundles over $M$ are defined following the same pattern as above for the tangent bundle: assume $M$ is modelled on $E$ and let $N$ be some manifold modelled on $F \in \mathcal{M}$. Assume that for all triples $(i, j, x)$ such that $(i, x), (j, x) \in S$ (i.e. $x \in \phi_i(U_i) \cap \phi_j(U_j)$) a diffeomorphism $g_{ij}(x) : N \to N$ is given such that the cocycle relations

$$g_{ij}(\phi_{jk}(x))g_{jk}(x) = g_{ik}(x), \quad g_{ii}(x) = \mathrm{id}_N,$$

are satisfied and such that $(x, w) \mapsto g_{ij}(x)w$ is smooth wherever defined. Then we define an equivalence relation on $S \times N$ by

$$(i, x, v) \sim (j, y, w) \quad \text{if and only if} \quad \phi_{ij}(x) = y, \, g_{ij}(x)v = w.$$

By the cocycle relations, this is indeed an equivalence relation, and by the smoothness assumption, $B := S \times F/\sim$ can be turned into a manifold with atlas modelled on $E \times F$ and locally isomorphic to $U_i \times N$ and such that the projection $\pi \colon B \to M$, $[i, x, w] \mapsto [i, x]$ is a well-defined smooth map whose fibers are all diffeomorphic to $N$.

If $N$ carries an additional structure ($\mathbb{K}$-module, affine space, projective space, ...) and the $g_{ij}$ respect this structure, then each fiber also carries this structure, and homomorphisms are required to respect this structure. In particular, if $N$ is a $\mathbb{K}$-module, then $b$ is called a *vector bundle*. *Direct sums* of vector bundles are defined as usual; it is also possible to define *tensor products* (cf. [8]) in a suitable way, but dual bundles should be strictly avoided. Nevertheless, *differential forms* can be defined, but should not be seen as sections of vector bundles (cf. [5], [51]).

## 8.4 Vector fields and derivations

A *section* of a vector bundle $B$ over $M$ is a smooth map $\xi \colon M \to B$ such that $\pi \circ \xi = \mathrm{id}_M$. The sections of $B$ form a module $\Gamma^\infty(B)$ over $C^\infty(M)$. Sections of $TM$ are also called *vector fields* and are often denoted by $X, Y, Z, \ldots$, and we use the notation $\mathfrak{X}(M)$ for $\Gamma^\infty(TM)$. In a chart $(U_i, \phi_i)$, vector fields can be identified with smooth maps $X_i \colon E \supseteq \phi_i(U_i) \to E$, given by

$$X_i := \mathrm{pr}_2 \circ T\phi_i \circ X \circ \phi_i^{-1} \colon \phi_i^{-1}(U_i) \to U_i \to TU_i \cong U_i \times E \to E.$$

Similarly, sections of an arbitrary vector bundle are locally represented by smooth maps $\xi_i \colon \phi_i^{-1}(U_i) \to F$.

For a vector field $X : M \to TM$ and a smooth function $f : M \to \mathbb{K}$, recall that $df \colon TM \to \mathbb{K}$ is smooth and hence we can define a smooth function $L_X f$ by $L_X f := df \circ X$. Then we have the Leibniz rule: $L_X(fg) = d(fg) \circ X = (fdg + gdf) \circ X = gL_X f + fL_X g$. Thus $X \mapsto L_X$ is a $\mathbb{K}$-linear map into the space of derivations of $C^\infty(M) := C^\infty(M, \mathbb{K})$.

**Remark 8.1** The map

$$\mathfrak{X}(M) \to \mathrm{Der}(C^\infty(M)), \quad X \mapsto L_X \tag{28}$$

will in general neither be injective nor surjective, not even locally; it therefore cannot be used to define the Lie algebra structure on $\mathfrak{X}(M)$.

**Example 8.2** Given $p \in \,]0,1[$, the map (28) fails to be injective if $M$ is a non-empty open subset of the real topological vector space $L^p(\mathbb{R})$, which does not have non-zero continuous linear functionals; in this case $C^\infty(M)$ consists of locally constant functions only, whence $L_X = 0$ for all $X \in \mathfrak{X}(M)$.

**Example 8.3** To see that the map (28) need not be surjective, we vary [36, La. 28.4]. Let $M$ be a non-empty open subset of an infinite-dimensional real Hilbert space $E$. Then the closure of $L(E,\mathbb{R}) \vee L(E,\mathbb{R})$ in the Banach space $L^2_{\text{sym}}(E,\mathbb{R})_b$ of continuous, symmetric bilinear forms on $E$ is the proper subspace $K_{\text{sym}}$ of compact symmetric bilinear forms (cf. [56], Theorem III.9.5 and III.9.2, Corollary 1). By the Hahn-Banach theorem, there is a non-zero continuous linear functional $\lambda \in L^2_{\text{sym}}(E,\mathbb{R})'_b$ vanishing on $K_{\text{sym}}$. Given $f \in C^\infty(M)$, the map $d^2 f \colon M \to L^2_{\text{sym}}(E,\mathbb{R})_b$ is smooth, using that mappings between Banach spaces are smooth in the Michal-Bastiani sense if and only if they are smooth in the Fréchet sense [33]. Hence

$$D \colon C^\infty(M) \to C^\infty(M), \quad Df := \lambda \circ d^2 f$$

maps into $C^\infty(M)$ indeed. The linear map $D$ is a derivation since $d^2(fg)(x)(v,w) = f(x)d^2 g(x)(v,w) + g(x)d^2 f(x)(v,w) + df(x)(v)dg(x)(w) + df(x)(w)dg(x)(v)$ for $x \in M$, $v,w \in E$, where the sum of the final two terms corresponds to the element $2\, df(x) \vee dg(x)$ of $L(E,\mathbb{R}) \vee L(E,\mathbb{R}) \subseteq K_{\text{sym}}$, and thus $D(fg)(x) = \lambda(d^2(fg)(x)) = f(x)\lambda(d^2 g(x)) + g(x)\lambda(d^2 f(x)) = f(x)Dg(x) + g(x)Df(x)$. There exists $\beta \in L^2_{\text{sym}}(E,\mathbb{R})$ such that $\lambda(\beta) \neq 0$. Define $f_\beta(x) := \frac{1}{2}\beta(x,x)$ for $x \in M$. Then $Df_\beta(x) = \lambda(\beta) \neq 0$ for all $x \in M$. Hence $D \neq 0$. Given $X \in \mathfrak{X}(M)$, we either have $X = 0$ (whence $L_X \neq D$), or there exists $x_0 \in M$ such that $X(x_0) \neq 0$; we may assume that $x_0 \neq 0$. Define $\alpha(v,w) := \langle v, x_0 \rangle \langle w, X(x_0) \rangle + \langle w, x_0 \rangle \langle v, X(x_0) \rangle$ for $v,w \in E$. Then $\alpha \in K_{\text{sym}}$ and thus $Df_\alpha = \lambda(\alpha) = 0$, but $(L_X f_\alpha)(x_0) = \|X(x_0)\|^2 \cdot \|x_0\|^2 + \langle x_0, X(x_0) \rangle^2 > 0$ and thus $L_X \neq D$ also in this case.

**Theorem 8.4** *There is a unique structure of a Lie algebra over $\mathbb{K}$ on $\mathfrak{X}(M)$ such that for all $X,Y \in \mathfrak{X}(M)$ and $(i,x) \in S$,*

$$[X,Y]_i(x) = dY_i(x)X_i(x) - dX_i(x)Y_i(x). \tag{29}$$

**Proof.** The uniqueness is clear. Let us show that, on the intersection of two chart domains, the bracket $[X,Y]$ is independent of the choice of chart. To this end, assume that $(i,x) \sim (j,y)$, i.e. $y = \phi_j \phi_i^{-1} x = \phi_{ji}(x)$; then $X_j(y) = d\phi_{ji}(x) X_i(x)$ or, considering now $d\phi_{ij}$ as a function of two arguments,

$$X_j \circ \phi_{ji} = d\phi_{ji} \circ (\text{id}, X_i). \tag{30}$$

We have to show that $[X,Y]_i$, defined by (29), has the same transformation property under changes of charts. We abbreviate $\phi := \phi_{ji}$. By (30), we have $Y_j = d\phi \circ (\mathrm{id}, Y_i) \circ \phi^{-1}$ and thus $dY_j = d(d\phi \circ (\mathrm{id}, Y_i) \circ \phi^{-1})$. Using this formula and (30), we get

$$\begin{aligned} dY_j(\phi(x))X_j(\phi(x)) &= d(d\phi \circ (\mathrm{id}, Y_i) \circ \phi^{-1})(T\phi(x, X_i(x))) \\ &= d(d\phi \circ (\mathrm{id}, Y_i))(x)X_i(x) \\ &= d^2\phi(x)(X_i(x), Y_i(x)) + d\phi(x)dY_i(x)X_i(x). \end{aligned}$$

We exchange the roles of $X$ and $Y$ and take the difference of the two equations thus obtained: we get, using Schwarz' lemma (Lemma 4.6),

$$dY_j(\phi(x))X_j(\phi(x)) - dX_j(\phi(x))Y_j(\phi(x)) = d\phi(x)(dY_i(x)X_i(x) - dX_i(x)Y_i(x))$$

which had to be shown. Summing up, the bracket operation $\mathfrak{X}(M) \times \mathfrak{X}(M) \to \mathfrak{X}(M)$ is well-defined, and it clearly is $\mathbb{K}$-bilinear and satisfies the identity $[X,X] = 0$.

All that remains to be proved is the Jacobi identity. This is done by a direct computation which involves only the chain rule and Schwarz' lemma: define a (chart dependent) "product" of $X_i$ and $Y_i$ by

$$(X_i \cdot Y_i)(x) := dY_i(x)X_i(x).$$

Then, by a direct calculation, one shows that this product is a *left symmetric* or *Vinberg algebra* (cf. [34]):

$$X_i \cdot (Y_i \cdot Z_i) - (X_i \cdot Y_i) \cdot Z_i = Y_i \cdot (X_i \cdot Z_i) - (Y_i \cdot X_i) \cdot Z_i.$$

But it is immediately verified that for every left symmetric algebra, the commutator $[X_i, Y_i] = X_i \cdot Y_i - Y_i \cdot X_i$ satisfies the Jacobi identity. □

The Lie bracket is natural in the following sense: assume $\phi\colon M \to N$ is a smooth map and $X \in \mathfrak{X}(M)$, $Y \in \mathfrak{X}(N)$. We say that the pair $(X,Y)$ is $\phi$-*related* if

$$Y \circ \phi = T\phi \circ X.$$

**Lemma 8.5** *If $(X,Y)$ and $(X',Y')$ are $\phi$-related, then so is $([X,X'],[Y,Y'])$. In particular, the diffeomorphism group of $M$ acts by automorphisms on the Lie algebra $\mathfrak{X}(M)$.*

**Proof.** This is the same calculation as the one after Eqn. (30). □

It easily follows from the definitions that the map $\mathfrak{X}(M) \to \mathrm{Der}(C^\infty(M))$, $X \mapsto L_X$ is a homomorphism of Lie algebras.

# 9 Lie groups

## 9.1 Manifolds with multiplication

As before, we fix a $C^0$-concept over $\mathbb{K}$ and consider smooth manifolds $M$ (in general with atlas, if $\mathbb{K}$ is a ring, and with maximal atlas if $\mathbb{K}$ is a field). A *product* or *multiplication* *map* on $M$ is a smooth binary map $m\colon M \times M \to M$, and *homomorphisms of manifolds with multiplication* are smooth maps that are compatible with the respective multiplication maps. (Note that the $C^0$-concept enters explicitly via the manifold structure on the product $M \times M$.) *Left and right multiplication operators*, defined by $l_x(y) = m(x,y) = r_y(x)$, are partial maps of $m$ and hence smooth self-maps of $M$. Applying the tangent functor to this situation, we see that $(TM, Tm)$ is again a manifold with multiplication, and tangent maps of homomorphisms are homomorphisms of the respective tangent spaces. The tangent map $Tm$ is given by the formula

$$T_{(x,y)}m(\delta_x, \delta_y) = T_{(x,y)}m((\delta_x, 0_y) + (0_x, \delta_y)) = T_x(r_y)\delta_x + T_y(l_x)\delta_y. \tag{31}$$

Formula (31) is nothing but the rule on partial derivatives (**3.8**), written in the language of manifolds. In particular, (31) shows that the canonical projection and the zero section,

$$\pi\colon TM \to M, \quad \delta_p \mapsto p, \qquad z\colon M \to TM, \quad p \mapsto 0_p \tag{32}$$

are homomorphisms of manifolds with multiplication. We will always identify $M$ with the subspace $z(M)$ of $TM$. Then (31) implies that the operator of left multiplication by $p = 0_p$ in $TM$ is nothing but $T(l_p)\colon TM \to TM$, and similarly for right multiplications.

## 9.2 Lie groups

A *Lie group* is a smooth manifold $G$ carrying a group structure such that the multiplication map $m\colon G \times G \to G$ and the inversion map $i\colon G \to G$ are smooth. Homomorphisms of Lie groups are smooth group homomorphisms. Clearly, Lie groups and their homomorphisms form a category in which direct products exist.

Applying the tangent functor to the defining identities of the group structure $(G, m, i, e)$, it is immediately seen that then $(TG, Tm, Ti, 0_{T_eG})$ is again a Lie group such that $\pi\colon TG \to G$ becomes a homomorphism of Lie groups and such that the zero section $z\colon G \to TG$ also is a homomorphism of Lie groups.

## 9.3 The Lie functor

A vector field $X \in \mathfrak{X}(G)$ is called *left invariant* if, for all $g \in G$, $X \circ l_g = Tl_g \circ X$. In particular, $X(g) = X(l_g(e)) = T_e l_g X(e)$; thus $X$ is uniquely determined by the value $X(e)$, and thus the map

$$\mathfrak{X}(G)^{l_G} \to T_eG, \quad X \mapsto X(e) \tag{33}$$

from the space of left-invariant vector fields into $T_eG$ is injective. It is also surjective: if $v \in T_eG$, then right multiplication with $v$ in $TG$, $Tr_v\colon TG \to TG$ preserves fibers (by (31)) and hence defines a vector field

$$\tilde{v} = Tr_v \circ z\colon G \to TG, \quad g \mapsto T_g r_v(0_g) = Tm(g,v) = T_e l_g(v) \tag{34}$$

which is left invariant since right multiplications commute with left multiplications. Now, the space $\mathfrak{X}(G)^{l_G}$ is a Lie subalgebra of $\mathfrak{X}(M)$, as readily follows from Lemma 8.5 because $X$ is left invariant if and only if the pair $(X,X)$ is $l_g$-related for all $g \in G$. The space $\mathfrak{g} := T_e G$ with the Lie bracket defined by $[v,w] := [\tilde{v},\tilde{w}]_e$ is called *the Lie algebra of $G$*.

**Theorem 9.1** *The Lie bracket $\mathfrak{g} \times \mathfrak{g} \to \mathfrak{g}$ is $C^0$, and every Lie group homomorphism $f : G \to H$ induces a homomorphism $T_e f : \mathfrak{g} \to \mathfrak{h}$ of $C^0$-Lie algebras.*

**Proof.** In order to prove the first statement, pick a chart $\phi : U \to V$ of $G$ such that $\phi(e) = 0$. Since $\tilde{w}(x) = Tm(x,w)$ depends smoothly on $(x,w)$, it is represented in the chart by a smooth map (which again will be denoted by $\tilde{w}(x)$). But this implies that $[\tilde{v},\tilde{w}](x) = d(\tilde{w})(x)\tilde{v}(x) - d(\tilde{v})(x)\tilde{w}(x)$ depends smoothly on $v, w$ and $x$ and hence $[v,w]$ depends smoothly on $v,w$.

In order to prove the second statement, one first has to check that the pair of vector fields $(\tilde{v}, (T_e\phi v)\tilde{\ })$ is $f$-related, and then one applies Lemma 7.2 in order to conclude that $T_e f[v,w] = [T_e fv, T_e fw]$. □

The functor from Lie groups over $\mathbb{K}$ into $C^0$-Lie algebras over $\mathbb{K}$ will be called the *Lie functor* (in a given $C^0$ concept over $\mathbb{K}$). At this point, the fundamental problem arises to describe the image of this functor: when can a $C^0$-Lie algebra over $\mathbb{K}$ be integrated to a Lie group? A similar problem arises in the context of symmetric spaces (see below). A more conceptual version of the Lie functor (also for symmetric spaces) will be described in [8]. In [27] the basic theory of infinite dimensional Lie groups over topological fields is developed; see Section 13 for a list of the main examples.

### 9.4 Symmetric spaces

A *symmetric space over $\mathbb{K}$* is a smooth manifold with a multiplication map $m : M \times M \to M$ such that, for all $x, y, z \in M$,

(M1) $m(x,x) = x$

(M2) $m(x, m(x,y)) = y$, i.e. $l_x^2 = \mathrm{id}_M$,

(M3) $m(x, m(y,z)) = m(m(x,y), m(x,z))$, i.e. $l_x \in \mathrm{Aut}(M,m)$,

(M4) $T_x(l_x) = -\mathrm{id}_{T_x M}$.

The left multiplication operator $l_x$ is, by (M1)–(M3), an automorphism of order two fixing $x$; it is called the *symmetry around $x$* and is usually denoted by $\sigma_x$. In the finite dimensional case over $\mathbb{K} = \mathbb{R}$, (M4) implies by the implicit function theorem that $x$ is an *isolated fixed point* of $\sigma_x$ and hence our definition contains the one from [40] as a special case.

The theory of symmetric spaces will be developed in [8] and in [9]: in analogy with the Lie functor for Lie groups, following the approach from [40], one can define a *Lie functor for symmetric spaces* associating to a symmetric space with base point a $C^0$-Lie triple system. Examples (based on Jordan theory) will be given in [9]. For a good theory of symmetric spaces one has to assume that 2 is invertible in $\mathbb{K}$.

# Part III: Differential Calculus over Topological Fields

In Parts III and IV, we focus on the case where $\mathbb{K}$ is a Hausdorff, non-discrete topological field, where $\mathcal{M}$ is the class of all Hausdorff topological $\mathbb{K}$-vector spaces, $C^0$-maps are the continuous maps, and where the product topology is used on products.

All topological spaces are assumed Hausdorff now. Topological fields are assumed Hausdorff and non-discrete. The term "ultrametric field" refers to a field $\mathbb{K}$, equipped with a non-trivial ultrametric absolute value $|.|: \mathbb{K} \to [0, \infty[$. Totally disconnected, non-discrete locally compact topological fields will be referred to as "local fields."

The goal of Part III is to provide specific results concerning the differential calculus over topological fields. The main result, prepared in Sections 10 and 11, is the following theorem, to be proved in Section 12:

**Theorem.** *Let $(\mathbb{K}, |.|)$ be either $\mathbb{R}$, equipped with the usual absolute value, or an ultrametric field. Let $E$ and $F$ be topological $\mathbb{K}$-vector spaces and $f: U \to F$ be a mapping, defined on a non-empty open subset $U \subseteq E$. Let $k \in \mathbb{N}_0$. If $E$ is metrizable, then $f$ is of class $C_\mathbb{K}^k$ if and only if the composition $f \circ c: \mathbb{K}^{k+1} \to F$ is of class $C_\mathbb{K}^k$, for every smooth mapping $c: \mathbb{K}^{k+1} \to U$.*

This theorem substantially simplifies differential calculus on metrizable topological vector spaces. Note that neither $E$ nor $F$ need to be locally convex here (for the real locally convex case, compare [35]).

## 10 Basic results

In this section, we prove various basic facts of differential calculus over topological fields. We also prove some more specific results, which are the basis for the considerations in Sections 11 and 12.

<div align="center">**Various useful facts**</div>

**Lemma 10.1** *Suppose that $E$ and $F$ are topological $\mathbb{K}$-vector spaces, $U$ an open subset of $E$, and $f: U \to F$ a mapping of class $C^k$, where $k \in \mathbb{N}_0 \cup \{\infty\}$. Let $F_0$ be a vector subspace of $F$ containing the image of $f$. If $F_0$ is closed in $F$ or if $F_0$ is sequentially closed and $\mathbb{K}$ is metrizable, then the co-restriction $f|^{F_0}: U \to F_0$ is of class $C^k$ as a mapping into $F_0$.*

**Proof.** We may assume that $k \in \mathbb{N}_0$; the proof is by induction on $k$. The case $k = 0$ is trivial. Now suppose the assertion holds for $C^k$-maps, where $k \in \mathbb{N}_0$, and let $f: U \to F$ be a $C^{k+1}$-map such that $f(U) \subseteq F_0$. Then apparently $f^{[1]}(x, y, t) \in F_0$ for all $(x, y, t) \in U^{[1]}$ such that $t \neq 0$. If $(x, y, 0) \in U^{[1]}$, we find a net (resp., sequence if $\mathbb{K}$ is metrizable) $(t_\alpha)_\alpha$ in $\mathbb{K}^\times$ such that $t_\alpha \to 0$ in $\mathbb{K}$ and $(x, y, t_\alpha) \in U^{[1]}$ for all $\alpha$. Then the net (resp., sequence) $(f^{[1]}(x, y, t_\alpha))_\alpha$ of elements of $F_0$ converges to $f^{[1]}(x, y, 0)$, entailing that $f^{[1]}(x, y, 0) \in F_0$. We readily deduce that $\operatorname{im} f^{[1]} \subseteq F_0$, whence $f|^{F_0}$ is of class $C^1$, with $(f|^{F_0})^{[1]} = f^{[1]}|^{F_0}$, which is a $C^k$-map by the induction hypothesis. Hence $f|^{F_0}$ is $C^{k+1}$ (Remark 4.2). $\square$

**Lemma 10.2** *Suppose that $E$ is a topological $\mathbb{K}$-vector space, $(F_i)_{i \in I}$ a family of topological $\mathbb{K}$-vector spaces, $U \subseteq E$ an open subset, and $f \colon U \to P$ a mapping, where $P := \prod_{i \in I} F_i$, with canonical projections $\operatorname{pr}_i \colon P \to F_i$. Let $k \in \mathbb{N}_0 \cup \{\infty\}$. Then $f$ is of class $C^k$ if and only if $\operatorname{pr}_i \circ f$ is of class $C^k$ for each $i \in I$.*

**Proof.** We may assume that $k \in \mathbb{N}_0$; the proof is a by trivial induction on $k$, which we leave to the reader (cf. **3.7** and Lemma 4.4). $\square$

**Lemma 10.3** *Suppose that $E$ is a topological $\mathbb{K}$-vector space, $((F_i)_{i \in I}, (\phi_{ij})_{i \leq j})$ a projective system of topological $\mathbb{K}$-vector spaces $F_i$ and continuous $\mathbb{K}$-linear maps $\phi_{ij} \colon F_j \to F_i$ (where $I$ is a directed set), $U \subseteq E$ an open subset, and $f \colon U \to F$ a mapping, where $F = \varprojlim_{i \in I} F_i$, with limit maps $\pi_i \colon F \to F_i$. Let $k \in \mathbb{N}_0 \cup \{\infty\}$. Then $f$ is of class $C^k$ if and only if $\pi_i \circ f$ is of class $C^k$ for each $i \in I$.*

**Proof.** Each $\pi_i$ being a continuous linear map and thus smooth, $\pi_i \circ f$ is of class $C^k$ if $f$ is of class $C^k$ (Proposition 4.5). Now suppose that $\pi_i \circ f$ is of class $C^k$ for each $i \in I$. We may assume that $F$ is realized as a closed vector subspace of $P := \prod_{i \in I} F_i$, and $\pi_i := \operatorname{pr}_i|_F$ (where $\operatorname{pr}_i \colon P \to F_i$ is the corresponding canonical projection). By Lemma 10.2, $f$ is of class $C^k$ as a mapping into $P$, and thus $f$ is $C^k$ also as a map into $F$, by Lemma 10.1. $\square$

**Lemma 10.4** *Let $E$ and $F$ be topological $\mathbb{K}$-vector spaces, $N$ be a closed vector subspace of $E$, and $q \colon E \to E/N =: E_1$ be the quotient map. Let $k \in \mathbb{N}_0 \cup \{\infty\}$, and $f_1 \colon U_1 \to F$ be a mapping, defined on an open subset $U_1 \subseteq E_1$. Let $U \subseteq E$ be an open subset such that $q(U) = U_1$. Then $f_1$ is of class $C^k$ if and only if $f := f_1 \circ q|_U^{U_1} \colon U \to F$ is of class $C^k$.*

**Proof.** The continuous linear map $q$ being smooth, $f$ will be of class $C^k$ if so is $f_1$ (Proposition 4.5). For the converse direction, we may assume that $k \in \mathbb{N}_0$; the proof is by induction.

*The case $k = 0$.* The mapping $q|_U^{U_1} \colon U \to U_1$ is a continuous open surjection and thus a quotient map. Hence, if $f = f_1 \circ q|_U^{U_1}$ is continuous, then so is $f_1$.

*Induction step.* Suppose that the assertion holds for some $k \in \mathbb{N}_0$, and suppose that $f = f_1 \circ q|_U^{U_1}$ is of class $C^{k+1}$. Then, by induction, $f_1$ is a mapping of class $C^k$, entailing that also
$$f_1^{]1[} \colon U_1^{]1[} \to F, \quad f_1^{]1[}(x', y', t) := \tfrac{1}{t}(f_1(x' + ty') - f_1(x'))$$

is a mapping of class $C^k$, where $U_1^{]1[} := \{(x', y', t) \in U_1^{[1]} : t \neq 0\}$. Next, suppose that $(x, y, t) \in U^{[1]} =: W$ and $(u, v, t) \in W$ such that $q(x) = q(u)$ and $q(y) = q(v)$. If $t \neq 0$, then

$$f^{[1]}(x, y, t) = \tfrac{1}{t}(f(x+ty) - f(x)) = \tfrac{1}{t}(f_1(q(x+ty)) - f_1(q(x))) = f^{[1]}(u, v, t)\,.$$

If $t = 0$, we pick a net $(t_\alpha)_\alpha$ in $\mathbb{K}^\times$ converging to 0 in $\mathbb{K}$ such that $(x, y, t_\alpha), (u, v, t_\alpha) \in W$ for all $\alpha$, and deduce from the preceding that $f^{[1]}(x, y, t) = \lim f^{[1]}(x, y, t_\alpha) = \lim f^{[1]}(u, v, t_\alpha) = f^{[1]}(u, v, 0)$. By the preceding, the mapping $g := f^{[1]} : W \to F$ factors over the open continuous surjection (and thus quotient map) $p := (q \times q \times \mathrm{id}_\mathbb{K})|_W^{W_1}$, where $W_1 := (q \times q \times \mathrm{id}_\mathbb{K})(W) \subseteq E_1 \times E_1 \times \mathbb{K}$. Thus, there is a mapping $g_1 : W_1 \to F$ such that $g_1 \circ p = g$. Note that $g_1$ is of class $C^k$, by induction. It is easy to see that $U_1^{[1]} = U_1^{]1[} \cup W_1$ and $f_1^{]1[}|_{U_1^{]1[} \cap W_1} = g_1|_{U_1^{]1[} \cap W_1}$. Thus $h : U_1^{[1]} \to F$, defined via $h|_{U_1^{]1[}} := f_1^{]1[}$, $h|_{W_1} := g_1$, is a mapping of class $C^k$. As $h(x', y', t) = f_1^{]1[}(x', y', t) = \tfrac{1}{t}(f_1(x' + ty') - f_1(x'))$ for all $(x', y', t) \in U_1^{]1[}$, we see that $f_1$ is of class $C^1$, with $f_1^{[1]} = h$. The mapping $h$ being of class $C^k$, we deduce with Remark 4.2 that $f_1$ is of class $C^{k+1}$. □

## Specific results for later use

We now make a number of observations of a more specialized nature, which are vital for the following sections.

The first lemma generalizes [58], Lemma 29.7.

**Lemma 10.5** *Let $E$ be a topological $\mathbb{K}$-vector space, $U \subseteq \mathbb{K}$ be a non-empty, open subset, $n \in \mathbb{N}$, and $f : U \to E$ be a $C_\mathbb{K}^{n-1}$-map. Let $\Delta := \{(x, x, \ldots, x) \in U^{n+1} : x \in U\}$ be the diagonal of $U^{n+1}$. Then there exists a continuous function $h : U^{n+1} \setminus \Delta \to E$ extending $f^{>n<}$.*

**Proof.** By Proposition 6.9, $f$ is a $C_{Sch}^{n-1}$-map. Now copy [58], proof of Lemma 29.7. □

**Lemma 10.6** *Let $X$ be a metrizable topological space, $Y$ be a regular topological space, and $f : D \to Y$ be a continuous mapping, defined on a dense subset $D \subseteq X$. If $\lim_{n \to \infty} f(x_n)$ exists for all sequences $(x_n)_{n \in \mathbb{N}}$ in $D$ which converge in $X$, then $f$ extends to a continuous mapping $\overline{f} : X \to Y$.*

**Proof.** Suppose that $x \in X \setminus D$. If $(x_n)_{n \in \mathbb{N}}$ and $(y_n)_{n \in \mathbb{N}}$ are sequences in $D$ which converge to $x$, then so does $x_1, y_1, x_2, y_2, \ldots$, whence the sequence $f(x_1), f(y_1), f(x_2), f(y_2), \ldots$ converges to some limit $p$ by hypothesis. Then $\lim_{n \to \infty} f(x_n) = \lim_{n \to \infty} f(y_n) = p$. As a consequence, $g(x) := \lim_{n \to \infty} f(x_n)$ is well-defined, independent of the choice of sequence $(x_n)_{n \in \mathbb{N}}$ converging to $x$. Set $g(y) := f(y)$ for $y \in D$. Then clearly $g : D \cup \{x\} \to Y$ is continuous. The assertion is therefore a special case of [19], Exerc. 3.2.B (see also [18]). □

**Lemma 10.7** *Let $\mathbb{K}$ be a metrizable, non-discrete topological field, $U \subseteq \mathbb{K}$ be a non-empty open subset, $E$ be a topological $\mathbb{K}$-vector space, and $f\colon U \to E$ be a map. Let $n \in \mathbb{N}$. If $f$ is continuous and the limit $\lim_{k\to\infty} f^{>n<}(x_k)$ exists for all sequences $(x_k)_{k\in\mathbb{N}}$ in $U^{>n<}$ which converge in $U^{n+1}$, then $f$ is a $C^n_{Sch}$-map (and thus of class $C^n_{\mathbb{K}}$).*

**Proof.** Since $f$ is continuous, a simple induction shows that $f^{>m<}$ is a continuous map, for all $m \in \mathbb{N}_0$. In particular, $f^{>n<}\colon U^{>n<} \to E$ is continuous. Since $U^{>n<}$ is dense in $U^{n+1}$ and $E$, like any topological group, is regular as a topological space, in view of the hypotheses we can apply Lemma 10.6 and obtain a continuous extension $f^{<n>}\colon U^{n+1} \to E$ of $f^{>n<}$. Thus $f$ is a $C^n_{Sch}$-map and hence of class $C^n_{\mathbb{K}}$ (Proposition 6.9). $\square$

In this connection, we shall find the following topological fact very useful:

**Lemma 10.8** *Suppose that $X$ is a metrizable topological space, $Y$ a topological space, and $f\colon S \to Y$ a map, defined on some subset $S$ of $X$. Let $x \in \overline{S}$ be in the closure of $S$. Then the following conditions are equivalent:*

(a) *The limit $\lim_{n\to\infty} f(x_n)$ exists, for every sequence $(x_n)_{n\in\mathbb{N}}$ in $S$ which converges to $x$.*

(b) *For every pair of sequences $(x_n)_{n\in\mathbb{N}}$ and $(y_n)_{n\in\mathbb{N}}$ in $S$ which converge to $x$, there exist subsequences $(x_{n_k})_{k\in\mathbb{N}}$, resp., $(y_{m_k})_{k\in\mathbb{N}}$ such that $f(x_{n_1}), f(y_{m_1}), f(x_{n_2}), f(y_{m_2}), \ldots$ is a convergent sequence in $Y$.*

**Proof.** (a)$\Rightarrow$(b): Supposing that (a) holds, let $(x_n)_{n\in\mathbb{N}}$ and $(y_n)_{n\in\mathbb{N}}$ be sequences in $S$ which converge to $x$. Then also the sequence $x_1, y_1, x_2, y_2, \ldots$ converges to $x$, and hence $f(x_1), f(y_1), f(x_2), f(y_n), \ldots$ converges in $Y$, by (a). Thus (a) entails (b).

(b)$\Rightarrow$(a): Supposing that (b) holds, let $(x_n)_{n\in\mathbb{N}}$ be a sequence in $S$ which converges to $x$. Then (b), applied with $(y_n)_{n\in\mathbb{N}} = (x_n)_{n\in\mathbb{N}}$, entails that $(f(x_n))_{n\in\mathbb{N}}$ has a convergent subsequence $(f(x_{n_k}))_{k\in\mathbb{N}}$, with limit $p$, say. If the sequence $(f(x_n))_{n\in\mathbb{N}}$ did not converge to $p$, then we could find a neighbourhood $U \subseteq Y$ of $p$ and a subsequence $(x_{m_k})_{k\in\mathbb{N}}$ of $(x_n)_{n\in\mathbb{N}}$ such that $f(x_{m_k}) \notin U$ for all $k \in \mathbb{N}$. By the argument used before, we then find a subsequence $(x_{m_{k_j}})_{j\in\mathbb{N}}$ of $(x_{m_k})_{k\in\mathbb{N}}$ such that $f(z_j)$ converges in $Y$ as $j \to \infty$, where we abbreviated $z_j := x_{m_{k_j}}$. Let $q := \lim_{j\to\infty} f(z_j)$. Then $q \neq p$, since $f(z_j) \notin U$ for all $j$. However, by (b), we can find subsequences of $(f(x_{n_k}))_{k\in\mathbb{N}}$ and of $(f(z_j))_{j\in\mathbb{N}}$ which converge to the same limit. Thus $p = q$, and we have reached a contradiction. $\square$

## 11 Continuity versus continuity along curves

The following discussions were inspired by the fact that the topology on a metrizable, real locally convex space $E$ is final with respect to the set $C^\infty(\mathbb{R}, E)$ of smooth curves ([36], Theorem 4.11). In this section, we prove an analogous result for metrizable topological vector spaces over ultrametric fields. Furthermore, in the real case, we are able to remove the hypothesis of local convexity.

**Lemma 11.1 (Special Curve Lemma: Ultrametric Case)** *Let $\mathbb{K}$ be an ultrametric field, and $0 \neq \rho \in \mathbb{K}$ be an element of absolute value $|\rho| < 1$. Let $E$ be a metrizable topological $\mathbb{K}$-vector space, and $U \subseteq E$ be an open subset. Let $x \in U$, and suppose that $(x_n)_{n \in \mathbb{N}}$ and $(y_n)_{n \in \mathbb{N}}$ are sequences in $U$ which converge to $x$. Then there exists an injective, monotonically increasing mapping $m \colon \mathbb{N} \to \mathbb{N}$, $k \mapsto m_k$ and a smooth curve $c \colon \mathbb{K} \to U$, with the following properties:*

(a) $c(\rho^k) = x_{m_k}$ *for all* $k \in 2\mathbb{N} - 1$, *and* $c(\rho^k) = y_{m_k}$ *for all* $k \in 2\mathbb{N}$;

(b) $c(0) = x$;

(c) $c|_{\mathbb{K}^\times}$ *is locally constant;*

(d) *If* $t, s \in \mathbb{K}$ *such that* $|t| = |s|$, *then* $c(t) = c(s)$;

(e) $c(\mathbb{K}^\times) = \{c(\rho^k) \colon k \in \mathbb{N}\}$.

*Note that if $(x_n)_{n \in \mathbb{N}} = (y_n)_{n \in \mathbb{N}}$ in particular, then $(c(\rho^k))_{k \in \mathbb{N}}$ is a subsequence of $(x_n)_{n \in \mathbb{N}}$.*

**Proof.** Without loss of generality, we may assume that $x = 0 \in U$ (the general case follows via translations). Let $\mathbb{O} := \{t \in \mathbb{K} \colon |t| \leq 1\}$ be the valuation ring of $\mathbb{K}$, which is an open, closed subset of $\mathbb{K}$. We recall that a zero-neighbourhood $V \subseteq E$ is called *balanced* if $\mathbb{O}V \subseteq V$. Since $E$ is a metrizable topological $\mathbb{K}$-vector space, there exists a sequence $(W_n)_{n \in \mathbb{N}}$ of balanced, open zero-neighbourhoods such that $\{W_n \colon n \in \mathbb{N}\}$ is a basis for the filter of zero-neighbourhoods in $E$, and furthermore $W_n \subseteq U$ for all $n \in \mathbb{N}$, $W_n + W_n \subseteq W_{n-1}$ if $n \geq 2$, and $\rho^{-2n} W_n \subseteq W_{n-1}$.

We pick $m_1 \in \mathbb{N}$ such that $x_{m_1} \in W_1$ and $y_{m_1} \in W_1$. This is possible since $\lim_{n \to \infty} x_n = \lim_{n \to \infty} y_n = 0$. Inductively, having defined natural numbers $m_1 < m_2 < \cdots < m_k$ for some $k \in \mathbb{N}$ such that $x_{m_j} \in W_j$ and $y_{m_j} \in W_j$ for all $j \in \{1, \ldots, k\}$, we pick $m_{k+1} \in \mathbb{N}$ such that $m_{k+1} > m_k$ and such that $x_{m_{k+1}} \in W_{k+1}$ and $y_{m_{k+1}} \in W_{k+1}$. In this way, we obtain a monotonically increasing, injective mapping $m \colon \mathbb{N} \to \mathbb{N}$, $k \mapsto m_k$ such that, for all $k \in \mathbb{N}$, we have $x_{m_k} \in W_k$ and $y_{m_k} \in W_k$.

Note that $U_k := \rho^k \mathbb{O} \setminus \rho^{k+1} \mathbb{O} = \{t \in \mathbb{K} \colon |\rho|^{k+1} < |t| \leq |\rho|^k\}$ is an open and closed subset of $\mathbb{K}$ for each $k \in \mathbb{N}$, and

$$\rho \mathbb{O} = \{0\} \cup \bigcup_{k \in \mathbb{N}} U_k,$$

where the union is disjoint. We define $c \colon \mathbb{K} \to U$ via

$$c(t) := \begin{cases} 0 & \text{if } t = 0; \\ x_{m_k} & \text{if } t \in U_k \text{ and } k \in 2\mathbb{N} - 1; \\ y_{m_k} & \text{if } t \in U_k \text{ and } k \in 2\mathbb{N}; \\ x_{m_1} & \text{if } |t| > |\rho|. \end{cases}$$

Then $c$ is constant on each of the open sets $U_k$ and also on $\mathbb{K} \setminus \rho \mathbb{O}$, entailing that $c|_{\mathbb{K}^\times}$ is locally constant and thus a mapping of class $C_\mathbb{K}^\infty$. Furthermore, $c$ satisfies (a)–(e) by construction. It only remains to show that $c$ is smooth. Given $r \in ]0, \infty[$, there exists a unique number $\nu(r) \in \mathbb{Z}$ such that $|\rho|^{\nu(r)+1} < r \leq |\rho|^{\nu(r)}$.

**Claim.** *For each $n \in \mathbb{N}_0$, the following holds*:

(i) $c$ *is of class* $C_{\mathbb{K}}^n$.

(ii) *We have* $c^{>n<}(t) \in W_{\nu(\|t\|_\infty)-2n}$, *for all* $0 \neq t = (t_1, \ldots, t_{n+1}) \in \mathbb{K}^{>n<} \subseteq \mathbb{K}^{n+1}$ *of supremum norm* $\|t\|_\infty \leq |\rho|^{4n+1}$ *(meaning that $\nu(\|t\|_\infty) \geq 4n+1$).*

(iii) *If $n \geq 1$ and $t = (t_1, \ldots, t_{n+1}) \in \mathbb{K}^{>n<}$ such that $|t_1| = |t_2| = \cdots = |t_{n+1}|$, then $c^{>n<}(t) = 0$.*

The proof is by induction on $n \in \mathbb{N}_0$. To prove the case $n = 0$, suppose that $0 \neq t \in \mathbb{K}^{>0<} = \mathbb{K}$ such that $|t| \leq |\rho|$ and thus $\nu(|t|) \geq 1$. Then $t \in U_{\nu(|t|)}$ and thus, by definition of $c$, we have $c^{>0<}(t) = c(t) \in \{x_{\nu(|t|)}, y_{\nu(|t|)}\} \subseteq W_{\nu(|t|)}$, showing that (ii) holds. Since $\nu(|t|) \to \infty$ as $t \to 0$ and the $W_j$'s form a descending fundamental sequence of zero-neighbourhoods, we deduce that $c(t) \to 0 = c(0)$ as $t \to 0$. Thus $c$ is continuous at 0 and thus continuous (being smooth on $\mathbb{K}^\times$). Thus $c$ is a $C^0$-map, establishing (i).

*Induction step.* Suppose that the claim holds for some $n \in \mathbb{N}_0$. Let $t = (t_1, t_2, \ldots, t_{n+2}) \in \mathbb{K}^{>n+1<}$. Then $t \neq 0$. If $|t_1| = |t_2| = \cdots = |t_{n+2}|$, then

$$c^{>n+1<}(t) = \frac{c^{>n<}(t_1, t_3, \cdots, t_{n+2}) - c^{>n<}(t_2, t_3, \ldots, t_{n+2})}{t_1 - t_2} = 0,$$

using that $c^{>n<}(t_1, t_3, \cdots, t_{n+2}) = c^{>n<}(t_2, t_3, \ldots, t_{n+2})$ by (d) in the case $n = 0$, resp., by (iii) for the map $c^{>n<}$ if $n \geq 1$ (valid by the induction hypothesis). Thus (iii) also holds for $c^{>n+1<}$.

To see that (ii) holds for $c^{>n+1<}$, let $t = (t_1, \ldots, t_{n+2}) \in \mathbb{K}^{>n+1<}$ such that $\|t\|_\infty \leq |\rho|^{4(n+1)+1}$. If $|t_1| = |t_2| = \cdots = |t_{n+2}|$, then $c^{>n+1<}(t) = 0 \in W_{\nu(\|t\|_\infty)-2(n+1)}$. It remains to consider the case when not all absolute values $|t_i|$ are equal. Since $c^{>n+1<}$ is symmetric in its $n+2$ variables, we may assume that $|t_1| = \|t\|_\infty$ and $|t_1| > |t_2|$ without loss of generality. Then $|t_1 - t_2| = |t_1| = \|t\|_\infty$ and thus, abbreviating $r := (t_1, t_3, \ldots, t_{n+2}) \in \mathbb{K}^{>n<}$ and $s := (t_2, t_3, \ldots, t_{n+2}) \in \mathbb{K}^{>n<}$, we have $\|r\|_\infty, \|s\|_\infty \leq \|t\|_\infty \leq |\rho|^{4(n+1)+1} \leq |\rho|^{4n+1}$, whence $c^{>n<}(r) \in W_{\nu(\|r\|_\infty)-2n} \subseteq W_{\nu(\|t\|_\infty)-2n}$ by the induction hypothesis and similarly[4] $c^{>n<}(s) \in W_{\nu(\|t\|_\infty)-2n}$. Thus

$$\begin{aligned}
c^{>n+1<}(t) &= \frac{c^{>n<}(r) - c^{>n<}(s)}{t_1 - t_2} \in \frac{1}{t_1 - t_2} \cdot \left(W_{\nu(\|t\|_\infty)-2n} + W_{\nu(\|t\|_\infty)-2n}\right) \\
&\subseteq \frac{1}{t_1 - t_2} \cdot W_{\nu(\|t\|_\infty)-2n-1} = \frac{\rho^{\nu(\|t\|_\infty)+1}}{t_1 - t_2} \cdot \rho^{-(\nu(\|t\|_\infty)+1)} \cdot W_{\nu(\|t\|_\infty)-2n-1} \\
&\subseteq \rho^{-(\nu(\|t\|_\infty)+1)} \cdot W_{\nu(\|t\|_\infty)-2n-1} \\
&\subseteq W_{\nu(\|t\|_\infty)-2n-2} = W_{\nu(\|t\|_\infty)-2(n+1)}.
\end{aligned}$$

---
[4]If $n = 0$, it may happen that $s = 0$. In this case, we cannot use the induction hypothesis, but observe directly that $c^{>0<}(s) = c^{>0<}(0) = c(0) = 0 \in W_{\nu(\|t\|_\infty)-2n}$.

Here, the penultimate inclusion holds since $W_{\nu(\|t\|_\infty)-2n-1}$ is balanced and $|\frac{\rho^{\nu(\|t\|_\infty)+1}}{t_1-t_2}| = \frac{|\rho|^{\nu(\|t\|_\infty)+1}}{|t_1-t_2|} = \frac{|\rho|^{\nu(\|t\|_\infty)+1}}{\|t\|_\infty} \leq 1$. To see that the final inclusion holds, abbreviate $N := \nu(\|t\|_\infty) - 2n - 1$. Then

$$N \geq 4(n+1) + 1 - 2n - 1 = 2n + 4 \tag{35}$$

by choice of $t$. Furthermore, $\nu(\|t\|_\infty) + 1 = N + 2n + 2 \leq 2N$, using (35). Thus

$$\rho^{-(\nu(\|t\|_\infty)+1)} \cdot W_{\nu(\|t\|_\infty)-2n-1} = \rho^{-(\nu(\|t\|_\infty)+1)} \cdot W_N \subseteq \rho^{-2N} \cdot W_N \subseteq W_{N-1} = W_{\nu(\|t\|_\infty)-2n-2},$$

as asserted. Thus (ii) is established for $c^{>n+1<}$.

As $c$ is of class $C_\mathbb{K}^n$ by induction, there exists a continuous map $h: \mathbb{K}^{n+2} \setminus \Delta \to E$ extending $c^{>n+1<}$, where $\Delta := \{(t,t,\ldots,t) \in \mathbb{K}^{n+2} : t \in \mathbb{K}\}$ (see Lemma 10.5). On the other hand, as $g := c|_{\mathbb{K}^\times}$ is smooth, the map $g^{>n+1<}$ extends to a continuous map $g^{<n+1>} : (\mathbb{K}^\times)^{n+2} \to E$. We define a mapping $c^{<n+1>}: \mathbb{K}^{n+2} \to E$ via

$$c^{<n+1>}(t) := \begin{cases} h(t) & \text{if } t \in \mathbb{K}^{n+2} \setminus \Delta \\ g^{<n+1>}(t) & \text{if } t \in (\mathbb{K}^\times)^{n+2} \\ 0 & \text{if } t = (0,0,\ldots,0) \end{cases}$$

for $t \in \mathbb{K}^{n+2}$. To see that $c^{<n+1>}$ is well-defined, note that $D := (\mathbb{K}^\times)^{>n+1<}$ is a dense subset of $\mathbb{K}^{n+2}$ and hence also of the open subset $V := (\mathbb{K}^\times)^{n+2} \cap (\mathbb{K}^{n+2} \setminus \Delta)$ of $\mathbb{K}^{n+2}$, which is the intersection of the domains of definition of $h$ and $g^{<n+1>}$. Since $g^{<n+1>}|_D = g^{>n+1<} = c^{>n+1<}|_D = h|_D$, the continuity of $g^{<n+1>}$ and $h$ entails that $g^{<n+1>}|_V = h|_V$. Thus $c^{<n+1>}$ is well-defined.

Note that $\mathbb{K}^{n+2}$ is a metrizable topological space and $E$, being a topological vector space, is a regular topological space. Since $c^{<n+1>}|_{(\mathbb{K}^\times)^{n+2}} = g^{<n+1>}$ and $c^{<n+1>}|_{\mathbb{K}^{n+2} \setminus \Delta} = h$ are continuous, $D$ is dense in $\mathbb{K}^{n+2}$, and $(\mathbb{K}^\times)^{n+2} \cup (\mathbb{K}^{n+2} \setminus \Delta) = \mathbb{K}^{n+2} \setminus \{0\}$, Lemma 10.6 shows that $c^{<n+1>}$ will be continuous if we can show that $c^{<n+1>}(t_j) \to c^{<n+1>}(0) = 0$ in $E$, for every sequence $(t_j)_{j\in\mathbb{N}}$ in $D = (\mathbb{K}^\times)^{>n+1<} \subseteq \mathbb{K}^{n+2}$ converging to 0 in $\mathbb{K}^{n+2}$. But this readily follows from (ii). In fact, given a zero-neighbourhood $Q$ in $E$, there exists $\ell \in \mathbb{N}$ such that $W_\ell \subseteq Q$. There is $j_0 \in \mathbb{N}$ such that $\|t_j\|_\infty \leq |\rho|^{4(n+1)+1}$ and $\|t_j\|_\infty \leq |\rho|^{\ell+2(n+1)}$ (and thus $\nu(\|t_j\|_\infty) - 2(n+1) \geq \ell$) for all $j \in \mathbb{N}$ such that $j \geq j_0$. For any such $j$, (ii) gives $c^{>n+1<}(t_j) \in W_{\nu(\|t_j\|_\infty)-2(n+1)} \subseteq W_\ell \subseteq Q$. Thus $c^{<n+1>}(t_j) = c^{>n+1<}(t_j) \to 0$ as $j \to \infty$ indeed. Thus $c^{<n+1>}$ is continuous and thus $c$ is of class $C_\mathbb{K}^{n+1}$ (see Proposition 6.9). This completes the inductive proof of the claim. The claim being true, $c$ is smooth. $\square$

We can also create smooth curves in metrizable real topological vector spaces (which need not be locally convex). Compare [21, Proposition 4.2.15] and [36, p. 18] for closely related results for real locally convex spaces.

**Lemma 11.2 (Special Curve Lemma: Real Case)** *Let $E$ be a metrizable topological real vector space and $U \subseteq E$ be an open subset. Let $x \in U$, and suppose that $(x_n)_{n\in\mathbb{N}}$ and $(y_n)_{n\in\mathbb{N}}$ are sequences in $U$ which converge to $x$. Let $\rho \in\, ]0,1[$. Then there exists an injective, monotonically increasing mapping $m: \mathbb{N} \to \mathbb{N}$, $k \mapsto m_k$ and a smooth curve $c: \mathbb{R} \to U$, with the following properties:*

(a) $c(\rho^k) = x_{m_k}$ for all $k \in 2\mathbb{N} - 1$, and $c(\rho^k) = y_{m_k}$ for all $k \in 2\mathbb{N}$;

(b) $c(0) = x$;

(c) $c(t) = c(-t)$ for all $t \in \mathbb{R}$;

(d) $c([\rho^{k+1}, \rho^k]) = \mathrm{conv}\{c(\rho^{k+1}), c(\rho^k)\}$ is the line segment joining $c(\rho^{k+1})$ and $c(\rho^k)$ in $E$, for all $k \in \mathbb{N}$;

(e) $c|_{[\rho, \infty[}$ is constant.

**Proof.** Without loss of generality, $x = 0 \in U$. There exists $a \in \mathbb{N}$ such that $\rho^a \leq 1 - \rho$. We choose a decending sequence $(W_n)_{n \in \mathbb{N}_0}$ of open, balanced zero-neighbourhoods $W_n \subseteq U$ of $E$ such that $\{W_n : n \in \mathbb{N}_0\}$ is a basis of zero-neighbourhoods, $W_n + W_n \subseteq W_{n-1}$ for all $n \in \mathbb{N}$, and $\rho^{-(n+1+a)n} W_n \subseteq W_{n-1}$. Then $\rho^{-2n} W_n \subseteq W_{n-1}$, using that $W_n$ is balanced and $n + 1 + a \geq 2$. We now define an injective, monotonically increasing function $m \colon \mathbb{N} \to \mathbb{N}$, $k \mapsto m_k$ as described in the proof of Lemma 11.1. Given $k \in \mathbb{N}$, we define $z_k := x_{m_k}$ if $k$ is odd, $z_k := y_{m_k}$ if $k$ is even. Fix $\varepsilon \in ]0, \tfrac{1}{2}[$. We pick a smooth, monotonically increasing function $\tau \colon \mathbb{R} \to \mathbb{R}$ such that $\tau|_{]-\infty, \varepsilon]} = 0$ and $\tau|_{[1-\varepsilon, \infty[} = 1$. We define $c(0) := 0$ and $c(t) := x_{m_1}$ for $t \in ]\rho, \infty[$. Given $t \in ]0, \rho[$, there is a unique $k \in \mathbb{N}$ such that $\rho^{k+1} < t \leq \rho^k$; we define $c(t)$ as the convex combination

$$c(t) = \left(1 - \tau\left(\frac{t - \rho^{k+1}}{\rho^k(1-\rho)}\right)\right) z_{k+1} + \tau\left(\frac{t - \rho^{k+1}}{\rho^k(1-\rho)}\right) z_k \in W_{k+1} + W_k \subseteq W_{k-1} \subseteq U. \quad (36)$$

Finally, set $c(t) := c(-t)$ for $t < 0$. We have defined a function $c \colon \mathbb{R} \to U$. It is clear that $c$ is smooth on $]-\infty, -\rho[$ and $]\rho, \infty[$, and also on some open neighbourhood of $-\rho^k$ and on some open neighbourhood of $\rho^k$, for each $k \in \mathbb{N}$ (since $c$ is constant there). Furthermore, $c$ is smooth on $]\rho^{k+1}, \rho^k[$ and $]-\rho^k, -\rho^{k+1}[$ for each $k \in \mathbb{N}$, thanks to the smoothness of $\tau$. As a consequence, $c|_{\mathbb{R}^\times}$ is smooth. Furthermore, $c$ satisfies (a)–(e) by construction. It only remains to show that $c$ is smooth on all of $\mathbb{R}$.

To this end, we observe first that $c$ is constantly $z_k$ on $[\rho^k - \varepsilon(\rho^k - \rho^{k+1}), \rho^k + \varepsilon(\rho^{k-1} - \rho^k)]$, entailing that (36) describes $c(t)$ for a given $k \in \mathbb{N}$ not only for $t \in ]\rho^{k+1}, \rho^k]$, but actually for all $t \in [\rho^{k+1} - \varepsilon \rho^{k+1}(1-\rho), \rho^k + \varepsilon \rho^{k-1}(1-\rho)]$.

Given $n \in \mathbb{N}$, there exists $C_n \in [0, \infty[$ such that

$$|\tau^{<n>}(t)| \leq C_n \quad \text{and} \quad |(1-\tau)^{<n>}(t)| \leq C_n \quad \text{for all } t \in [-\varepsilon\rho, 1 + \tfrac{1}{\rho}\varepsilon]^{n+1}, \quad (37)$$

the set $[-\varepsilon\rho, 1 + \tfrac{1}{\rho}\varepsilon]^{n+1}$ being compact.

Next, given $k \in \mathbb{N}$, consider the affine-linear bijection $\alpha_k \colon \mathbb{R} \to \mathbb{R}$, $\alpha_k(t) = \frac{t - \rho^{k+1}}{\rho^k(1-\rho)}$. For $n \in \mathbb{N}_0$, let $\alpha_{k,n} \colon \mathbb{R}^{>n<} \to \mathbb{R}^{>n<}$ be the map $(t_1, \ldots, t_{n+1}) \mapsto (\alpha_k(t_1), \ldots, \alpha_k(t_{n+1}))$. Then

$$(\tau \circ \alpha_k)^{>1<}(t_1, t_2) = \frac{(\tau \circ \alpha_k)(t_1) - (\tau \circ \alpha_k)(t_2)}{t_1 - t_2} = \frac{1}{\rho^k(1-\rho)} \frac{\tau(\alpha_k(t_1)) - \tau(\alpha_k(t_2))}{\alpha_k(t_1) - \alpha_k(t_2)}$$

for $(t_1, t_2) \in \mathbb{R}^{>1<}$ and thus $(\tau \circ \alpha_k)^{>1<} = (\rho^k(1-\rho))^{-1} \cdot (\tau^{>1<} \circ \alpha_{k,1})$. A trivial induction shows that $(\tau \circ \alpha_k)^{>n<} = (\rho^k(1-\rho))^{-n} \cdot (\tau^{>n<} \circ \alpha_{k,n})$ for all $n \in \mathbb{N}_0$, and similarly $((1-\tau) \circ \alpha_k)^{>n<} = (\rho^k(1-\rho))^{-n}((1-\tau)^{>n<} \circ \alpha_{k,n})$. As a consequence, for all $n \in \mathbb{N}_0$ and $t = (t_1, \ldots, t_{n+1}) \in [\rho^{k+1} - \varepsilon\rho^{k+1}(1-\rho), \rho^k + \varepsilon\rho^{k-1}(1-\rho)]^{n+1} \cap \mathbb{R}^{>n<}$, we have

$$
\begin{aligned}
c^{>n<}(t) &= ((1-\tau) \circ \alpha_k)^{>n<}(t) \cdot z_{k+1} + (\tau \circ \alpha_k)^{>n<}(t) \cdot z_k \\
&\in (\rho^k(1-\rho))^{-n}[-C_n, C_n] \cdot z_{k+1} + (\rho^k(1-\rho))^{-n}[-C_n, C_n] \cdot z_k \\
&\subseteq (\rho^k(1-\rho))^{-n} C_n W_{k+1} + (\rho^k(1-\rho))^{-n} C_n W_k \\
&\subseteq (\rho^k(1-\rho))^{-n} C_n W_{k-1}.
\end{aligned} \tag{38}
$$

A very similar argument shows that this conclusion is also valid for all $t = (t_1, \ldots, t_{n+1}) \in (-[\rho^{k+1} - \varepsilon\rho^{k+1}(1-\rho), \rho^k + \varepsilon\rho^{k-1}(1-\rho)]^{n+1}) \cap \mathbb{R}^{>n<}$.

Given $r \in ]0, \infty[$, there is a unique $\nu(r) \in \mathbb{Z}$ such that $\rho^{\nu(r)+1} < r \leq \rho^{\nu(r)}$. We claim:

**Claim.** *For each $n \in \mathbb{N}_0$, the following holds*:

(i) *$c$ is of class $C^n_\mathbb{R}$.*

(ii) *There exists $L_n > 0$ such that, for all $0 \neq t = (t_1, \ldots, t_{n+1}) \in \mathbb{R}^{>n<}$ of supremum norm $\|t\|_\infty \leq \rho^{4n+2}$, we have $c^{>n<}(t) \in L_n \cdot W_{\nu(\|t\|_\infty) - 2n - 1}$.*

The proof is by induction on $n \in \mathbb{N}_0$. To prove the case $n = 0$, suppose that $0 \neq t \in \mathbb{K}^{>0<} = \mathbb{K}$ such that $|t| \leq \rho^2$. Set $k := \nu(|t|) \geq 2$. Then $t \in ]\rho^{k+1}, \rho^k] \cup [-\rho^k, -\rho^{k+1}[$ and thus $c^{>0<}(t) = c(t) \in W_{k-1} = W_{k-2 \cdot 0 - 1}$ (see (36)), showing that (ii) holds with $L_0 := 1$. Since $\nu(|t|) \to \infty$ as $t \to 0$ and the $W_j$'s (for $j \in \mathbb{N}$) form a descending fundamental sequence of zero-neighbourhoods, we deduce that $c(t) \to 0 = c(0)$ as $t \to 0$. Thus $c$ is continuous at $0$ and thus continuous. We have shown that $c$ is a $C^0$-map: (i) is established.

*Induction step.* Suppose that the claim holds for some $n \in \mathbb{N}_0$. Let $t = (t_1, t_2, \ldots, t_{n+2}) \in \mathbb{R}^{>n+1<}$ such that $\|t\|_\infty \leq \rho^{4n+6}$. Then $t \neq 0$. We abbreviate $k := \nu(\|t\|_\infty)$. There are two cases. First, suppose that

$$
\begin{aligned}
\{t_1, \ldots, t_{n+2}\} &\subseteq [\rho^{k+1} - \varepsilon\rho^{k+1}(1-\rho), \rho^k] \quad \text{or} \\
\{t_1, \ldots, t_{n+2}\} &\subseteq -[\rho^{k+1} - \varepsilon\rho^{k+1}(1-\rho), \rho^k]
\end{aligned} \tag{39}
$$

holds. Then $c^{>n+1<}(t) \in (\rho^k(1-\rho))^{-(n+1)} C_{n+1} W_{k-1} \subseteq \rho^{-(k+a)(n+1)} C_{n+1} W_{k-1}$ by (38). Note that $n+1 \leq k-1$ here since $k \geq 4n+6$. Hence $\rho^{-(k+a)(n+1)} W_{k-1} \subseteq \rho^{-((k-1)+a+1)(k-1)} W_{k-1} \subseteq W_{k-2} \subseteq W_{k-2(n+1)-1}$, whence indeed (ii) holds, provided we choose $L_{n+1} \geq C_{n+1}$.

If, on the other hand, (39) is not satisfied, then there exist $i, j \in \{1, \ldots, n+2\}$ such that $i \neq j$ and $|t_i - t_j| > \varepsilon\rho^{k+1} \cdot (1-\rho)$. Since $c^{>n+1<}$ is symmetric in its $n+2$ variables, we may assume that $i = 1$ and $j = 2$. Abbreviate $r := (t_1, t_3, \ldots, t_{n+2}) \in \mathbb{R}^{>n<}$ and

$s := (t_2, t_3, \ldots, t_{n+2}) \in \mathbb{R}^{>n<}$. Then $\|r\|_\infty, \|s\|_\infty \leq \|t\|_\infty \leq \rho^{4n+6} \leq \rho^{4n+2}$ and hence[5]

$$
\begin{aligned}
c^{>n+1<}(t) &= \frac{c^{>n<}(r) - c^{>n<}(s)}{t_1 - t_2} \\
&\in \frac{1}{t_1 - t_2} \left( L_n W_{\nu(\|r\|_\infty)-2n-1} + L_n W_{\nu(\|s\|_\infty)-2n-1} \right) \\
&\subseteq \frac{L_n}{t_1 - t_2} (W_{k-2n-1} + W_{k-2n-1}) \subseteq \frac{L_n}{t_1 - t_2} W_{k-2n-2} \\
&\subseteq \frac{L_n}{\varepsilon(1-\rho)} \rho^{-(k+1)} W_{k-2n-2} \subseteq \frac{L_n}{\varepsilon(1-\rho)} W_{k-2n-2-1} = \frac{L_n}{\varepsilon(1-\rho)} W_{k-2(n+1)-1}.
\end{aligned}
$$

To see that the final inclusion holds, note that $k+1 \leq 2(k-2n-2)$ because $k = \nu(\|t\|_\infty) \geq 4n+6$. Hence $\rho^{-(k+1)} W_{k-2n-2} \subseteq \rho^{-2(k-2n-2)} W_{k-2n-2} \subseteq W_{k-2n-2-1}$ by choice of the fundamental sequence of balanced zero-neighbourhoods.

Thus, in both possible cases, (ii) holds for $c^{>n+1<}$ if we define $L_{n+1} := \max\{C_{n+1}, \frac{L_n}{\varepsilon \cdot (1-\rho)}\}$. Using (ii) and the induction hypothesis that $c$ is of class $C_\mathbb{R}^n$, we now see as in the ultrametric case (proof of Lemma 11.1) that $c^{>n+1<}$ extends to a continuous mapping $\mathbb{R}^{n+2} \to E$, whence $c$ is of class $C_\mathbb{R}^{n+1}$.

Thus the claim is established. Consequently, $c$ is smooth. The proof is complete. $\square$

**Theorem 11.3** *Let $(\mathbb{K}, |.|)$ be either $\mathbb{R}$, equipped with the usual absolute value, or an arbitrary ultrametric field. Let $E$ be a metrizable topological $\mathbb{K}$-vector space and let $U \subseteq E$ be a non-empty open subset. Then we have:*

(a) *The given topology on $U$ (induced by $E$) is the final topology with respect to the set $C^\infty(\mathbb{K}, U) := C^\infty(\mathbb{K}, E) \cap U^\mathbb{K}$ of smooth curves $c \colon \mathbb{K} \to U$.*

(b) *A mapping $f \colon U \to Z$ from $U$ to a topological space $Z$ is continuous if and only if $f \circ c \colon \mathbb{K} \to Z$ is continuous for each smooth curve $c \colon \mathbb{K} \to U$.*

**Proof.** Let $\mathcal{O}$ be the given topology on $U$, and $\mathcal{O}_f$ be the final topology on $U$ with respect to the set $C^\infty(\mathbb{K}, U)$. Every smooth curve $c \colon \mathbb{K} \to U$ being continuous, we have $\mathcal{O} \subseteq \mathcal{O}_f$. The converse inclusion $\mathcal{O}_f \subseteq \mathcal{O}$ (and thus the validity of (a)) will follow if we can prove (b) (as we may take $Z = (U, \mathcal{O}_f)$ and $f \colon U \to Z$, $x \mapsto x$ in particular).

Thus, let $f \colon U \to Z$ be a map from $U$ to a topological space $Z$. If $f$ is continuous, then $f \circ c$ is continuous for each $c \in C^\infty(\mathbb{K}, U)$. Conversely, suppose that $f \circ c$ is continuous for each $c \in C^\infty(\mathbb{K}, U)$. Note that a sequence $(z_n)_{n\in\mathbb{N}}$ in the topological space $Z$ converges to $z \in Z$ if and only if every subsequence of $(z_n)_{n\in\mathbb{N}}$ has a subsequence which converges to $z$. Thus, $U$ being metrizable, the map $f$ will be continuous if we can show that every convergent sequence $(x_n)_{n\in\mathbb{N}}$ in $U$, with limit $x$, say, has a subsequence $(x_{n_k})_{k\in\mathbb{N}}$ such that $(f(x_{n_k}))_{k\in\mathbb{N}}$ converges to $f(x)$. Now, given $(x_n)_{n\in\mathbb{N}}$, the Special Curve Lemma 11.1 (resp., 11.2) applied with $(y_n)_{n\in\mathbb{N}} = (x_n)_{n\in\mathbb{N}}$ provides a smooth curve $c\colon \mathbb{K} \to U$ and an injective, monotonically increasing function $\mathbb{N} \to \mathbb{N}$, $k \mapsto m_k$ such that $c(0) = x$ and $c(\rho^k) = x_{m_k}$ for all $k \in \mathbb{N}$, for some element $0 \neq \rho \in \mathbb{K}$ such that $|\rho| < 1$. Since $f \circ c$ is continuous by hypothesis, and $\rho^k \to 0$ in $\mathbb{K}$, we deduce that $f(x_{m_k}) = (f \circ c)(\rho^k) \to (f \circ c)(0) = f(x)$ as $k \to \infty$. Thus $(x_{m_k})_{k\in\mathbb{N}}$ is a subsequence of $(x_n)_{n\in\mathbb{N}}$ such that $f(x_{m_k}) \to f(x)$ as $k \to \infty$. In view of the above considerations, we conclude that $f$ is continuous. □

---
[5]If $n = 0$ and $s = 0$, pass directly to the third line.

Compare also [3] and [4] for related results in the real case.

## 12 Testing whether $f$ is $C^k$ via smooth maps on $\mathbb{K}^{k+1}$

It is well-known that a mapping $f\colon U \to F$ on an open subset $U$ of a real Fréchet space $E$, with values in a Mackey complete locally convex real topological vector space $F$, is smooth if and only if it is smooth in the sense of convenient differential calculus, *viz.* $f \circ c$ is smooth for every smooth curve $c\colon \mathbb{R} \to U$ (cf. [36], Theorem 4.11 (1); if $E$ is finite-dimensional and $F = \mathbb{R}$, this is a classical theorem by Boman [11]). In this section, we show that variants of this fact are valid for *not necessarily locally convex* topological vector spaces over $\mathbb{K} = \mathbb{R}$, and also for topological vector spaces over ultrametric fields $\mathbb{K}$, when smooth curves $\mathbb{K} \to U$ are replaced with smooth mappings $\mathbb{K}^k \to U$ for suitable $k \in \mathbb{N}$.

We start with some preparatory definitions and considerations.

**12.1** Given a non-discrete topological field $\mathbb{K}$, topological $\mathbb{K}$-vector space $E$, and open subset $U \subseteq E$, we define $U^{]0[} := U$, $U^{]1[} := U^{[1]} \cap (E \times E \times \mathbb{K}^\times)$, and inductively $U^{]k[} := (U^{]k-1[})^{[1]} \cap (E^{[k-1]} \times E^{[k-1]} \times \mathbb{K}^\times)$. If $F$ is a topological $\mathbb{K}$-vector space and $f\colon U \to F$ a map, we define $f^{]0[} := f$. Recursively, given $k \in \mathbb{N}$ and having defined $f^{]j[}\colon U^{]j[} \to F$ for $j = 0, \ldots, k-1$, we define $f^{]k[}\colon U^{]k[} \to F$ via

$$f^{]k[}(x, y, t) := \frac{1}{t}(f^{]k-1[}(x + ty) - f^{]k-1[}(x)) \quad \text{for all } (x, y, t) \in U^{]k[} \subseteq E^{[k-1]} \times E^{[k-1]} \times \mathbb{K}^\times.$$

It is clear that $U^{]k[}$ is a dense, open subset of $U^{[k]}$ for each $k \in \mathbb{N}_0$. Furthermore, if $f$ is continuous, then so is $f^{]k[}$ for each $k \in \mathbb{N}_0$.

**12.2** In the preceding situation, the map $f$ is of class $C^k_\mathbb{K}$ (for $k \in \mathbb{N}_0$) if and only if $f$ is continuous and $f^{]k[}$ extends to a continuous map $g\colon U^{[k]} \to F$; in this case, $g = f^{[k]}$.

[In fact, if $f$ is $C^k_\mathbb{K}$, then $g := f^{[k]}$ is a continuous extension of $f^{]k[}$, and the only such by density of $U^{]k[}$ in $U^{[k]}$. The converse direction can be proved by induction on $k$. The case $k = 0$ is trivial. Now suppose that the assertion holds for $k$ and suppose that $f^{]k+1[}$ has a continuous extension $g\colon U^{[k+1]} \to F$. Pick $x_0 \in U^{]k[}$. Then $f^{]k[}(x) = f^{]k[}(x_0) + f^{]k+1[}(x_0, x - x_0, 1)$ for all $x \in U^{]k[}$ shows that the continuous map $h\colon U^{[k]} \to F$, $h(x) := f^{]k[}(x_0) + g(x_0, x - x_0, 1)$ extends $f^{]k[}$. By induction, $f$ is $C^k$ with $h = f^{[k]}$. Then $(f^{[k]})^{]1[} = g|_{(U^{[k]})^{]1[}}$, since both functions are continuous and coincide with $f^{]k+1[}$ on the dense subset $U^{]k+1[}$ of their domain of definition. Hence $f^{[k]}$ is $C^1$ with $(f^{[k]})^{[1]} = g$, and thus $f$ is $C^{k+1}$ with $f^{[k+1]} = (f^{[k]})^{[1]} = g$.]

**Lemma 12.3** *Let $E$ and $F$ be topological vector spaces over a non-discrete topological field $\mathbb{K}$ and $f\colon U \to F$ be a map, defined on an open subset $U \subseteq E$. Let $m \in \mathbb{N}$, $n \in \mathbb{N}_0$, and $\gamma\colon V \to U^{[n]}$ be a smooth map, defined on an open subset $V \subseteq \mathbb{K}^m$. Let $D \subseteq V$ be a subset such that $\gamma(D) \subseteq U^{]n[}$, and suppose that $X_0 \subseteq \overline{D}$ is a non-empty finite subset of the closure of $D$ in $V$. Then there exists a smooth map $\Gamma\colon W \to U$, defined on an open subset $W \subseteq \mathbb{K}^{m+n}$, an open neighbourhood $Y$ of $X_0$ in $V$, and a smooth map $g\colon Y \to W^{[n]}$ such that $g(D \cap Y) \subseteq W^{]n[}$ and*

$$(\forall x \in D \cap Y) \quad f^{]n[}(\gamma(x)) = (f \circ \Gamma)^{]n[}(g(x)).$$

**Proof.** We prove the assertion for all $n \in \mathbb{N}_0$, by induction.

The case $n = 0$ is trivial: we can take $W := Y := V$, $\Gamma := \gamma$, and $g := \mathrm{id}_V$.

*Induction step.* Suppose that $n \in \mathbb{N}$, and suppose that the assertion of the lemma is correct if $n$ is replaced with $n-1$. Let $\gamma\colon V \to U^{[n]}$, $D$, and $X_0 \subseteq \overline{D}$ be as described in the lemma.

According to $E^{[n]} = E^{[n-1]} \times E^{[n-1]} \times \mathbb{K}$, we have $\gamma = (\gamma_1, \gamma_2, \gamma_3)$ with smooth coordinate functions $\gamma_1$, $\gamma_2$, and $\gamma_3$. For all $x \in D$, we have $\gamma(x) \in U^{]n[}$ and therefore $\gamma_1(x) \in U^{]n-1[}$, $\gamma_1(x) + \gamma_3(x)\gamma_2(x) \in U^{]n-1[}$, and $\gamma_3(x) \in \mathbb{K}^\times$. Furthermore, by definition of $f^{]n[}$,

$$(\forall x \in D) \quad f^{]n[}(\gamma(x)) = \frac{1}{\gamma_3(x)} \left( f^{]n-1[}(\gamma_1(x) + \gamma_3(x)\gamma_2(x)) - f^{]n-1[}(\gamma_1(x)) \right). \qquad (40)$$

Define $X_1 := \{(x, \gamma_3(x)) : x \in X_0\} \cup \{(x, 0) : x \in X_0\} \subseteq V \times \mathbb{K}$. There exists an open neighbourhood $V_1$ of $X_1$ in $V \times \mathbb{K}$ such that $\eta(x,t) := \gamma_1(x) + t\gamma_2(x) \in U^{[n-1]}$ for all $(x,t) \in V_1$. Then $\eta\colon V_1 \to U^{[n-1]}$ is smooth. We define $D_1 := (\{(x,0) : x \in D\} \cup \{(x, \gamma_3(x)) : x \in D\}) \cap V_1$. Then $\eta(D_1) \subseteq \gamma_1(D) \cup \{\gamma_1(x) + \gamma_3(x)\gamma_2(x) : x \in D\} \subseteq U^{]n-1[}$ by the above observations. If $x \in X_0$ and $(x_\alpha)$ is a net in $D$ converging to $x$, then $(x_\alpha, 0) \in V_1$ and $(x_\alpha, \gamma_3(x_\alpha)) \in V_1$ eventually (as $V_1$ is an open neighbourhood of $(x,0)$ and $(x, \gamma_3(x))$, where $\gamma_3$ is continuous). Thus $(x_\alpha, 0) \in D_1$ and $(x_\alpha, \gamma_3(x_\alpha)) \in D_1$ eventually, entailing that $(x, 0) \in \overline{D_1}$ and $(x, \gamma_3(x)) \in \overline{D_1}$. Thus $X_1 \subseteq \overline{D_1}$.

By the induction hypothesis, there exists an open subset $W \subseteq \mathbb{K}^{(n-1)+(m+1)} = \mathbb{K}^{n+m}$ and smooth map $\Gamma\colon W \to U$, an open neighbourhood $Y_1$ of $X_1$ in $V_1$, and a smooth map $h\colon Y_1 \to W^{[n-1]}$ such that $h(D_1 \cap Y_1) \subseteq W^{]n-1[}$ and

$$(\forall x \in D_1 \cap Y_1) \quad f^{]n-1[}(\eta(x)) = (f \circ \Gamma)^{]n-1[}(h(x)). \qquad (41)$$

There exists an open neighbourhood $Y \subseteq V$ of $X_0$ such that $Y \times \{0\} \subseteq Y_1$ and $\{(x, \gamma_3(x)) \colon x \in Y\} \subseteq Y_1$. Then $(x, 0) \in D_1$ and $(x, \gamma_3(x)) \in D_1$ for all $x \in Y \cap D$ and thus

$$\begin{aligned}
f^{]n[}(\gamma(x)) &= \frac{1}{\gamma_3(x)} \left( f^{]n-1[}(\eta(x, \gamma_3(x))) - f^{]n-1[}(\eta(x, 0)) \right) \\
&= \frac{1}{\gamma_3(x)} \left( (f \circ \Gamma)^{]n-1[}(h(x, \gamma_3(x))) - (f \circ \Gamma)^{]n-1[}(h(x, 0)) \right) \\
&= \frac{1}{\gamma_3(x)} \left( (f \circ \Gamma)^{]n-1[}(h(x, 0) + \gamma_3(x)\gamma_3(x)^{-1}(h(x, \gamma_3(x)) - h(x, 0))) \right. \\
&\qquad \left. - (f \circ \Gamma)^{]n-1[}(h(x, 0)) \right) \\
&= (f \circ \Gamma)^{]n[}(h(x, 0), h^{]1[}((x, 0), (0, 1), \gamma_3(x)), \gamma_3(x)),
\end{aligned}$$

using (40) to obtain the first equality and (41) to obtain the second. Note that $g \colon Y \to W^{[n]}$, $g(x) := (h(x, 0), h^{[1]}((x, 0), (0, 1), \gamma_3(x)), \gamma_3(x))$ actually maps into $W^{[n]}$, is smooth, and takes $Y \cap D$ into $W^{]n[}$. Furthermore, by the preceding $f^{]n[}(\gamma(x)) = (f \circ \Gamma)^{]n[}(g(x))$ for all $x \in Y \cap D$. Thus the assertion is established for $f$. This completes the proof. $\square$

We are now in the position to prove the theorem announced at the beginning of Part III.

**Theorem 12.4** *Let $(\mathbb{K}, |.|)$ be either $\mathbb{R}$, equipped with the usual absolute value, or an arbitrary ultrametric field. Let $E$ and $F$ be topological $\mathbb{K}$-vector spaces and $f \colon U \to F$ be a mapping, defined on a non-empty open subset $U \subseteq E$. If $E$ is metrizable, then for each $k \in \mathbb{N}_0$, the following conditions are equivalent:*

(a) *$f$ is a mapping of class $C_\mathbb{K}^k$.*

(b) *The composition $f \circ c \colon \mathbb{K}^{k+1} \to F$ is of class $C_\mathbb{K}^k$, for every smooth mapping $c \colon \mathbb{K}^{k+1} \to U$.*

*In particular, $f$ is smooth if and only if $f \circ c$ is smooth, for every $k \in \mathbb{N}$ and every smooth map $c \colon \mathbb{K}^k \to U$.*

**Proof.** By the Chain Rule, (a) implies (b). We prove that (b) implies (a) by induction on $k \in \mathbb{N}_0$. For $k = 0$, (b) entails (a) by Theorem 11.3. Thus, assume that $k \in \mathbb{N}$ and assume that (b) implies (a) when $k$ is replaced with $k - 1$. Let $f \colon U \to F$ be a map satisfying (b). If $c \colon \mathbb{K}^k \to U$ is smooth, then so is $c_0 \colon \mathbb{K}^{k+1} \to U$, $c_0(t_1, \ldots, t_{k+1}) := c(t_1, \ldots, t_k)$, and thus $f \circ c_0$ is of class $C_\mathbb{K}^k$ by (b). Since $(f \circ c)(t_1, \ldots, t_k) = (f \circ c_0)(t_1, \ldots, t_k, 0)$, we deduce that also $f \circ c$ is of class $C_\mathbb{K}^k$, whence $f$ is of class $C_\mathbb{K}^{k-1}$ by the induction hypothesis and thus continuous. Therefore $f$ will be of class $C^k$ if we can show that $f^{]k[}$ extends to a continuous map $U^{[k]} \to F$ (see **12.2**). In view of Lemmas 10.6 and 10.8, to obtain the continuous extension, we only need to prove the following claim:

*Claim. For every pair of convergent sequences $(x_n)_{n \in \mathbb{N}}$ and $(y_n)_{n \in \mathbb{N}}$ in $U^{]k[}$ with the same limit, there are subsequences of $(f^{]k[}(x_n))_{n \in \mathbb{N}}$ and $(f^{]k[}(y_n))_{n \in \mathbb{N}}$ converging to the same limit $x \in U^{[k]}$.*

To establish the claim, note first that by Lemma 11.1, resp., 11.2, there is a smooth curve $\gamma\colon \mathbb{K} \to U^{[k]}$ such that $(\gamma(\rho^{2j-1}))_{j\in\mathbb{N}}$ is a subsequence of $(x_n)_{n\in\mathbb{N}}$ and $(\gamma(\rho^{2j}))_{j\in\mathbb{N}}$ is a subsequence of $(y_n)_{n\in\mathbb{N}}$, for a suitable element $0 \neq \rho \in \mathbb{K}$ such that $|\rho| < 1$. Applying Lemma 12.3 with $m = 1$, $D := \{\rho^j : j \in \mathbb{N}\}$, and $X_0 := \{0\}$, we find a smooth map $\Gamma\colon W \to U$, defined on an open subset $W \subseteq \mathbb{K}^{k+1}$, an open zero-neighbourhood $Y$ in $\mathbb{K}$, and a smooth map $g\colon Y \to W^{[k]}$ such that $g(D \cap Y) \subseteq W^{]k[}$ and

$$f^{]k[}(\gamma(\rho^j)) = (f \circ \Gamma)^{]k[}(g(\rho^j)) \quad \text{for all } j \geq j_0, \tag{42}$$

where we have chosen $j_0 \in \mathbb{N}$ so large that $\rho^j \in Y$ for all $j \geq j_0$. Here $f \circ \Gamma$ is of class $C^k$. In fact, if $z \in W$, there exists a smooth function $\chi_z\colon \mathbb{K}^{k+1} \to W$ such that $\chi_z(w) = w$ for all $w$ in some open neighbourhood $W_z \subseteq W$ of $z$.[6] Then $\Gamma_z := \Gamma \circ \chi_z\colon \mathbb{K}^{k+1} \to U$ is a smooth mapping defined on all of $\mathbb{K}^{k+1}$, whence $f \circ \Gamma_z$ is of class $C^k$ by hypothesis, and hence so is $f \circ \Gamma|_{W_z} = f \circ \Gamma_z|_{W_z}$. Being locally $C^k$, the map $f \circ \Gamma$ is of class $C^k$. Now $(f \circ \Gamma)^{[k]}$ being continuous, we observe that the right hand side of (42) converges to $(f \circ \Gamma)^{[k]}(g(0))$ as $j \to \infty$. The left hand side of (42) provides subsequences of $(f^{]k[}(x_n))_{n\in\mathbb{N}}$ (for odd $j$) and $(f^{]k[}(y_n))_{n\in\mathbb{N}}$ (for even $j$). By the preceding, both of them converge to $(f \circ \Gamma)^{[k]}(g(0))$. Thus the claim is established, and the proof is complete. $\square$

Compositions with smooth maps $\mathbb{R}^k \to E$, for *locally convex* $E$ and $F$, have also been considered in [35], for similar purposes. Cf. also Souriau's theory of diffeological spaces [60], [42].

**Remark 12.5** Let $\mathbb{K}$ be as before. Given a topological $\mathbb{K}$-vector space $E$, let us write $c^\infty(E)$ for $E$, equipped with the final topology with respect to the set $C^\infty(\mathbb{K}, E)$ of smooth curves. If $E$ and $F$ are topological $\mathbb{K}$-vector spaces and $U$ an open subset of $c^\infty(E)$, let us call a mapping $f\colon U \to F$ *conveniently smooth* or a $c^\infty$-*map* if $f \circ c$ is smooth, for every smooth curve $c\colon \mathbb{K} \to E$ with image in $U$. For metrizable $E$, we have $E = c^\infty(E)$ by Theorem 11.3, and it is natural to wonder whether a mapping $f$ as before is smooth if and only if it is conveniently smooth (as in the real locally convex case). We hope to explore this question, as well as a potential "ultrametric convenient differential calculus" and its category-theoretical properties, in subsequent research.

# Part IV: Examples of Lie Groups over Topological Fields

## 13 Some classes of examples

In this section, we describe concrete examples of Lie groups over topological fields, to illustrate the abstract theory presented so far. As we shall see, all of the major construction principles for infinite-dimensional real or complex Lie groups carry over to more general topological fields. For proofs, the reader is referred to [27].

---

[6]If $\mathbb{K}$ is an ultrametric field, simply define $\chi(w) := w$ on some open and closed neighbourhood of $z$ in $\mathbb{K}^{k+1}$ which is contained in $W$, and define $\chi(w) := z$ elsewhere. The real case is standard.

## 13.1 Linear Lie groups

Paradigms of real or complex Lie groups are linear Lie groups, *i.e.*, unit groups of unital Banach algebras (or other well-behaved topological algebras) and their Lie subgroups. If $\mathbb{K}$ is a general topological field, the right class of topological algebras to look at are the *continuous inverse algebras* (or CIAs), *i.e.*, unital associative topological $\mathbb{K}$-algebras $A$ such that the group of units $A^\times$ is open in $A$ and such that the inversion map $\iota \colon A^\times \to A$, $a \mapsto a^{-1}$ is continuous. Then $\iota$ is smooth, and thus $A^\times$ is a $\mathbb{K}$-Lie group. For example, $\mathbb{K}$ is a CIA. If $A$ is a CIA over $\mathbb{K}$, then so is the matrix algebra $M_n(A)$, for each $n \in \mathbb{N}$. Every finite-dimensional unital associative $\mathbb{K}$-algebra is a CIA when equipped with the canonical vector topology ($\cong \mathbb{K}^d$). If $K$ is a compact topological space and $A$ a CIA over $\mathbb{K}$, then the algebra $C(K, A)$ of continuous $A$-valued maps is a CIA, with respect to pointwise operations and the topology of uniform convergence. See [27] for further examples of CIAs over general topological fields $\mathbb{K}$, [24] for a detailed discussion of real and complex locally convex CIAs from the point of view of Lie theory. Further examples have been compiled in [30, 1.15]; cf. also [61].

## 13.2 Mapping groups

The second widely studied class of infinite-dimensional real Lie groups are the mapping groups, for example, loop groups $C(\mathbb{S}^1, G)$ and $C^\infty(\mathbb{S}^1, G)$, where $\mathbb{S}^1$ is the unit circle and $G$ a finite-dimensional real Lie group ([45], [53]).

The classical constructions of mapping groups can be generalized to a large extent to the case of Lie groups over topological fields. The following results can be obtained:

**Proposition 13.1 (Groups of continuous mappings)** *Let $X$ be a topological space, $K \subseteq X$ a compact subset, $\mathbb{K}$ a topological field, and $G$ a $\mathbb{K}$-Lie group. Consider*

$$C_K(X, G) := \{\gamma \in C(X, G) \colon \gamma|_{X \setminus K} = 1\},$$

*the group of continuous $G$-valued maps supported in $K$ (with pointwise group operations). Then there is a uniquely determined smooth manifold structure on $C_K(X, G)$ making it a $\mathbb{K}$-Lie group, and such that*

$$C_K(X, U) \to C_K(X, V) \subseteq C_K(X, L(G)), \quad \gamma \mapsto \phi \circ \gamma$$

*defines a chart of $C_K(X, G)$ around $1$, for a chart $\phi \colon G \supseteq U \to V \subseteq L(G)$ of $G$ such that $\phi(1) = 0$. Here $C_K(X, L(G))$ carries the topology of uniform convergence, $C_K(X, U) := \{\gamma \in C_K(X, G) \colon \operatorname{im} \gamma \subseteq U\}$, and $C_K(X, V) := \{\gamma \in C_K(X, L(G)) \colon \operatorname{im} \gamma \subseteq V\}$.* □

In particular, $C(K, G) = C_K(K, G)$ is a $\mathbb{K}$-Lie group for each $\mathbb{K}$-Lie group $G$ and compact topological space $K$.

**Proposition 13.2 (Groups of differentiable mappings)** *Let $\mathbb{K}$ be a locally compact topological field, $\mathbb{L}$ a topological extension field of $\mathbb{K}$, $r \in \mathbb{N}_0 \cup \{\infty\}$, $M$ a finite-dimensional* $C^r$-manifold over $\mathbb{K}$, $K \subseteq M$ be a compact subset, and $G$ an $\mathbb{L}$-Lie group. Let $C_K^r(M,G)$ be the group of $G$-valued $C_{\mathbb{K}}^r$-maps $\gamma$ on $M$ such that $\gamma|_{M\setminus K} = 1$, with pointwise operations. Then there is a uniquely determined smooth $\mathbb{L}$-manifold structure on $C_K^r(M,G)$ making it an $\mathbb{L}$-Lie group and such that

$$C_K^r(M,U) \to C_K^r(M,V) \subseteq C_K^r(M, L(G)), \quad \gamma \mapsto \phi \circ \gamma$$

defines a chart of $C_K^r(M,G)$ around $1$, for some chart $\phi \colon G \supseteq U \to V \subseteq L(G)$ of $G$ such that $\phi(1) = 0$, using the natural vector topology on $C_K^r(M, L(G))$.

If $M$ is $\sigma$-compact here, $\mathbb{K} \neq \mathbb{C}$, and if the topology on $\mathbb{L}$ arises from an absolute value, then furthermore there is a natural $\mathbb{L}$-Lie group structure on the "test function group" $C_c^r(M,G) := \bigcup_K C_K^r(M,G)$, modelled on the direct limit of topological $\mathbb{L}$-vector spaces $C_c^r(M, L(G)) = \varinjlim C_K^r(M, L(G))$, uniquely determined by a condition analogous to the preceding ones (see [27]; cf. [1] and [49] for the case where $G$ is a finite-dimensional real Lie group, [23] for the case where $\mathbb{K} = \mathbb{R}$, $\mathbb{L} \in \{\mathbb{R}, \mathbb{C}\}$ and $L(G)$ is locally convex).

## 13.3 Diffeomorphism groups

If $\mathbb{K}$ is a local field (of arbitrary characteristic) and $M$ a $\sigma$-compact finite-dimensional smooth $\mathbb{K}$-manifold, then the group $\mathrm{Diff}(M)$ of all $C^\infty$-diffeomorphisms of $M$ can be made a $\mathbb{K}$-Lie group, and in fact in two ways, either modelled on the space

$$C_c^\infty(M, TM) = \varinjlim_K C_K^\infty(M,TM) = \varinjlim_K \varprojlim_{p \in \mathbb{N}_0} C_K^p(M,TM)$$

of compactly supported smooth vector fields, equipped with the LF-topology, or modelled on the same vector space, equipped however with the coarser topology

$$C_c^\infty(M,TM) = \bigcap_{p \in \mathbb{N}_0} C_c^p(M,TM) = \varprojlim_{p \in \mathbb{N}_0} \varinjlim_K C_K^p(M,TM)$$

(see [27]). Based on different notions of smooth and $C^k$-maps (which involve additional estimates and boundedness conditions), groups of special types of diffeomorphisms of manifolds over local fields of characteristic zero have also been considered in the works of S. V. Ludkovsky (see, e.g., [43]). Cf. also [39], [52], [44], [31], [45] for the classical constructions of Lie group structures on diffeomorphism groups of finite-dimensional real manifolds.

## 13.4 Direct limit groups

Consider an ascending sequence $G_1 \subseteq G_2 \subseteq \cdots$ of finite-dimensional Lie groups over $\mathbb{R}$ or a local field $\mathbb{K}$, each a closed submanifold of the next. Then $G := \bigcup_{n \in \mathbb{N}} G_n$ can be given a $\mathbb{K}$-Lie group structure modelled on the locally convex direct limit $\varinjlim L(G_n)$, making it the direct limit of the given directed sequence in the category of $\mathbb{K}$-Lie groups (modelled on not necessarily locally convex topological $\mathbb{K}$-vector spaces) and smooth homomorphisms [27]. This result extends earlier studies in [48], [49], and [25].

**Remark 13.3** It would be interesting to know whether, analogous to the real case, every finite-dimensional smooth $p$-adic Lie group (in the sense considered in this article) has a compatible $p$-adic analytic Lie group structure (and hence is a $p$-adic Lie group in the usual sense, as in [16] or [59]). For a finite-dimensional smooth $p$-adic Lie group $G$, it is not even clear whether every tangent vector at the identity element is tangent to a (local) $p$-adic one-parameter subgroup of $G$, since in general we cannot hope for existence (nor uniqueness) of solutions to differential equations. This closely resembles the situation familiar from infinite-dimensional real Lie theory (beyond the Banach case), where it is still unknown whether every Lie group modelled on a Mackey complete locally convex space possesses an exponential function (for non-Mackey complete counterexamples, see [24]).

**Remark 13.4** Part III and IV, Proposition 6.9, Section 7.2 and much of Section 7.1 were contributed solely by the second author, and may become part of his Habilitationsschrift.

# A  Vector-valued forms and polynomial mappings

In this appendix, $\mathbb{K}$ denotes a commutative ring with unit 1, and $E, F$ are $\mathbb{K}$-modules. For $v \in E$, the (*first order*) *difference operator* $\Delta_v$, acting on maps $q \colon E \to F$, is defined by

$$\Delta_v q(x) := q(x+v) - q(x);$$

then $\Delta_v q$ is again a map $E \to F$; iterating, we define the *k-th order linearization* of $q$ by

$$L_k q(x; v_1, \ldots, v_k) := (\Delta_{v_1} \ldots \Delta_{v_k} q)(x).$$

Since the difference operators $\Delta_v$ and $\Delta_w$ commute, this is symmetric in $v_1, \ldots, v_k$.

**Definition A.1** An *$F$-valued form of degree $k$* is a map $q \colon E \to F$ which is homogeneous of degree $k$ (i.e., $q(tx) = t^k q(x)$ for all $t \in \mathbb{K}$, $x \in E$) and such that

$$L^k q(0; v_1, \ldots, v_k)$$

is $\mathbb{K}$-multilinear in $v_1, \ldots, v_k$.

For $k = 2$ our definition coincides with the usual definition of an $F$-valued quadratic form (cf. [14]). We will derive some consequences of the last condition. In a first step, let us forget the multiplication by scalars from $\mathbb{K}$ and consider $E, F$ just as $\mathbb{Z}$-modules, i.e., as abelian groups.

**Lemma A.2** *For a map $f \colon E \to F$ the following are equivalent:*

(0) *$f$ is affine over $\mathbb{Z}$, i.e., $f(x) = a(x) + b$ with $a \colon E \to F$ additive and $b \in F$.*

(1) *$L_2 f = 0$, i.e., for all $u, v, x \in E$, $\Delta_u \Delta_v f(x) = 0$.*

(2) *$L_2 f(0; \cdot) = 0$, i.e., for all $u, v \in E$, $\Delta_u \Delta_v f(0) = 0$.*

(3) *For all $x \in E$, the map $E \to F$, $v \mapsto \Delta_v f(x)$ is additive.*

(4) *The map $E \to F$, $v \mapsto \Delta_v f(0)$ is additive.*

**Proof.** Assume $f$ is affine over $\mathbb{Z}$. Then $\Delta_v f(x) = f(x+v) - f(x) = a(v)$, hence $\Delta_v f$ is additive as a function of $v$ and constant as a function of $x$. Thus (3) and (4) hold. Since $\Delta_w c = 0$ for all constants, (1) and (2) follow.

Assume, conversely, that (4) holds. Set $a(x) := f(x) - f(0)$ and $b := f(0)$. Then $\Delta_v a = \Delta_v f$. By assumption, $v \mapsto \Delta_v f(0) = \Delta_v a(0) = a(v) - a(0) = a(v)$ is additive. Thus $f(x) = a(x) + b$ is affine over $\mathbb{Z}$.

Assume (2) holds. As above, let $a(x) := f(x) - f(0)$ and $b := f(0)$. By (2), $\Delta_u \Delta_v a(0) = 0$ for all $u, v$:

$$0 = a(u+v) - a(u) - a(v) + a(0) = a(u+v) - a(u) - a(v).$$

It follows that $a$ is additive, and hence $f(x) = a(x) + b$ is affine over $\mathbb{Z}$.

Since clearly (3) implies (4) and (1) implies (2), the lemma is proved. $\square$

**Proposition A.3** *For a map $q \colon E \to F$ the following are equivalent:*

(1) $L_{k+1} q = 0$, *i.e., for all $v_0, \ldots, v_k, x \in E$, $\Delta_{v_0} \ldots \Delta_{v_k} q(x) = 0$.*

(2) $L_{k+1} q(0; \cdot) = 0$, *i.e., for all $v_0, \ldots, v_k \in E$, $\Delta_{v_0} \ldots \Delta_{v_k} q(0) = 0$.*

(3) *For all $x \in E$, the map $E^k \to F$, $(v_1, \ldots, v_k) \mapsto \Delta_{v_1} \ldots \Delta_{v_k} q(x)$ is $k$-additive.*

(4) *The map $E^k \to F$, $(v_1, \ldots, v_k) \mapsto \Delta_{v_1} \ldots \Delta_{v_k} q(0)$ is $k$-additive.*

**Proof.** Trivially, (3) implies (4) and (1) implies (2).

Assume (2). Let $f := \Delta_{v_2} \ldots \Delta_{v_k} q$. Then $f$ fulfils Condition (2) of Lemma A.2 and hence is affine. By Condition (3) of Lemma A.2, $v_1 \mapsto \Delta_{v_1} f(x) = \Delta_{v_1} \ldots \Delta_{v_k} q(x)$ is additive in $v_1$. By symmetry, it is then additive in each argument $v_1, \ldots, v_k$, and hence (3) holds.

Assume (4). Let $f := \Delta_{v_2} \ldots \Delta_{v_k} q$. Then $\Delta_{v_1} f(0)$ is additive in $v_1$ and hence $f$ is affine by Lemma A.2. But then $\Delta_{v_0} \Delta_{v_1} f = 0$ by part (1) of the lemma, and hence (1) holds. $\square$

**Lemma A.4** *Let $q \colon E \to F$ be homogeneous of degree $k$ over $\mathbb{K}$. Then, for all $x \in E$,*

$$L_k q(0; x, \ldots, x) = k! q(x).$$

**Proof.** Since both sides depend only on the values of $q$ on the line $\mathbb{K}x$, we may replace $q$ by the map $\tilde{q}(t) = q(tx) = t^k q(x)$, i.e., we may assume that $E = \mathbb{K}$, $x = 1$ and $q(t) = t^k$. But then $\Delta_1 q(t) = (t+1)^k - t^k = kt^{k-1} +$ lower order terms, and the claim is proved by a straightforward induction. $\square$

**Definition A.5** A map $f : E \to F$ between $\mathbb{K}$-modules $E, F$ is called a *homogeneous polynomial mapping of degree $k$* if, for every $\mathbb{K}$-linear surjection $\phi : \tilde{E} \to E$ from a free $\mathbb{K}$-module $\tilde{E}$ onto $E$, there is a $k$-multilinear map $m : \tilde{E}^k \to F$ with $f(\phi(x)) = m(x, \ldots, x)$. The map $f$ is called a *polynomial mapping* if it is a sum of homogeneous polynomial mappings. Note that $\tilde{E} \to F$, $x \mapsto m(x, \ldots, x)$ is a homogeneous polynomial mapping in the sense of [13] ch. 4, par. 5, no. 9 (where it is always assumed that the departure space is a free module). In view of [13] ch. 4, par. 5, Prop. 13, this means that, given any system of generators $(e_i)_{i \in I}$ of $E$, the function $f$ is given by an expression of the form

$$f(\sum_{i \in I} t_i e_i) = \sum_{|\alpha|=k} t^\alpha a_\alpha,$$

using the usual multi-index notation $t^\alpha := \prod_{i \in I} t_i^{\alpha_i}$ and $|\alpha| := \sum_i \alpha_i$ for a multi-index $\alpha$ which vanishes almost everywhere. The coefficients $a_\alpha \in F$ of course depend on $(e_i)_{i \in I}$.

If $E$ itself is free, it is easily seen that our definition coincides with the one from [13]. Comparing with Definition A.1, we see that every homogeneous polynomial map of degree $k$ is a homogeneous form. For $k = 0, 1$ the converse is clear, and for $k = 2$ we have:

**Lemma A.6** *If $q : E \to F$ is a quadratic $F$-valued form, then it is a homogeneous polynomial map of degree $2$.*

**Proof.** In case $E$ is free, this is proved in [14] Ch. 9, par. 3, Prop. 2. The general case follows by applying this to the quadratic form $q \circ \phi$ where the linear map $\phi$ is as in A.5. □

Note that the proof of the quoted result from [14] does not immediately carry over to general $k$; thus it is not clear whether the analogue of A.6 holds for general $k$. For *smooth* forms we prove this in 5.5. Finally, we remark that a more formal concept of polynomial mappings between modules has been defined by N. Roby in [54] (see also the appendix of [41] for the basic facts). Quadratic maps in the sense of [54] give rise to quadratic maps as defined here; this is proved in [54], but for maps of higher degree no analogues are given. It would be desirable to have general algebraic results clarifying the relations between the various concepts.

# References

[1] Albeverio, S., J. Høegh-Krohn, J. A. Marion, D. H. Testard, and B. S. Torrésani, "Non-commutative Distributions," Marcel-Dekker, New York, 1993.

[2] Averbukh, V. I. and O. G. Smolyanov, *The various definitions of the derivative in linear topological spaces*, Russ. Math. Surv. **23** (1968), 67–113.

[3] Balanzat, M. M., *La différential en les espacios métricas affines*, Math. Notae **9** (1949), 29–51.

[4] ——, *La différentielle d'Hadamard-Fréchet dans les espaces vectoriels topologiques*, C. R. Acad. Sci. Paris **251** (1960), 2459–2461.

[5] Beggs, E., *De Rham's theorem for infinite-dimensional manifolds*, Quart. J. Math. **38** (1987), 131–154.

[6] Bertram, W., "The Geometry of Jordan and Lie Structures," Springer LNM **1754**, Berlin, 2000.

[7] ——, *Generalized projective geometries: general theory and equivalence with Jordan structures*, Adv. Geom. **2** (2002), 329–369.

[8] ——, *Differential geometry, Lie groups and symmetric spaces over general base fields and rings* (work in progress).

[9] Bertram, W. and K.-H. Neeb, *Projective completions of Jordan pairs* (work in progress).

[10] Bochnak, J. and J. Siciak, *Analytic functions in topological vector spaces*, Studia Math. **39** (1971), 77–112.

[11] Boman, J., *Differentiability of a function and of its compositions with functions of one variable*, Math. Scand. **20** (1967), 249–268.

[12] Bourbaki, N., "Topological Vector Spaces," Chapters 1–5, Springer-Verlag, 1987.

[13] ——, "Algèbre. Chapitres 4–5," Hermann, Paris, 1971.

[14] ——, "Algèbre. Chapitres 8–9," Hermann, Paris, 1971.

[15] ——, "Variétés différentielles et analytiques. Fascicule de résultats," Hermann, Paris, 1967.

[16] ——, "Lie Groups and Lie Algebras," Chapters 1–3, Springer-Verlag, 1989.

[17] Cartan, H., "Calcul Différentiel", Herrmann, Paris 1971.

[18] Dieudonné, J., *Sur les espaces uniformes complets*, Ann. Sci. École Normale Sup. **56** (1939), 277–291.

[19] Engelking, R., "General Topology," Heldermann Verlag, 1989.

[20] van Est, W. T. and T. J. Korthagen, *Non-enlargible Lie algebras*, Indag. Math. **26** = Proc. Kon. Ned. Akad. v. Wet. A **67** (1964), 15–31.

[21] Frölicher, A. and A. Kriegl, "Linear Spaces and Differentiation Theory," John Wiley, 1988.

[22] Glöckner, H., *Infinite-dimensional Lie groups without completeness restrictions*, pp. 43–59 in "Geometry and Analysis on Finite- and Infinite-Dimensional Lie Groups," Strasburger, A. et al. (Eds.), Banach Center Publications **55**, 2002.

[23] ——, *Lie group structures on quotient groups and universal complexifications for infinite-dimensional Lie groups*, J. Funct. Analysis **194** (2002), 347–409.

[24] ——, *Algebras whose groups of units are Lie groups*, Studia Math. **153** (2002), 147–177.

[25] ——, *Direct limit Lie groups and manifolds*, Math. J. Kyoto Univ. **43**, No. 1, to appear.

[26] ——, *Implicit functions from topological vector spaces to Banach spaces*, Preprint, TU Darmstadt, 2003.

[27] ——, *Lie groups over non-discrete topological fields*, Preprint, TU Darmstadt, 2003.

[28] ——, *Examples of differentiable mappings into non-locally convex spaces*, Preprint, TU Darmstadt, 2003.

[29] ——, "Infinite-Dimensional Analysis," Lecture Notes, TU Darmstadt, 2003.

[30] Gramsch, B., *Relative Inversion in der Störungstheorie von Operatoren und $\Psi$-Algebren*, Math. Ann. **269** (1984), 22–71.

[31] Hamilton, R., *The inverse function theorem of Nash and Moser*, Bull. Amer. Math. Soc. **7** (1982), 65–222.

[32] Holdgrün, H. S., "Analysis I," Leins Verlag, Göttingen, 2000.

[33] Keller, H. H., "Differential Calculus in Locally Convex Spaces," Springer Verlag, 1974.

[34] Koecher, M., *Gruppen und Lie-Algebren von rationalen Funktionen*, Math. Z. **109** (1969), 349–392.

[35] Kriegl, A., "Eine Theorie glatter Mannigfaltigkeiten und Vektorbündel," Doctoral Dissertation, University of Vienna, 1980.

[36] Kriegl, A., and P. W. Michor, "The Convenient Setting of Global Analysis," AMS, Providence, 1997.

[37] Lang, S., "Fundamentals of Differential Geometry," Graduate Texts in Math. **191**, Springer-Verlag, Berlin, 1999.

[38] Laugwitz, D., *Ist Differentialrechnung ohne Grenzwertbegriff möglich?*, Math. Phys. Semesterberichte **20** (1973), 182–201.

[39] Leslie, J., *On a differential structure for the group of diffeomorphisms*, Topology **6** (1967), 263–271.

[40] Loos, O., "Symmetric Spaces I: General Theory," Benjamin, New York, 1969.

[41] ——, "Jordan Pairs," Springer Lecture Notes in Mathematics **460**, Berlin, 1975.

[42] Losik, M. V., *Fréchet manifolds as diffeologic spaces*, Russ. Math. **36** (1992), No. 5, 31–37.

[43] Ludkovsky, S. V., *Measures on groups of diffeomorphisms of non-archimedian Banach manifolds*, Russian Math. Surv. **51** (1996), No. 2, 338–340.

[44] Michor, P. W., "Manifolds of Differentiable Mappings," Shiva Publishing, 1980.

[45] Milnor, J., *Remarks on infinite-dimensional Lie groups*, pp. 1008–1057 in: DeWitt, B., and R. Stora (Eds.), "Relativity, Groups and Topology II," North Holland, 1983.

[46] Monna, A. F., "Analyse Non-Archimédienne," Springer-Verlag, 1979.

[47] Moerdijk, I. and G.E. Reyes, "Models for Smooth Infinitesimal Analysis," Springer-Verlag, 1991.

[48] Natarajan, L., E. Rodríguez-Carrington and J. A. Wolf, *Differentiable structure for direct limit groups*, Lett. Math. Phys. **23** (1991), 99–109.

[49] ——, *Locally convex Lie groups*, Nova J. Algebra Geom. **2** (1993), 59–87.

[50] Neeb, K.-H., *Central extensions of infinite-dimensional Lie groups*, Ann. Inst. Fourier (Grenoble) **52** (2002), 1365–1442.

[51] ——, *Nancy lectures on infinite-dimensional Lie groups*, TU Darmstadt, Preprint **2203**, March 2002.

[52] Omori, H., "Infinite-Dimensional Lie Groups," Amer. Math. Soc., 1997.

[53] Pressley, A. and G. Segal, "Loop Groups," Oxford University Press, 1986.

[54] Roby, N., *Lois polynômes et lois formelles en théorie de modules*, Ann. Sci. Ec. Norm. Sup. 3e. série t. **80** (1963), 213–348.

[55] van Rooij, A. C. M., "Non-Archimedian Functional Analysis," Marcel Dekker, 1978.

[56] Schaefer, H. H. and M. P. Wolff, "Topological Vector Spaces," Springer-Verlag, 1999.

[57] Schikhof, W. H., *Differentiation in non-archimedian valued fields*, Nederl. Akad. Wet., Proc., Ser. A **73** (1970), 35–44.

[58] ——, "Ultrametric Calculus – An Introduction to *p*-adic Analysis," Cambridge University Press, Cambridge, 1984.

[59] Serre, J.-P., "Lie Algebras and Lie Groups," Benjamin, New York, 1965.

[60] Souriau, J.-M., *Groupes différentiels de physique mathématique*, pp. 73–119 in: P. Dazord and N. Desolneux-Moulis, "Feuilletages et quantification géometrique," Journ. lyonnaises Soc. math. France 1983, Sémin. sud-rhodanien de géom. II, Hermann, Paris, 1984.

[61] Yoshioka, A., Maeda, Y., Omori, H., and Kobayashi, O., *On regular Fréchet-Lie groups VII; the group generated by pseudodifferential operators of negative order,* Tokyo J. Math. **7** (1984), 315–336.

W. Bertram, Institut Elie Cartan, Université Nancy I, Faculté des Sciences, B.P. 239, 54506 Vandœuvre-lès-Nancy, Cedex, France; bertram@iecn.u-nancy.fr

H. Glöckner, TU Darmstadt, FB Mathematik AG 5, Schlossgartenstr. 7, 64289 Darmstadt, Germany; gloeckner@mathematik.tu-darmstadt.de

K.-H. Neeb, TU Darmstadt, FB Mathematik AG 5, Schlossgartenstr. 7, 64289 Darmstadt, Germany; neeb@mathematik.tu-darmstadt.de